\documentclass[12pt]{article}
\usepackage[utf8]{inputenc}
\usepackage{mathtools}
\usepackage{booktabs}
\usepackage[english]{babel} 
\usepackage[protrusion=true,expansion=true]{microtype} 
\usepackage{amsmath,amsfonts,amsthm}
\usepackage{ amssymb }
\usepackage{graphicx}
\usepackage{array}
\usepackage{multirow}
\usepackage{hhline}
\usepackage[titletoc]{appendix}
\usepackage{extarrows}

\usepackage[margin=1in,hmarginratio=1:1,top=20mm,columnsep=20pt]{geometry} 
\usepackage{booktabs} 
\usepackage{natbib}
\usepackage{setspace}
\usepackage{smile}
\usepackage[table]{xcolor}

\usepackage{subfig}
\usepackage{tabularx}
\usepackage[font = small,labelfont=bf,textfont=it]{caption} 
\usepackage{footnote}
\usepackage{algorithm}

\newtheorem{assum}{Assumption}[subsection]

\newcommand{\argmin}{\operatorname*{argmin}}

\begin{document}

\title{Gaussian process regression: Optimality, robustness, and relationship with kernel ridge regression}

\author{ Wenjia Wang$^a$ and Bing-Yi Jing$^b$ \\
       $^a$The Hong Kong University of Science and Technology (Guangzhou)\\
      and\\
      The Hong Kong University of Science and Technology\\
        $^b$Department of Statistics and Data Science\\ 
        Southern University of Science and Technology}

\maketitle

\begin{abstract}%
Gaussian process regression is widely used in many fields, for example, machine learning, reinforcement learning and uncertainty quantification. One key component of Gaussian process regression is the unknown correlation function, which needs to be specified. In this paper, we investigate what would happen if the correlation function is misspecified. We derive upper and lower error bounds for Gaussian process regression with possibly misspecified correlation functions.  We find that when the sampling scheme is quasi-uniform, the optimal convergence rate can be attained even if the smoothness of the imposed correlation function exceeds that of the true correlation function. We also obtain convergence rates of kernel ridge regression with misspecified kernel function, where the underlying truth is a deterministic function. Our study reveals a close connection between the convergence rates of Gaussian process regression and kernel ridge regression, which is aligned with the relationship between sample paths of Gaussian process and the corresponding reproducing kernel Hilbert space. This work establishes a bridge between Bayesian learning based on Gaussian process and frequentist kernel methods with reproducing kernel Hilbert space.
\end{abstract}

\section{Introduction}
Gaussian process regression is widely applied in machine learning \citep{rasmussen2006gaussian}, including reinforcement learning \citep{rasmussen2003gaussian} and Bayesian optimization \citep{shahriari2016taking,frazier2018tutorial,bull2011convergence,klein2016fast}; spatial statistics \citep{cressie2015statistics,stein2012interpolation,matheron1963principles}; computer experiments \citep{santner2013design,sacks1989design}, and many others, to capture the intrinsic randomness of the underlying function. The goal of Gaussian process regression is to recover an underlying function based on noisy observations. As a Bayesian machine learning method, the key idea of Gaussian process regression is to impose a probabilistic structure, which is a Gaussian process, on the underlying truth. Based on this assumption, the conditional distribution at each unobserved point in the interested domain is normal with explicit mean and variance. The conditional mean provides a natural predictor of the function value, and the pointwise confidence interval constructed based on the conditional variance can be used for statistical uncertainty quantification. 

In this work, we establish error bounds on Gaussian process regression, where the smoothness of the correlation function can be misspecified, and the observations have noise. We consider that the underlying truth is a Gaussian process, which is a standard setting in Gaussian process modeling \citep{stein2012interpolation,santner2013design,gramacy2020surrogates}. Although the noisy observations have been extensively considered in the setting that the underlying truth is a deterministic function \citep[and references therein]{wynne2021convergence,steinwart2009optimal,fischer2017sobolev,vaart2011information},  (see Section \ref{sec:rlworks} for more discussions), it is somewhat surprising that there has been little study on this in the literature when the underlying truth is a Gaussian process. The difference is that, the convergence results for a deterministic function usually depend on some quantities related to the underlying function (e.g., the norm of the underlying function in some function space), while for a Gaussian process, these quantities themselves may be random. Thus, the convergence for a Gaussian process regression needs to be analyzed with a different approach. 
We derive prediction lower and upper error bounds under $L_2$ metric and with fixed design. Specifically, we show that if the smoothness of the true correlation function is $m_0$ and the smoothness of the imposed correlation function lies in $[m_0,\infty)$, with an appropriate regularization parameter and quasi-uniform design points, the convergence rate under $L_2$ metric is $n^{-(m_0-d/2)/(2m_0)}$, where $d$ is the dimension and $n$ is the sample size. Furthermore, we prove that this convergence rate is optimal under certain assumptions. Our theory can be applied to justifying the use of space-filling designs, where the design points spread approximately evenly in the region of interest, since quasi-uniform designs are space-filling designs. If the smoothness of the imposed correlation function, denoted by $m$, is less than $m_0$, we show that the convergence rate of upper error bound is $n^{-(m-d/2)/(2m)}$. 

Here, we should keep in mind not to confuse the setting of Gaussian process regression with the settings of other fields, including nonparametric regression \citep{gu2013smoothing,geer2000empirical} and posterior contraction of Gaussian process priors \citep{van2008rates, vaart2011information}. The hypothesis spaces are different in the later two fields. In particular, the underlying function is assumed to be \textit{deterministic}, which leads to different notions of \textit{smoothness} and convergence rates \citep{kanagawa2018gaussian,tuo2020kriging}. 

We also consider one popular kernel method: kernel ridge regression, where the reproducing kernel Hilbert space can be misspecified. This is a frequentist approach, where the underlying truth is assumed to be a deterministic function. The reason for considering kernel ridge regression is two-fold. 

First, the study paves the way to establish the intriguing relationship between Gaussian process regression and kernel ridge regression with more details given later. At first sight, the two areas are very different, for example, completely different approaches have been employed to investigate their convergence rates respectively. On the other hand, the two areas share some striking similarities in certain aspects, for example, their predictors take rather similar forms, and also their model assumptions bear strong resemblance. A thorough review on the differences and connection of Gaussian process and reproducing kernel Hilbert space can be found in \cite{kanagawa2018gaussian}. Therefore, it is natural to ask whether there are some deep relationships between Gaussian process regression and kernel ridge regression. \cite{kanagawa2018gaussian} provides a positive answer. Remark 5.5 of \cite{kanagawa2018gaussian} states a theoretical equivalence between Gaussian process regression and kernel ridge regression, where the Gaussian process regression model and the convergence rate \citep[Theorem 5.1]{kanagawa2018gaussian} is based on the posterior contraction of Gaussian process priors in \cite{vaart2011information}. Although the underlying truth in \cite{vaart2011information} is still a deterministic function, Remark 5.5 of \cite{kanagawa2018gaussian} reveals a relationship between Gaussian process regression and kernel ridge regression. Based on the constructed convergence rate in Gaussian process regression, we conduct a further investigation and establish a relationship based on the situations where ``the underlying truth in Gaussian process regression is a Gaussian process'' and ``the underlying truth in kernel ridge regression is a deterministic function''.

\begin{figure}[h!]
\centering
\includegraphics[width=0.9\textwidth]{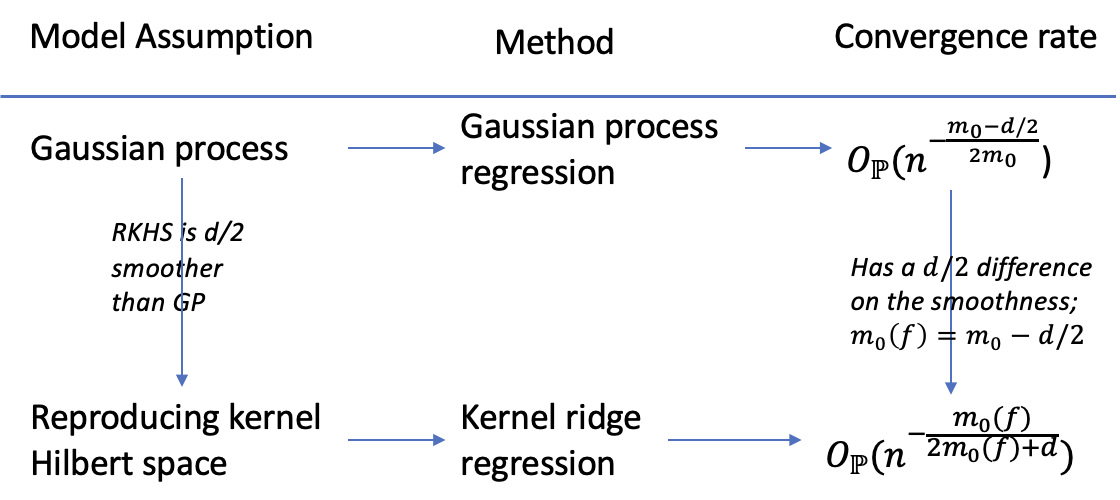}
\caption{\label{fig:fig1}Relationship between the convergence rates of oversmoothed Gaussian process regression ($m_0\leq m <\infty$) and kernel ridge regression. We use the following abbreviation. RKHS: Reproducing kernel Hilbert space; GP: Gaussian process.}
\end{figure}

We now describe briefly their relationship, which is summarized in Figure \ref{fig:fig1}. If the true correlation function has smoothness $m_0$, then the sample paths of the Gaussian process have a smoothness $m_0-d/2$, but do not lie in the Sobolev space $H^{m_0-d/2}$ with a strictly positive probability \citep{steinwart2019convergence,kanagawa2018gaussian}. For a deterministic function $f$ with smoothness $m_0(f)=m_0-d/2$, the optimal convergence rate is $n^{-m_0(f)/(2m_0(f)+d)}=n^{-(m_0-d/2)/(2m_0)}$, which coincides with the optimal convergence rate of Gaussian process regression. Furthermore, the optimal value of the regularization parameter in kernel ridge regression coincides with that of the regularization parameter in Gaussian process regression. In other words, we can regard Gaussian process regression as kernel ridge regression with an oversmoothed kernel function shifted by a smoothness $d/2$, from prediction perspective. This coincidence reveals an interesting relationship between kernel ridge regression and Gaussian process regression, and provides more insights on the connection between these two methods.

Second, we will derive some new and interesting results on convergence rates, which complements the existing literature on this topic. Specifically, suppose $f$ has smoothness $m_0(f)$, and the corresponding Sobolev space with smoothness $m_0(f)$ is denoted by $H^{m_0(f)}$. We show that the kernel ridge regression can still achieve the optimal convergence rate, if the corresponding reproducing kernel Hilbert space is associated with a smoothness lying in $[m_0(f)/2,\infty)$. If $f\in H^{m_0(f)}$, this recovers the convergence results in the misspecified kernel ridge regression literature (e.g., \cite{blanchard2018optimal,dicker2017kernel,guo2017learning,lin2017distributed,steinwart2009optimal,fischer2017sobolev}), while the settings are different. See Section \ref{subsec:crinDF} for detailed discussion. Nevertheless, we note that if a function has smoothness $m_0(f)$, it may not lie in the corresponding Sobolev space with smoothness $m_0(f)$; examples include triangle function and Mat\'ern kernel functions; see Section \ref{subsec:smoothness}. 
We are not aware of any work related to the convergence rate under the scenario $f\notin H^{m_0(f)}$ but has smoothness $m_0(f)$\footnote{
This implies $f\in H^{m}$ for all $m<m_0(f)$ but $f\notin H^{m_0(f)}$. Note that this is \textit{not} misspecification. The condition $f\notin H^{m_0(f)}$ is only related to $f$ itself, not the prediction method we use. In previous works, it is often assumed that $f\in H^{m}$ for all $m\leq m_0(f)$.}. Table \ref{tab:results} summarizes the results obtained in this work. For the briefness, we assume the design is quasi-uniform in Table \ref{tab:results}, and present general results in the main text.

\begin{table}
\centering
\begin{tabular}{|l|l|l|l|l|}
\hline
Method &  Model Assumption & \multicolumn{2}{|c|}{Convergence rate}\\\hline
 Gaussian process & $f$ is a realization of  & $O_{\mathbb
 {P}}(n^{-\frac{2m_0-d}{2m_0}})$ (optimal rate), & $O_{\mathbb
 {P}}(n^{-\frac{2m-d}{2m}})$,  \\
regression & a Gaussian process $Z$. & for $m_0 \leq m < \infty $. & for $m\leq m_0$. \\\hline

Kernel ridge & $f$ is a deterministic 
&  \cellcolor[gray]{0.9} $f\in H^{m_0(f)}$: $O_{\mathbb
 {P}}(n^{-\frac{m_0(f)}{2m_0(f)+d}})$, & \cellcolor[gray]{0.9} $O_{\mathbb
 {P}}(n^{-\frac{2m}{4m+d}})$, \\
regression  & function.   & $f\notin H^{m_0(f)}$: $O_{\mathbb
 {P}}(n^{-\frac{m_0(f)}{2m_0(f)+d}}Q(n))$, & \multirow{2}{*}{for $m < m_0(f)/2$.}  \\
& & for $m_0(f)/2 \leq m < \infty$. & \\\hline
\end{tabular}
\caption{\label{tab:results} Summary of the $L_2$ convergence rates of misspecified Gaussian process regression and kernel ridge regression. The function $Q$ satisfies $\lim_{s\rightarrow +\infty}(\log Q(s))/(\log s) = 0$. The two rates on the shaded row were presented in previous literature, while our settings and mathematical development are different. 
}
\end{table}

The rest of this paper is arranged as follows. 
We first make comparison to related works in Section \ref{sec:rlworks}.
In Section \ref{sec:pre}, we introduce some preliminaries. In Section \ref{subsec:GPmodelresults}, we provide convergence rates of misspecified Gaussian process regression. In Section \ref{sec:relation}, we discuss the relationship between the convergence rates of misspecified Gaussian process regression and kernel ridge regression, where we also present convergence rates of misspecified kernel ridge regression. Numerical experiments are conducted in Section \ref{sec:simulation}. Conclusions and discussion are made in Section \ref{sec:conclu}. Technical proofs are provided in Appendix.

\section{Related literature}\label{sec:rlworks}
In this subsection, we first remark some differences between our results and previous works. The previous works can be roughly divided into two fields: Gaussian process modeling, where the underlying truth is assumed to be a Gaussian process, and deterministic function reconstruction, where the underlying truth is modeled as a deterministic function. The difference between the convergence rate analysis in these two settings is that, in deterministic function reconstruction, the convergence rate usually involves some kind of function norm of the underlying true function, while for Gaussian process modeling, this norm itself is also random, which needs to be further considered. Although our work focuses on the Gaussian process modeling, we also consider the kernel ridge regression and obtain some interesting results. We utilize mathematical tools from both fields in the present work. For example, Lemma \ref{thm211} comes from Gaussian process modeling, and the rate of convergence of kernel ridge regression is established based on the previous works \cite{tuo2020improved, geer2000empirical}. Moreover, mathematical tools in scattered data approximation \citep{wendland2004scattered} play an important role in our analysis.

\subsection{Gaussian process modeling}\label{subsec:crinGP}
The rate of convergence of Gaussian process regression without noise has been studied in much literature, see \cite{seleznjev1999certain,yakowitz1985comparison,stein1990uniform} for example, where the convergence rate is pointwise or the input points are not general scattered data points. Recent works \cite{wang2019prediction,tuo2020kriging} study the rate of convergence of Gaussian process regression in the $L_p(\Omega)$ norm, with $1\leq p\leq \infty$ under different designs and misspeficied correlation functions in the noiseless case. To the best of our knowledge, the only work that studies Gaussian process regression with noisy observations is \cite{lederer2019uniform}. In \cite{lederer2019uniform}, a uniform error bound of Gaussian process regression has been provided, where the unknown realization $f$ and the correlation function are assumed to have a Lipschitz continuity, and the noise is normal. Furthermore, the correlation function in \cite{lederer2019uniform} is well-specified. 

In this work, we study the rate of convergence of Gaussian process regression in the $L_2(\Omega)$ norm, under different designs and misspeficied correlation functions, but we take the noise into consideration. These settings differentiate our work with the previous works in Gaussian process modeling.

\subsection{Deterministic function reconstruction}\label{subsec:crinDF}

Comparing with Gaussian process modeling, there are much more literature studying the deterministic function reconstruction. The most related fields to the present work are kernel ridge regression, posterior contraction of Gaussian process priors, and scattered data approximation.

\paragraph{Kernel ridge regression} Although we focus on the Gaussian process regression, we also consider the kernel ridge regression and obtain some interesting results. We consider that $f$ has smoothness $m_0(f)$, in the sense that is to be introduced later in Section \ref{subsec:smoothness}. If the underlying true function $f\in H^{m_0(f)}$, we recover the convergence rates obtained by \cite{blanchard2018optimal,dicker2017kernel,guo2017learning,lin2017distributed,steinwart2009optimal}, while the model settings are different. Specifically, the design points in the above works are random. Moreover, the assumptions are different. The aforementioned works impose conditions on the eigenvalues and eigenfunctions \citep{blanchard2018optimal,dicker2017kernel,guo2017learning,lin2017distributed,steinwart2009optimal,fischer2017sobolev}.

The aforementioned works have different model settings from our work,
which provides some additional insights on the study of kernel ridge regression. Specifically, we adopt model settings similar to \cite{tuo2020improved}, where the widely used Mat\'ern kernel functions can be used and the design points are fixed. These model settings allow us to consider the case that the underlying function $f\notin H^{m_0(f)}$ but has smoothness $m_0(f)$. We employ the empirical process technique together with Fourier transform to derive the convergence rates. Following this approach, we do not need to make assumptions on the eigenvalues and eigenfunctions, but we need additional conditions on the interested region. 
Moreover, our results show the advantage of space-filling designs. 

\paragraph{Posterior contraction of Gaussian process priors} 
In this field, despite the use of Gaussian process priors, the underlying function is still assumed to be deterministic. An incomplete list of papers in this area includes \cite{castillo2008lower,castillo2014bayesian,giordano2019consistency,nickl2017nonparametric,pati2015optimal,vaart2011information,van2008rates,van2016gaussian}. We are not aware of any error bounds in this area in terms of our settings, i.e., fixed designs, fill and separation distances.

\paragraph{Scattered data approximation} In the field of scattered data approximation, the goal is to approximate or interpolate an underlying deterministic function. Examples include \cite{wendland2004scattered,wendland2005approximate,rieger2009deterministic,narcowich2006sobolev}, which cover the noiseless case, and \cite{wynne2021convergence,rieger2009deterministic}, which cover the case that the observations have noise. The misspecification case is considered in \cite{narcowich2006sobolev,wynne2021convergence}. Although the observations are corrupted by noise in \cite{wynne2021convergence,rieger2009deterministic} as we considered in the present work, the convergence rates are different. If one plugs in our settings into their bounds, it can be seen that the prediction error bound does not converge to zero. This is the price for the more general noise assumption in \cite{wynne2021convergence,rieger2009deterministic}. We impose the sub-Gaussian assumption (see Condition (C5) in Section \ref{subsec:conditions}) and obtain a sharper error bound.

\section{Preliminaries}\label{sec:pre}
In this section, we introduce problem settings and conditions. 

\subsection{Problem settings}\label{subsec:GPintro}
Suppose that our observations $(x_k,y_k)$ satisfy the following model
\begin{align}
y_k = f(x_k) + \epsilon_k, k=1,...,n, \label{recoveringGP}
\end{align}
where $x_k \in \Omega\subset \RR^d$ and $\epsilon_k\overset{\text{i.i.d.}}{\sim}(0,\sigma_\epsilon^2)$, i.e., independent and identically distributed random noise with mean zero and variance $\sigma_\epsilon^2$. In Gaussian process regression, the underlying function $f$ is assumed to be a realization of a Gaussian process $Z$. From this point of view, we shall not differentiate $f$ and $Z$ in Gaussian process regression. We assume $Z$ is a zero-mean stationary Gaussian process, denoted by $Z\sim GP(0,\sigma^2 \Psi)$, with $\text{Cov}(Z(x),Z(x'))=\sigma^2 \Psi(x - x')$ for $x,x'\in \RR^d$. Here $\sigma^2$ is the variance, and $\Psi$ is the ture but typically unknown correlation function which is stationary, positive definite and integrable on $\RR^d$.

For a moment, assume $\epsilon_k$'s are normal. Given the correlation function $\Psi$ and conditional on $Y=(y_1,...,y_n)^T$, $Z(x)$ is normally distributed at an unobserved point $x$. The conditional expectation and variance of $Z(x)$ are given by 
\begin{align}
\mathbb{E}[Z(x)|Y]= &r(x)^T (R+\mu I_n)^{-1} Y,\label{mean}\\
\text{Var}[Z(x)|Y]= & \sigma^2(\Psi(x-x)-r(x)^T (R+\mu I_n)^{-1} r(x)),\label{var}
\end{align}
where $r(x)=(\Psi(x-x_1),\ldots,\Psi(x-x_n))^T, R=(\Psi(x_j-x_k))_{j k}$, $I_n$ is an identity matrix, and $\mu = \sigma_\epsilon^2/\sigma^2$. The conditional expectation \eqref{mean} is a natural predictor of $Z(x)$, and it can be shown that the conditional expectation is indeed the best linear predictor \citep{ankenman2010stochastic}, in the sense that it has the minimal mean squared prediction error (MSPE), which equals $\text{Var}[Z(x)|Y]$.

In this work, we investigate what happens if another correlation function $\Phi$, referred to as the \textit{imposed correlation function}, is used in Gaussian process regression in place of the \textit{true correlation function} $\Psi$. The resulting Gaussian process regression predictor after using $\Phi$ is
\begin{align}\label{eq:GPpredic}
\hat f_{G}(x) = r_{m}(x)^T (R_{m}+\mu_{m} I_n)^{-1} Y, \quad x\in \Omega,
\end{align}
where $r_{m}(x)=(\Phi(x-x_1),\ldots,\Phi(x-x_n))^T$ and $R_{m}=(\Phi(x_j-x_k))_{j k}$. We suppose $\mu_m$ is chosen according to our will and call it the \textit{regularization parameter}.  
Clearly, $\hat f_{G}$ in \eqref{eq:GPpredic} is no longer the best linear unbiased predictor. In this work, we are interested in the $L_2$ prediction error using the imposed correlation function $\Phi$, i.e.,
\begin{align}\label{eq:errorGP}
\|Z - \hat f_{G}\|_{L_2(\Omega)}.
\end{align}
Similar problem without the influence of noise has been considered in \cite{tuo2020kriging}. Other convergence results of Gaussian process regression with misspecified correlation functions can be found in \cite{stein1988asymptotically,stein1990bounds,stein1990uniform,tuo2020kriging,wang2019prediction,yakowitz1985comparison}, where the observations are noiseless. However, the appearance of noise can significantly change the analysis of convergence when the underlying truth is a Gaussian process\footnote{This is different with the settings in \cite{wynne2021convergence}, 
who also consider the noisy observations case with fixed designs. In \cite{wynne2021convergence}, the underlying truth is a deterministic function.}, as we will see later. 

\subsection{Notation and conditions}\label{subsec:conditions}
In the rest of this work, the following definitions are used. For two positive sequences $a_n$ and $b_n$, we write $a_n\asymp b_n$ if, for some $C,C'>0$, $C\leq a_n/b_n \leq C'$. Similarly, we write $a_n\gtrsim b_n$ if $a_n\geq Cb_n$ for some constant $C>0$, and $a_n\lesssim b_n$ if $a_n\leq C'b_n$ for some constant $C'>0$. Also, $C,C',c_j,C_j, j\geq 0$ are generic positive constants, of which value can change from line to line. We use $Q(s)$ to denote an increasing positive function satisfying 
\begin{align}\label{condition_h2}
	\lim_{s\rightarrow +\infty} \frac{\log Q(s)}{\log s} = 0
\end{align}
and not depending on $n$, which may vary at each occurrence. The Euclidean metric is denoted by $\|\cdot\|_2$. The Fourier transform of $f\in L_1(\mathbb{R}^d)$ is given by
$$\mathcal{F}(f)(\omega)=(2\pi)^{-d/2}\int_{\mathbb{R}^d} f(x) e^{-i x^T\omega}d x.$$

The following conditions will be assumed throughout the paper, unless otherwise specified.

\begin{enumerate}
    \item[(C1)] The region of interest $\Omega\subset \mathbb{R}^d$ is a compact set with positive Lebesgue measure and Lipschitz boundary, and satisfies an interior cone condition, i.e., there exist $\alpha\in (0,\pi/2)$ and $ \mathcal{R} > 0$ such that for every $x\in\Omega$, a unit vector $\xi(x)$ exists such that the cone
	$\mathcal{C}(x,\xi(x),\alpha,\mathcal{R}):=\left\{x+\lambda y:y\in\mathbb{R}^d,\|y\|=1,y^T\xi(x)\geq \cos\alpha,\lambda\in[0,\mathcal{R}]\right\} $
	is contained in $\Omega$.
    \item[(C2)]  There exists $m_0 > d/2$ such that,
	$$ c_1(1+\|\omega\|_2^2)^{-m_0} \leq  \mathcal{F}(\Psi)(\omega)\leq c_2(1+\|\omega\|_2^2)^{-m_0}, \forall \omega\in\mathbb{R}^d.$$
	\item[(C3)]  There exists $m > d/2$ such that,
	$$ c_3(1+\|\omega\|_2^2)^{-m} \leq  \mathcal{F}(\Phi)(\omega)\leq c_4(1+\|\omega\|_2^2)^{-m}, \forall \omega\in\mathbb{R}^d.$$
	\item[(C4)] Let $\mathcal{X} = \{X_1,X_2,...\}$ be a sequence of designs. Without loss of generality, assume that ${\rm card}(X_n)=n$, where $n$ takes its value in an infinite subset of $\mathbb{N}$, and card$(X)$ denote the cardinality of set $X$. We call $\mathcal{X}$ a \textit{sampling scheme}. The fill distance of $X_n$, defined by 
	\begin{align}\label{filldist}
        h_{X_n,\Omega}= \sup_{x\in\Omega}\inf_{x_j\in X_n}\|x-x_j\|_2,
    \end{align}
	satisfies
	\begin{align*}
	h_{X_n,\Omega}\leq Cn^{-1/d}, \forall n\geq 1.
	\end{align*} 
	\item[(C5)](Sub-Gaussian) Suppose $\epsilon_k$'s in \eqref{recoveringGP} are independent and identically distributed random variables satisfying
	\begin{align*}
	C^2 (\mathbb{E} e^{|\epsilon_k|^2/C^2}-1)\leq C', \quad k=1,...,n.
	\end{align*}
	Such random variables are called sub-Gaussian \citep{geer2000empirical}.
\end{enumerate}

Condition (C1) is a geometric condition on the region $\Omega$. We believe it holds in most practical situations, because 
the compactness and convexity imply the interior cone condition; see \cite{hofmann2007geometric,niculescu2006convex}.

Conditions (C2) and (C3) imply that the Fourier transforms of the true correlation function and imposed correlation function have an algebraical decay. A prominent class of correlation functions that have an algebraical decay of their Fourier transforms is the (isotropic) Mat\'ern correlation functions. The isotropic Mat\'ern correlation functions \citep{stein2012interpolation} is defined by
\begin{align}\label{materngai}
\Psi_M(x;\phi,\nu)=\frac{1}{\Gamma(\nu)2^{\nu -1}}(2\sqrt{\nu }\phi \| x\|_2)^{\nu } K_{\nu }(2\sqrt{\nu }\phi\| x\|_2),
\end{align}
with the Fourier transform \citep{tuo2015theoretical}
\begin{align}\label{maternspec}
\mathcal{F}(\Psi_M)(\omega;\nu ,\phi)= 4^{\nu+\frac{d}{2}}\pi^{\frac{d}{2}}\frac{\Gamma(\nu +\frac{d}{2})}{\Gamma(\nu )}(\nu  \phi^2)^{\nu } (4\nu \phi^2+\|\omega\|_2^2)^{-(\nu +\frac{d}{2})},
\end{align}
where $\phi>0$ is the scale parameter, and $K_{\nu}$ is the modified Bessel function of the second kind. The parameter $\nu$ is the smoothness parameter, which is associated with the smoothness of the kernel function $\Psi_M$.

Another example of correlation functions with algebraically decayed Fourier transforms is the generalized Wendland correlation function \citep{wendland2004scattered,gneitingstationary,chernih2014closed, bevilacqua2019estimation,fasshauer2015kernel}, defined as
\begin{align*}
    \Psi_{GW}(x) = \left\{
    \begin{array}{lc}
         \frac{1}{B(2\kappa,\eta+1)}\int_{\|\phi x\|}^1 u(u^2-\|\phi x\|^2)^{\kappa - 1}(1-u)^\eta du, & 0\leq \|x\|<\frac{1}{\phi}, \\
         0, & \|x\|\geq \frac{1}{\phi},
    \end{array}\right.
\end{align*}
where $\phi,\kappa > 0$ and $\eta \geq (d+1)/2 + \kappa$, and $B$ denotes the beta function. See Theorem 1 of \cite{bevilacqua2019estimation}.

In this work, we consider fixed designs, where the design points $x_1,...,x_n$ are fixed and can be chosen according to our will. Fixed designs are widely used in the field of computer experiments \citep{santner2013design}. Such designs include quasi-uniform designs \citep{borodachov2007asymptotics,utreras1988convergence}, maximin Latin hypercube designs \citep{van2007maximin}, optimal Latin hypercube designs \citep{park1994optimal}, and grid points. Condition (C4) states that the fill distance of designs can be controlled at a certain rate. It can be seen that any quasi-uniform sampling scheme satisfies Condition (C4), as stated in the following example.

\begin{definition}[Separation radius]
For $X=\{x_1,...,x_n\}$, define the separation radius as
\begin{align}\label{eq:sepdist}
	q_X=\min_{1\leq j\neq k\leq n}\|x_j-x_k\|_2/2.
\end{align}
\end{definition}

\begin{example}[Quasi-uniform designs]\label{eg:quasi}
	It is easy to check that $h_{X,\Omega}\geq q_X$ \citep{wendland2004scattered} for any set of points $X$. A sampling scheme $\mathcal{X}=\{X_1,X_2,\ldots\}$ is called quasi-uniform if $h_{X_n,\Omega}/q_{X_n}\leq C$ for all $n$. For a quasi-uniform sampling scheme, $h_{X_n,\Omega}\asymp q_{X_n}\asymp n^{-1/d}$ \citep{muller2009komplexitat}. 
	
	Obviously, a sampling scheme satisfying Condition (C4) may not be quasi-uniform. For example, we can add a point which is very close to one design point of a quasi-uniform design such that the separation radius is close to zero, and $h_{X_n,\Omega}/q_{X_n}\leq C$ does not hold.
\end{example}
\begin{remark}
	Random samplings do not satisfy Condition (C4); see Example 1 of \cite{tuo2020kriging}.
\end{remark}

\section{Rates of convergence for misspecified Gaussian process regression}\label{subsec:GPmodelresults}

In this section, we present our results on the convergence rate of the prediction error of misspecified Gaussian process regression \eqref{eq:errorGP}. 

We start with the easiest case. If the imposed correlation function is the same as the true correlation function, i.e., $\Phi = \Psi$ and $\mu_m = \sigma_\epsilon^2/\sigma^2$, then $\hat f_{G}$ is the best linear predictor and achieves the minimal MSPE, which is $\text{Var}[Z(x)|Y]$; see Section \ref{subsec:GPintro}. Obviously, the best linear predictor achieves the optimal convergence rate for a sampling scheme $\mathcal{X} = \{X_1,X_2,...\}$. It can be shown that if $\mu_m$ is any fixed positive constant and $\Phi = \Psi$, the optimal convergence rate can still be achieved. 
\begin{proposition}\label{prop:gpcons}
	Let $\hat f_{G}(x)$ be as in \eqref{eq:GPpredic} with $\Phi = \Psi$ and $\mu_m= C$, where $C>0$ is any fixed constant. For any fixed design $X=\{x_1,...,x_n\}\subset \Omega$, 
	\begin{align*}
	\mathbb{E}(\hat f_{G}(x) - Z(x))^2\leq C_1 \sigma^2(\Psi(x-x)-r(x)^T (R+\mu I_n)^{-1} r(x))=C_1\text{Var}[Z(x)|Y]
	\end{align*}
	holds for all $x\in \Omega$, where $R$, $r(x)$ and $\mu$ are as in \eqref{mean}, and the constant $C_1$ only depends on $C$, $\sigma_\epsilon^2$ and $\sigma^2$.
\end{proposition}
The proof of Proposition \ref{prop:gpcons} can be found in Appendix \ref{subsec:pfpropgpcons}. Because $\text{Var}[Z(x)|Y]$ is the minimal MSPE, $\mathbb{E}(\hat f_{G}(x) - Z(x))^2\geq \text{Var}[Z(x)|Y]$. Proposition \ref{prop:gpcons} shows that if the true correlation function is used, the regularization parameter can be changed to any fixed constant and would not influence the optimal convergence rate. However, Proposition \ref{prop:gpcons} does not provide any assertion on the convergence rate. 

In the following, we provide several error bounds of misspecified Gaussian process regression under noisy observations. Suppose that $\Psi$ and $\Phi$ satisfy Condition (C2) and Condition (C3), respectively. If $m_0 < m < \infty$, we call this case \textit{oversmoothed} case and call the corresponding imposed correlation function $\Phi$ \textit{oversmoothed correlation function}. On the other hand, if $d/2 < m < m_0$, we call this case \textit{undersmoothed} case and call the corresponding imposed correlation function $\Phi$ \textit{undersmoothed correlation function}. If $m=m_0$, we call this case \textit{well-specified} case.

We first provide an upper bound on the term $\Psi(x-x) - r(x)^T(R + \mu_1 I_n)^{-1}r(x)$ in the following proposition, which is closely related to the conditional variance in \eqref{var}. The proof of Proposition \ref{LEMMAERRORWNUG} is provided in Appendix \ref{app:pflemwnug}. Proposition \ref{LEMMAERRORWNUG} plays a key role in the proofs of Theorems \ref{thm:GPover} and \ref{thm:GPunder}.
\begin{proposition}\label{LEMMAERRORWNUG}
Suppose Conditions (C1)-(C4) hold. Then we have for any positive constant $\mu_1\gtrsim n^{1-2m_0/d}$, \begin{align}\label{eq:upnuget}
    \Psi(x-x) - r(x)^T(R + \mu_1 I_n)^{-1}r(x) \lesssim (\mu_1/n)^{1 - \frac{d}{2m_0}},
\end{align}
where $r(x)$ and $R$ are as in \eqref{mean}.
\end{proposition}

\begin{remark}
Proposition \ref{LEMMAERRORWNUG} is a deterministic version of Lemma F.8 in \cite{wang2020inference}. In Proposition \ref{LEMMAERRORWNUG}, the design points are fixed, while in Lemma F.8 of \cite{wang2020inference}, the design points are uniformly distributed on $\Omega$.
\end{remark}

\begin{remark}
Note that when the observations are noisy, the convergence rate of the conditional variance \eqref{var} can be directly obtained by setting $\mu_1=\mu$, where $\mu$ is as in \eqref{mean}. This result is different with the existing results in scattered data approximation \citep{wendland2004scattered,wu1993local}, where the observations have no noise.
\end{remark}

We start with the oversmoothed case. In the following theorem, we assume that both the true correlation function $\Psi$ and the imposed correlation function $\Phi$ are Mat\'ern correlation functions as in \eqref{materngai}. Recall that $h_{X,\Omega}$ and $q_X$ are the fill distance and separation radius for a design $X$ as defined in \eqref{filldist} and \eqref{eq:sepdist}, respectively. The proof of Theorem \ref{thm:GPover} is in Appendix \ref{subsec:pfgpover}.

\begin{theorem}[Oversmoothed Mat\'ern correlation function]\label{thm:GPover}
	Let $\Psi$ and $\Phi$ be two Mat\'ern correlation functions as in \eqref{materngai}. Suppose Conditions (C1)-(C5) hold. Suppose $m_0\leq m < \infty$ and $\mu_m \gtrsim n^{1-\frac{2m}{d}}$. Then, for any $t_1,t_2\geq C_0$ and all $n$, with probability at least $1-\exp(-t_1^2)-\exp(-t_2^2)$, we have
	\begin{align*}
	\|Z - \hat f_{G}\|_{L_2(\Omega)}^2\leq C\left((1+t_1)^2T + t_2^2\mu_m^{-\frac{d}{2m}}n^{-\left(1-\frac{d}{2m}\right)}\right),
	\end{align*}
	where $T = \mu_m^{-\frac{m-m_0}{m}} q_{X_n}^{-\frac{(m-m_0)d}{m}}(\mu_m/n)^{1-\frac{d}{2m}}+ \mu_m^{\frac{2m_0-d}{2m}} q_{X_n}^{\frac{(2m_0-d)d}{2m}}$, and the constants $C,C_0$ do not depend on $n$, $t_1$ and $t_2$.
\end{theorem}

The following corollary states that, if a sampling scheme is quasi-uniform, then Gaussian process regression with an oversmoothed Mat\'ern correlation function can still lead to the error bound $O_{\mathbb{P}}(n^{-\frac{2m_0-d}{2m_0}})$. Recall that a sampling scheme $\mathcal{X}=\{X_1,X_2,\ldots\}$ is said to be quasi-uniform if $h_{X_n,\Omega}/q_{X_n}\leq C$ for all $n$ (see Example \ref{eg:quasi}). Corollary \ref{corr:GPoverquasi} is a direct result of Theorem \ref{thm:GPover} and the proof is omitted.

\begin{corollary}[Oversmoothed Mat\'ern correlation function and quasi-uniform design]\label{corr:GPoverquasi}
	Let $\Psi$ and $\Phi$ be two Mat\'ern correlation functions as in \eqref{materngai}. Suppose Conditions (C1)-(C3) and (C5) hold. Suppose $m_0\leq m < \infty$ and the sampling scheme $\mathcal{X}$ is quasi-uniform. Let $\mu_m \asymp n^{-m/m_0+1}$. Then, for all $t\geq C_0$ and $n$, with probability at least $1-\exp(-t)$, we have
	\begin{align*}
	\|Z - \hat f_{G}\|_{L_2(\Omega)}^2\leq C (1+t)n^{-\frac{2m_0 - d}{2m_0}},
	\end{align*}
	where $C_0$ and $C$ are constants not depending on $n$ and $t$. In particular, we have
	\begin{align*}
	\|Z - \hat f_{G}\|_{L_2(\Omega)}^2=O_{\mathbb{P}}( n^{-\frac{2m_0 - d}{2m_0}}).
	\end{align*}
\end{corollary}

Now we consider Gaussian process regression with undersmoothed correlation functions. The following theorem indicates that the convergence rate is slower than that of Gaussian process regression with oversmoothed correlation functions, whose proof is provided in Appendix \ref{subsec:pfgpunder}. Note that in Theorem \ref{thm:GPunder}, the true and imposed correlation functions are not necessarily Mat\'ern correlation functions.

\begin{theorem}[Undersmoothed correlation function]\label{thm:GPunder}
    Suppose Conditions (C1)-(C5) hold. Suppose $d/2 < m \leq m_0$ and $\mu_m \gtrsim n^{1-\frac{2m}{d}}$. Then, for any $t_1,t_2\geq C_0$ and $n$, with probability at least $1-\exp(-t_1^2) - \exp(-t_2^2)$, we have
	\begin{align*}
	\|Z - \hat f_{G}\|_{L_2(\Omega)}^2\lesssim t_1^2\mu_m^{-\frac{d}{2m}}n^{-\left(1-\frac{d}{2m}\right)}+ (1+t_2)^2(\mu_m/n)^{1 - \frac{d}{2m}}.
	\end{align*}
	In particular, if $\mu_m$ is a fixed constant, $\|Z - \hat f_{G}\|_{L_2(\Omega)}^2=O_{\mathbb{P}}( n^{-\frac{2m-d}{2m}})$.
\end{theorem}

The following corollary provides error bounds in the well-specified case. Corollary \ref{corr:GPoverws} is a direct result of Theorem \ref{thm:GPunder}, and the proof is omitted.

\begin{corollary}[Well-specified correlation function]\label{corr:GPoverws}
    Suppose Conditions (C1)-(C5) hold. Furthermore, suppose $m = m_0$ and $\mu_m \asymp 1$. Then, for all $t\geq C_0$ and $n$, with probability at least $1-\exp(-t)$, we have
	\begin{align*}
	\|Z - \hat f_{G}\|_{L_2(\Omega)}^2\leq C(1+t)n^{-\frac{2m_0-d}{2m_0}},
	\end{align*}
	where $C_0$ and $C$ are constants not depending on $n$ and $t$. In particular, $\|Z - \hat f_{G}\|_{L_2(\Omega)}^2 = O_{\mathbb{P}}(n^{-\frac{2m_0-d}{2m_0}}).$
\end{corollary}

Theorem \ref{thm:GPunderlb} provides a lower error bound of Gaussian process regression, whose proof is presented in Appendix \ref{subsec:pfGPunderlb}. 

\begin{theorem}[Lower error bounds of Gaussian process regression]\label{thm:GPunderlb}
    Suppose Conditions (C1)-(C5) hold. Assume $m_0 > d$ and Assumption \ref{ass:efhasb} holds. Then we have
	\begin{align*}
	    \mathbb{E}\|Z - \hat f_{G}\|_{L_2(\Omega)}^2 \gtrsim n^{-\frac{2m_0-d}{2m_0}}.
	\end{align*}
\end{theorem}

\begin{remark}
Theorem \ref{thm:GPunderlb} requires a technical assumption Assumption \ref{ass:efhasb} in Appendix \ref{subsec:pfGPunderlb}, which essentially requires that there exists a correlation function $K$ with uniformly bounded eigenfunctions such that $\mathcal{F}(K)/\mathcal{F}(\Psi)$ is uniformly bounded. This assumption is slightly weaker than the assumption that $\Psi$ has uniformly bounded eigenfunctions. The later assumption is typical in nonparametric regression literature. See \cite{mendelson2010regularization,steinwart2009optimal} for example. Unfortunately, to the best of our knowledge, whether Assumption \ref{ass:efhasb} holds for Mat\'ern correlation functions is not present in literature.
\end{remark}

\begin{remark}
Note that the convergence rate in Theorem \ref{thm:GPunderlb} is different with the \textit{minimax} convergence rate in nonparameteric regression, where the underlying truth is a deterministic function. Besides the different settings, another difference is that the minimax convergence rate is considered in the \textit{worst} case for a given function class, while Theorem \ref{thm:GPunderlb} can be treated as in an \textit{average} case. We also note that \cite{tuo2020kriging} provide lower error bounds of Gaussian process regression in the noiseless case. 
\end{remark}

Combining Theorem \ref{thm:GPunderlb} and Corollary \ref{corr:GPoverquasi}, it can be seen that Gaussian process regression with an oversmoothed Mat\'ern correlation function achieves the optimal convergence rate, if 
the sampling scheme is quasi-uniform, and the optimal convergence rate for Gaussian process regression is $O_{\mathbb{P}}(n^{-\frac{2m_0-d}{2m_0}})$. Corollary \ref{corr:GPoverws} states that the optimal convergence rate can also be achieved if the Gaussian process regression with the true correlation function is used, which is intuitively true.

Theorem \ref{thm:GPunder} provides an upper bound on the $L_2$ prediction error of Gaussian process regression with an undersmoothed correlation function. The upper bound is larger than that of the Gaussian process regression with the true correlation function. Note that in \cite{tuo2020kriging} and \cite{wang2019prediction}, if the observations have no noise, it has been shown that using an oversmoothed Mat\'ern correlation function and a quasi-uniform sampling scheme can achieve the optimal convergence rate, while using an undersmoothed correlation function leads to an upper bound that has a slower convergence rate. Combining their results and ours, we can conclude that if the sampling scheme is quasi-uniform, using oversmoothed correlation functions is not detrimental to the convergence rate, no matter the observations are corrupted by noise or not. Nevertheless, we still recommend practitioners to try to find the correlation function with smoothness closed to the true smoothness. This is because the constant in the convergence rate can be large if the imposed smoothness is too far away from the true smoothness. Moreover, our results suggest that it is important to choose good designs in practice.

\section{Relationship of the convergence rates between Gaussian process regression and kernel ridge regression}\label{sec:relation}

In this section, we discuss the relationship between the convergence rates of Gaussian process regression and kernel ridge regression. For the conciseness of this paper, we move the introduction to reproducing kernel Hilbert spaces, Sobolev spaces, and kernel ridge regression to Appendix \ref{subsec:introrkhs}.

\subsection{Rates of convergence for misspecified kernel ridge regression}\label{sec:KRRresults}

\subsubsection{Smoothness of a deterministic function}\label{subsec:smoothness}

We say that a deterministic function $g\in L_2(\mathbb{R}^d)$ has a finite degree of smoothness if the quantity
\begin{eqnarray}\label{smoothness}
m_0(g):=\sup\{k\geq 0:g\in H^k(\mathbb{R}^d)\}
\end{eqnarray}
is finite, where $H^k(\mathbb{R}^d)$ is the Sobolev space with smoothness $k$. Here $k$ can be a non-integer, and the corresponding Sobolev space is called the \textit{fractional Sobolev space.} We call the quantity \eqref{smoothness} the smoothness of $g$. The functions considered in this work are assumed to have smoothness greater than $d/2$, which implies such function are continuous. Since $\Omega$ is compact and has a Lipschitz boundary, there exists an extension operator from $L_2(\Omega)$ to $L_2(\mathbb{R}^d)$, such that the smoothness of each function is maintained \citep{devore1993besov,rychkov1999restrictions}. We define the smoothness of a function $g\in L_2(\Omega)$ by the smoothness of the extended function $g_e\in L_2(\mathbb{R}^d)$ using (\ref{smoothness}).

From (\ref{smoothness}), it can be seen that a function $g$ with smoothness $m_0(g)$ can be divided into two scenarios: 1) $g\in H^{m_0(g)}(\Omega)$ but $g\not\in H^{m}(\Omega)$ for any $m>m_0(g)$; 2) $g\in H^{m}(\Omega)$ for any $m<m_0(g)$ but $g\not\in H^{m_0(g)}(\Omega)$. As a simple example, by \eqref{maternspec}, it can be checked that Mat\'ern correlation functions fall into the second scenario. To the best of our knowledge, the existing results on the kernel ridge regression only investigate the functions in the first scenario.

The following lemma provides a characterization of function $g$ that has smoothness $m_0(g)$ but is not in the Sobolev space $H^{m_0(g)}(\RR^d)$.

\begin{lemma}\label{LEMMA1}
	Let $m_0(g)\in (d/2,+\infty)$ be the smoothness of $g$. If $g\notin H^{m_0(g)}(\RR^d)$, there exists an increasing positive function $Q:\RR_+ \mapsto \RR_+$ satisfying \eqref{condition_h2} such that
	\begin{align*}
	& \int_{\mathbb{R}^d} \frac{|\mathcal{F}(g)(\omega)|^2}{Q(\|\omega\|)}(1+\|\omega\|^2)^{m_0(g)}d\omega \leq 1,\nonumber\\
	& \int_{\mathbb{R}^d} \frac{|\mathcal{F}(g)(\omega)|^2}{Q(\|\omega\|)}(1+\|\omega\|^2)^{m_0(g)+\delta}d\omega = \infty, \forall \delta>0.
	\end{align*}
\end{lemma}
Note that \eqref{condition_h2} implies $Q(s)$ increases slower than any $s^\delta$ with any $\delta>0$. The proof of Lemma \ref{LEMMA1} can be found in Appendix \ref{pfLEMMA1}. We use the following example to illustrate the intuition behind Lemma \ref{LEMMA1}.

\begin{example}\label{eg:tf1}
	Consider the triangle function
	\begin{align*}
	    f(x) = \left\{\begin{array}{ll}
	        1-|x|, &  |x|\leq 1,\\
	        0, & |x| > 1. 
	    \end{array}\right.
	\end{align*}
	It can be checked that $f$ has smoothness $3/2$ but $f\notin H^{3/2}(\RR)$. One can choose $Q(t):=C\log^2 (1+t)$ defined on $\RR_+$ with $C$ an appropriate constant such that
	\begin{align*}
	    \int_\RR \frac{|\mathcal{F}(f)(\omega)|^2}{Q(|\omega|)}(1 + |\omega|^2)^{3/2}d\omega & \leq 1,\\
	    \int_\RR \frac{|\mathcal{F}(f)(\omega)|^2}{Q(|\omega|)}(1 + |\omega|^2)^{3/2+\delta}d\omega & =\infty, \forall \delta>0.
	\end{align*}
	For the proof of the above statements, see Appendix \ref{app:pfeg1}.
\end{example}

\subsubsection{Main results for misspecified kernel ridge regression}
Let $f$ be a deterministic function with smoothness $m_0(f)$. The corresponding function space of interest is the Sobolev space $H^{m_0(f)}(\Omega)$, because by the definition of smoothness, $m_0(f) = \sup\{k>d/2:f\in H^{k}(\Omega)\}$. Theorem 10.45 of \cite{wendland2004scattered} suggests that if the kernel function $\Psi$ satisfies Condition (C2) with $m_0=m_0(f)$, $\mathcal{N}_\Psi(\Omega)$ coincides with the Sobolev space $H^{m_0(f)}(\Omega)$, where $\mathcal{N}_\Psi(\Omega)$ is the reproducing kernel Hilbert space generated by $\Psi$. Suppose a kernel ridge regression with reproducing kernel Hilbert space $\mathcal{N}_\Phi(\Omega)$ is used to recover the function $f$. Furthermore, assume $\Phi$ satisfies Condition (C3), which implies that the corresponding reproducing kernel Hilbert space $\mathcal{N}_\Phi(\Omega)$ coincides with the Sobolev space $H^{m}(\Omega)$. We call $\Phi$ the \textit{imposed kernel function}, and $\Psi$ the \textit{true kernel function}. 

\begin{remark}
	For any constant $c>0$, it can be seen that $\mathcal{N}_\Psi(\Omega)$ coincides with $\mathcal{N}_{c\Psi}(\Omega)$, and two norms are equivalent. Therefore, we pick any fixed kernel function $\Psi$ satisfying Condition (C2) with $m_0 = m_0(f)$ and call it the true kernel function. Any other kernel function is called imposed kernel function if it is used in the kernel ridge regression.  
\end{remark}

With a slight abuse of terminology, we refer to the kernel ridge regression with reproducing kernel Hilbert space $\mathcal{N}_\Phi(\Omega)$ as the misspecified kernel ridge regression. The misspecified kernel ridge regression can be written as
\begin{align}\label{KRRest1}
\hat f_{m} = \operatorname*{argmin}_{\hat f\in \mathcal{N}_{\Phi}(\Omega)}\bigg( \frac{1}{n}\sum_{k=1}^n(y_k-\hat f(x_k))^2 + \lambda_m \|\hat f\|^2_{\mathcal{N}_{\Phi}(\Omega)}\bigg),
\end{align}
where $\lambda_m>0$ is a regularization parameter. Note that if $\lambda_m=\mu_m/n$, where $\mu_m$ is as in \eqref{eq:GPpredic}, the representer theorem implies that 
$\hat f_{m}$ has the same form as $\hat f_G$ in \eqref{eq:GPpredic}. 

There are two cases, according to the smoothness $m$ of the reproducing kernel Hilbert space $\mathcal{N}_{\Phi}(\Omega)$, or equivalently, the smoothness of the Sobolev space $H^{m}(\Omega)$. If $m_0(f) \leq m < \infty$, the corresponding Sobolev space $H^{m}(\Omega)\subset H^{m_0(f)}(\Omega)$. With a slight abuse of terminology, we call this case \textit{oversmoothed} case and call the corresponding kernel function $\Phi$ \textit{oversmoothed kernel function}, even if $m$ may equal to $m_0(f)$. On the other hand, if $d/2 < m < m_0(f)$, we call this case \textit{undersmoothed} case and call $\Phi$ \textit{undersmoothed kernel function}.

In this work, we are interested in the convergence rate of the $L_2$ prediction error $\|f-\hat f_m\|_{L_2(\Omega)}$.
The following theorem states that, using an oversmoothed kernel function can still lead to the optimal convergence rate, if the regularization parameter is appropriately chosen. The proof of Theorem \ref{thm:krrover} is presented in Appendix \ref{subsec:pfkrrover}.

\begin{theorem}[Kernel ridge regression with oversmoothed kernel function]\label{thm:krrover}
	Suppose $f$ has smoothness $m_0(f)$. Suppose Conditions (C1)-(C5) hold and $m_0(f)\leq m < \infty$. 
	If $\lambda_m \asymp n^{-\frac{2m}{2m_0(f)+d}}$, the following statements are true for all $n$.
	\begin{enumerate}
		\item If $f\in H^{m_0(f)}(\Omega)$, then
		\begin{align*}
		\|f-\hat f_m\|_{L_2(\Omega)} = O_{\mathbb{P}}(n^{-\frac{m_0(f)}{2m_0(f)+d}}),
		\|\hat f_m\|_{\mathcal{N}_{\Phi}(\Omega)} =O_{\mathbb{P}}(n^{\frac{m-m_0(f)}{2m_0(f)+d}}).
		\end{align*}
		\item If $f\notin H^{m_0(f)}(\Omega)$, then
		\begin{align*}
		\|f-\hat f_m\|_{L_2(\Omega)} = O_{\mathbb{P}}(n^{-\frac{m_0(f)}{2m_0(f)+d}}Q(n)),
		\|\hat f_m\|_{\mathcal{N}_{\Phi}(\Omega)} =O_{\mathbb{P}}(n^{\frac{m-m_0(f)}{2m_0(f)+d}}Q(n)).
		\end{align*}
	\end{enumerate}
	
\end{theorem}

The next theorem states the convergence rate of upper error bounds in the undersmoothed case,
whose proof is provided in Appendix \ref{subsec:pfkrrunder}.

\begin{theorem}[Kernel ridge regression with undersmoothed kernel function]\label{thm:krrunder}
	Suppose $f$ has smoothness $m_0(f)$. Suppose Conditions (C1)-(C5) hold and $d/2 < m < m_0(f)$. If $\lambda_m \asymp n^{-\frac{2m}{2m_0(f)+d}}$, 
the following statements are true  for all $n$.
	\begin{enumerate}
		\item If $m_0(f)/2 \leq m < m_0(f)$, then we have:
		\begin{align*}
		\|f-\hat f_m\|_{L_2(\Omega)} =\left\{
		\begin{array}{ll}
		O_{\mathbb{P}}(n^{-\frac{m_0(f)}{2m_0(f)+d}}), & \mbox{ if } f\in H^{m_0(f)}(\Omega),\nonumber\\
		O_{\mathbb{P}}(n^{-\frac{m_0(f)}{2m_0(f)+d}}Q(n)), & \mbox{ if } f\notin H^{m_0(f)}(\Omega).
		\end{array}\right.
		\end{align*}
		\item If $d/2 < m < m_0(f)/2$, then we have:
		\begin{align*}
		\|f-\hat f_m\|_{L_2(\Omega)} =
		O_{\mathbb{P}}(n^{-\frac{2m}{4m+d}}).
		\end{align*}
	\end{enumerate}
\end{theorem}
From Theorems \ref{thm:krrover} and \ref{thm:krrunder}, we can see that the misspecified kernel ridge regression can still achieve the optimal convergence rate, as long as the imposed kernel function satisfies Condition (C3) with $m \geq m_0(f)/2$. These results generalize the results in \cite{tuo2020improved}, where the imposed kernel function satisfies $m = m_0(f)/2$. Furthermore, this work establishes the convergence results 
under the case that $f\notin H^{m_0(f)}(\Omega)$ but has smoothness $m_0(f)$.

\subsection{Relationship of the convergence rates of kernel ridge regression and Gaussian process regression}\label{subsec:relation}

Although kernel ridge regression and Gaussian process regression have different model assumptions, and we have applied completely different approaches to obtain the convergence rates of error bounds,
there is an intimate relationship between the constructed convergence rates. This relationship, notably, is aligned with the relationship between the reproducing kernel Hilbert space and Gaussian process, as we will explain in this section. For the ease of mathematical treatment, we assume that the sampling scheme is quasi-uniform. 
We use $\Psi_K$, $\Phi_K$, $\Psi_G$, and $\Phi_G$ to denote the true kernel function, the imposed kernel function, the true correlation function, and the imposed correlation function, respectively. 

We first link the prediction error of kernel ridge regression and that of Gaussian process regression, as shown in the following proposition. The proof is presented in Appendix \ref{subsec:pfofpropkg}.
\begin{proposition}\label{prop:kgrela1}
	Suppose $f\in \mathcal{N}_{\Psi_K}(\Omega)$ is a deterministic function, and $Z\sim GP(0,\Psi_G)$ is a Gaussian process. Suppose $\Psi_K$ and $\Phi_K$ are stationary, positive definite and integrable on $\RR^d$, $\Psi_G=\Psi_K$, $\Phi_G=\Phi_K$ and $\lambda_m = \mu_m/n$. Then
	\begin{align*}
	\mathbb{E}(f(x)-\hat f_m(x))^2 \leq C\mathbb{E}(Z(x)-\hat f_G(x))^2, \forall x\in \Omega,
	\end{align*}
	where $\hat f_m$ and $\hat f_G$ are as in \eqref{KRRest1} and \eqref{eq:GPpredic}, respectively, and $C=\max (1,\|f\|_{\mathcal{N}_{\Psi}}^2)$.
\end{proposition}
Proposition \ref{prop:kgrela1} states that the MSPE of kernel ridge regression $\mathbb{E}(f(x)-\hat f_m(x))^2$ on any point $x$ can be bounded by the MSPE of Gaussian process regression $\mathbb{E}(Z(x)-\hat f_G(x))^2$, when the correlation functions are the same as the kernel functions, and $\lambda_m = \mu_m/n$. However,
Proposition \ref{prop:kgrela1} does not provide the optimal convergence rate of the MSPE $\mathbb{E}(f(x)-\hat f_m(x))^2$. To see this, let $\Psi_K=\Psi_G=\Phi_K=\Phi_G$.
Furthermore, assume $\Psi_K$ satisfies Condition (C2).
The optimal convergence rate in kernel ridge regression
is achieved if $\lambda_{m} \asymp n^{-\frac{2m_0}{2m_0+d}}$. However, Proposition \ref{prop:gpcons} suggests that the optimal convergence rate of $\|Z-\hat f_G\|_{L_2(\Omega)}$ is achieved if $\mu_{m}$ is a fixed constant, and $\mu_{m}$ does not have the same order of magnitude as $n\lambda_{m}$. On the other hand, if we set $\lambda_{m} \asymp 1/n$, Theorem 4.1 of \cite{wang2020inference} implies that 
the optimal convergence rate in kernel ridge regression
cannot be achieved. In other words, if we use Gaussian process regression with correlation function $\Psi_G=\Psi_K$ and a constant regularization parameter to make prediction on a deterministic function in $\mathcal{N}_{\Psi_K}(\Omega)$, the optimal convergence rate cannot be achieved. The difference between the convergence rates of $\mathbb{E}(f(x)-\hat f_m(x))^2$ and $\mathbb{E}(Z(x)-\hat f_G(x))^2$ can be interpreted by the difference of the support of a Gaussian process and the corresponding reproducing kernel Hilbert space, where the former is typically larger than the later \citep{van2008reproducing}.

In our results of Gaussian process regression, the smoothness $m_0$ showing in Condition (C2) for a stationary Gaussian process $Z$ should be interpreted as the mean squared differentiability \citep{stein2012interpolation} of the Gaussian process, which is determined by the smoothness of the correlation function $\Psi_G$. This is different with the smoothness of deterministic functions. Nonetheless, we can consider the smoothness of sample paths of $Z$, under the usual definition of smoothness for deterministic functions, which reveals an interesting connection between convergence rates of kernel ridge regression and Gaussian process regression. If $Z\sim GP(0,\Psi_G)$ is a stationary Gaussian process with correlation function $\Psi_G$ satisfying Condition (C2), it can be shown that the sample path smoothness is lower than $m_0$ with probability one \citep{driscoll1973reproducing,kanagawa2018gaussian,steinwart2019convergence}. The difference between the support of a Gaussian process and the corresponding reproducing kernel Hilbert space has been sharply characterized by \cite{steinwart2019convergence}. Specifically, \cite{steinwart2019convergence} shows that the sample paths of Gaussian process $Z$ lie in Sobolev space $H^\alpha(\Omega)$ with $\alpha\in (d/2, m_0-d/2)$ with probability one, and do not lie in the Sobolev space $H^{m_0-d/2}(\Omega)$ with a strictly positive probability. This implies that the sample paths of Gaussian process $Z$ have smoothness $m_0-d/2$ with a strictly positive probability. Consider a deterministic function $f$ with smoothness $m_0(f)=m_0-d/2$ but not lying in $H^{m_0-d/2}(\Omega)$. Theorem \ref{thm:krrover} suggests that if an oversmoothed kernel function with smoothness $m$ is used and $\lambda_m \asymp n^{-\frac{2m}{2m_0(f)+d}} =  n^{-\frac{m}{m_0}}$, then the convergence rate is $n^{-\frac{m_0(f)}{2m_0(f)+d}}=n^{-\frac{m_0-d/2}{2m_0}}$, up to a difference of $Q(n)$ with $Q(n) = o(n^\delta)$ for any $\delta>0$. This convergence rate coincides with the optimal convergence rate of Gaussian process regression, and the choice of the regularization parameter has the same order of magnitude as $n\lambda_m$, i.e., $\mu_m \asymp n^{1-\frac{m}{m_0}}\asymp n\lambda_m$. If we choose the optimal order of magnitude of $\lambda_m = Cn^{-\frac{m}{m_0}}$ for any fixed positive constant $C$ and $\mu_m = n\lambda_m$, the predictors of Gaussian process regression and kernel ridge regression are identical \citep{kimeldorf1970correspondence}, and both achieve the optimal convergence rate. In other words, we can regard Gaussian process regression as kernel ridge regression with an oversmoothed kernel function, from the prediction perspective, and the optimal convergence rates are almost the same, up to a small order of $n^{\delta}$ with any $\delta>0$.

\begin{remark}
    \citet[Section 5.1]{kanagawa2018gaussian} also discuss relationship between Gaussian process regression and kernel ridge regression. The relationship of the convergence rate of kernel ridge regression and the rate of posterior contraction of Gaussian process priors is established. Note that the Gaussian process regression model and the convergence rate \citep[Theorem 5.1]{kanagawa2018gaussian} is based on the posterior contraction of Gaussian process priors in \cite{vaart2011information}, where the underlying truth is still a deterministic function. We consider ``the underlying truth in Gaussian process regression is a Gaussian process'' and ``the underlying truth in kernel ridge regression is a deterministic function''. This differentiates our discussion with that in \cite{kanagawa2018gaussian}.
\end{remark}

\section{Numerical experiments}\label{sec:simulation}
In this section, we conduct numerical experiments to study whether the convergence rates given by Theorems \ref{thm:GPover} and \ref{thm:GPunder} are accurate. We consider the region of interest $\Omega = [0,1]$. It has been shown in Theorems \ref{thm:GPover} and \ref{thm:GPunder} that, if $m_0\leq m$, taking $\mu\asymp n^{-m/m_0+1}$ leads to the error bound $O_{\mathbb{P}}(n^{-\frac{2m_0-d}{2m_0}})$; on the other hand, if $m_0> m$, taking $\mu\asymp 1$ yields the error bound $O_{\mathbb{P}}(n^{-\frac{2m-d}{2m}})$.

Let $\mathcal{E}=\mathbb{E}\|Z - \hat f_{G}\|_{L_2(\Omega)}^2$. We consider grid designs, such that the fill distance has the same order of magnitude of the separation distance. If the convergence rates of $\mathcal{E}$ are sharp, then we have the approximation
\begin{align}\label{simu1}
\log \mathcal{E}&\approx \frac{2m_0-d}{2m_0}\log (1/n) +\log c_1, \mbox{ if } m_0\leq m,\nonumber\\
\log \mathcal{E}&\approx \frac{2m-d}{2m}\log (1/n) +\log c_2, \mbox{ if } m_0 > m,
\end{align}
where $c_1$ and $c_2$ are constants. Therefore, in the numerical experiments, we regress $\log \mathcal{E}$ on $\log (1/n)$ and check whether the estimated slope is close to the theoretical assertion $\frac{2m_0-d}{2m_0}$ and $\frac{2m-d}{2m}$, when $m_0\leq m$ and $m_0>m$, respectively. 

We consider the sample sizes $n=10k$, for $k=2,3,...,15$. For each $k$, we simulate 100 realizations of a Gaussian process, where the correlation function is a Mat\'ern correlation function given by \eqref{materngai}. We take $\mu = 0.1\times n^{-m/m_0+1}$ when $m_0\leq m$, and take $\mu=0.1$ when $m_0>m$. The noise is set to be normal with mean zero and variance $0.25$. For $i$-th realization of a Gaussian process, we generate $10k$ grid points as $X$, and use $\mathcal{E}_i=\frac{1}{200}\sum_{j=1}^{200}(Z(x_j)-\hat f_G(x_j))^2$ to approximate $\|Z - \hat f_{G}\|_{L_2(\Omega)}^2$, where $x_j$'s are the first $200$ points of the Halton sequence \citep{niederreiter1992random}. This should provide a good approximation since the points are dense enough. The expectation $\mathcal{E}$ is approximated by $\frac{1}{100}\sum_{i=1}^{100}\mathcal{E}_i$.

The results are presented in Table \ref{Tab:simuResults}. The first two columns of Table \ref{Tab:simuResults} show the true and imposed smoothness. We consider three scenarios: oversmoothed case (row 1 and row 2), well-specified case (row 3), and undersmoothed case (row 4). The third and the fourth columns show the convergence rates obtained from the numerical experiments and the theoretical analysis, respectively. The fifth column shows the difference between the fourth and the fifth columns, and the last column gives the $R$-squared values of the linear regression of the simulated data.
\begin{table}[h]
\centering 
\begin{tabular}{|c|c|c|c|c|c|}
\hline
$m_0$ & $m$ &  Estimated slope & Theoretical slope & Difference & $R^2$\\
\hline
1.6 & 3.3 & 0.7138  & 0.6875 & 0.0263 &  0.9846\\
2.0 & 3.0 & 0.7664  & 0.7500 & 0.0164 & 0.9810\\
2.0 & 2.0 & 0.7691 & 0.7500  & 0.0191 & 0.9817\\
3.0 & 2.0 & 0.7856  & 0.7500  & 0.0356 &  0.9787\\
\hline
\end{tabular}
\caption{{\rm Numerical studies on the convergence rates of $\mathbb{E}\|Z - \hat f_{G}\|_{L_2(\Omega)}^2$.}}
\label{Tab:simuResults}
\end{table}

\begin{figure}[h!]
    \centering
    \begin{subfigure}
        \centering
        \includegraphics[height=2.3in]{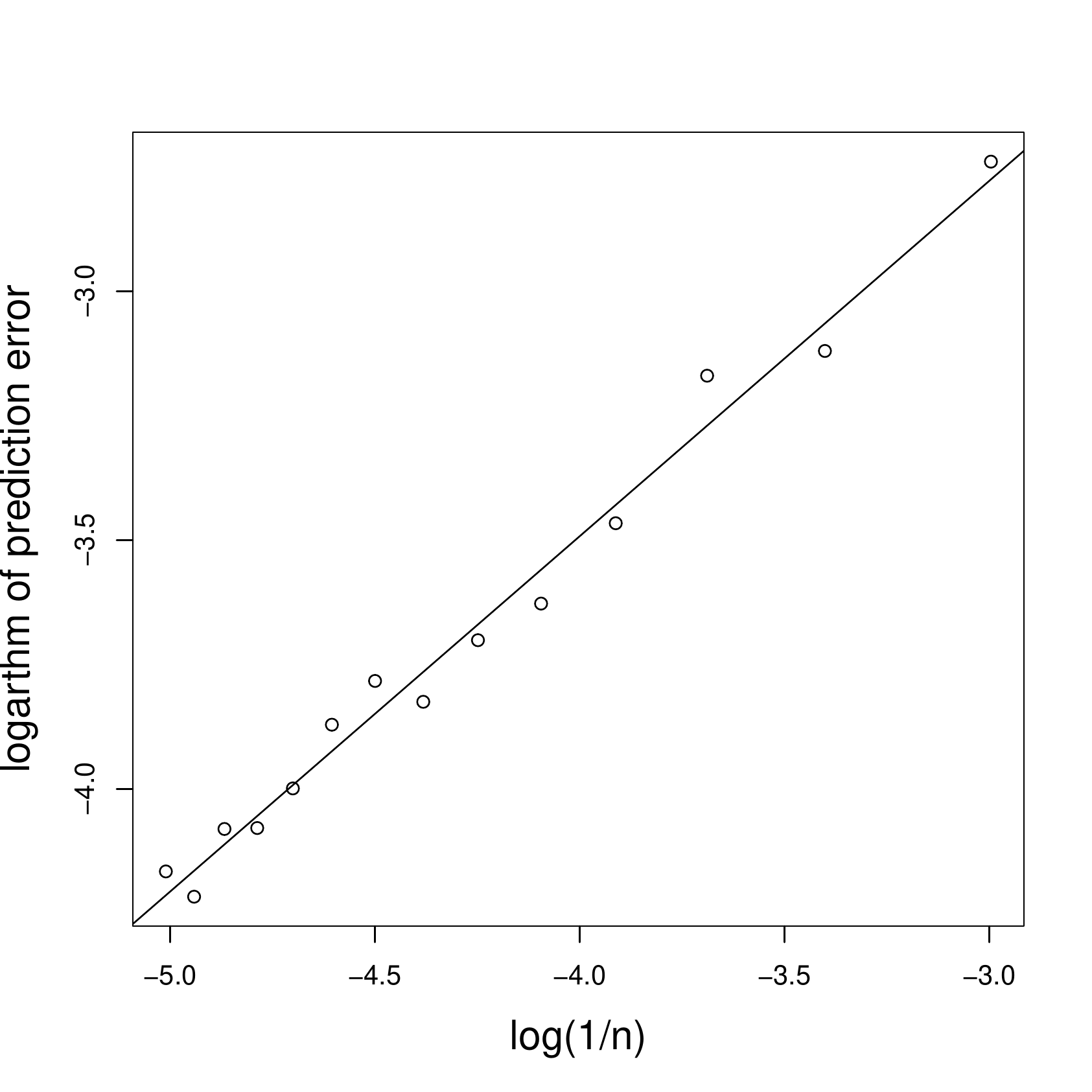}
    \end{subfigure}
    \begin{subfigure}
        \centering
        \includegraphics[height=2.3in]{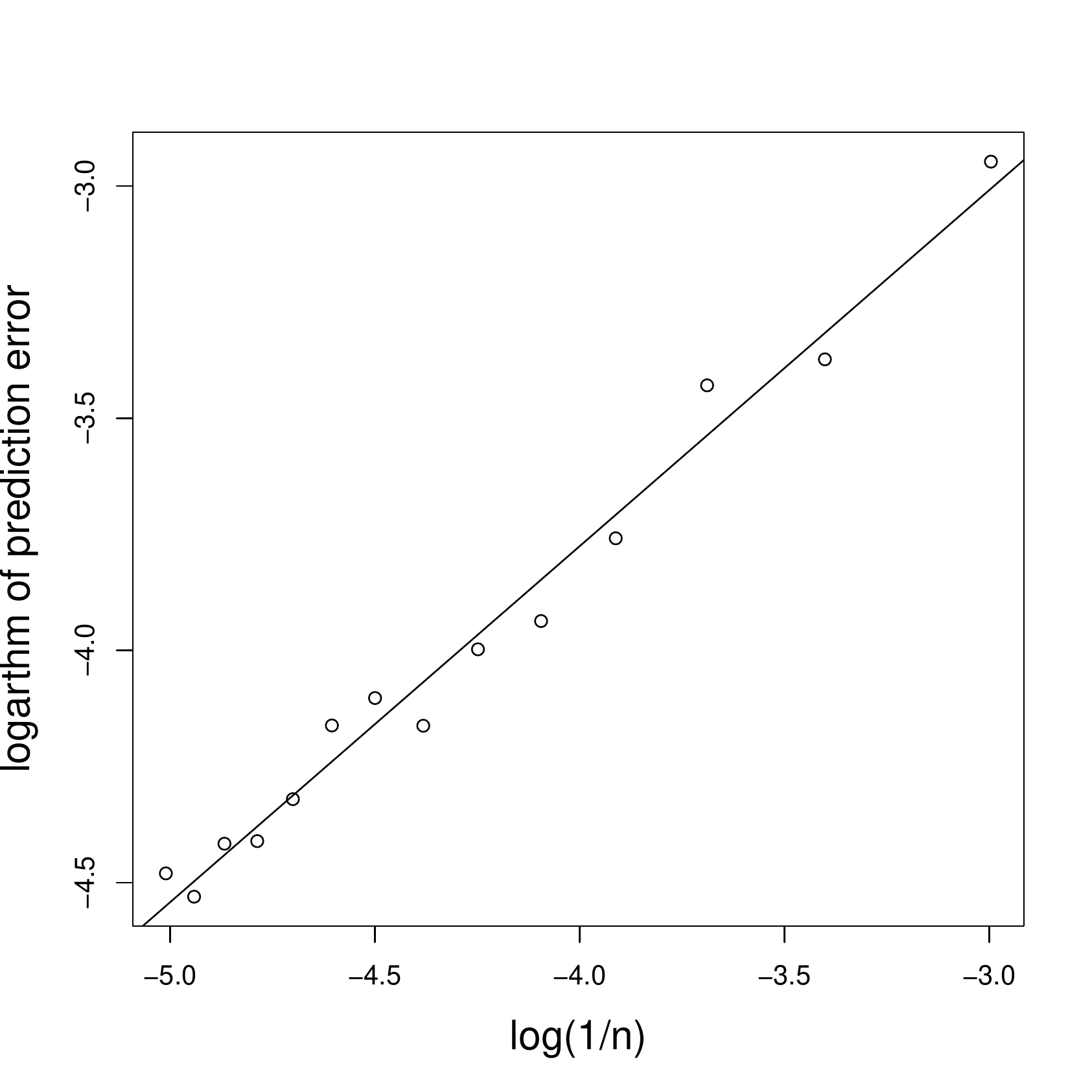}
    \end{subfigure}
    \begin{subfigure}
        \centering
        \includegraphics[height=2.3in]{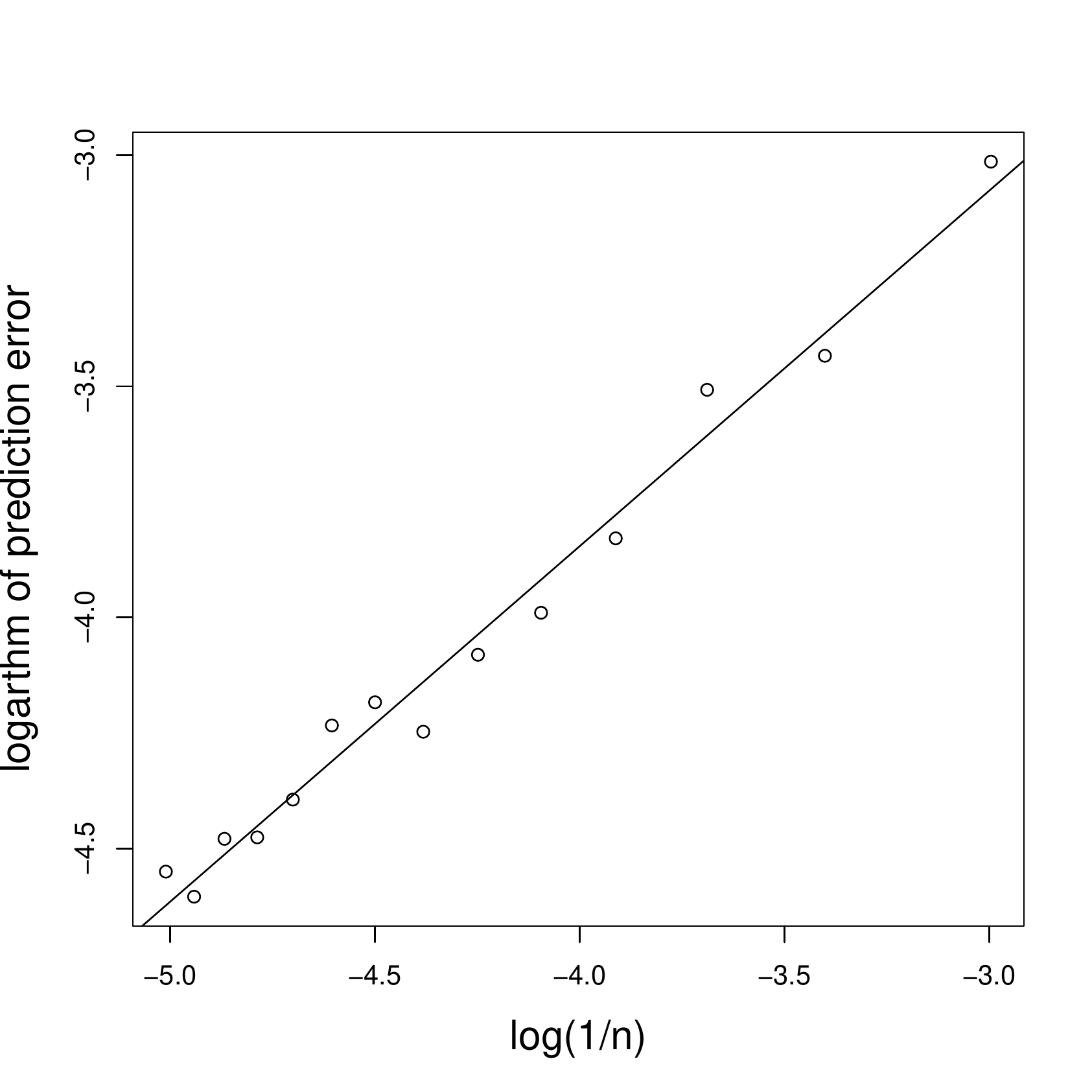}
    \end{subfigure}
    \begin{subfigure}
        \centering
        \includegraphics[height=2.3in]{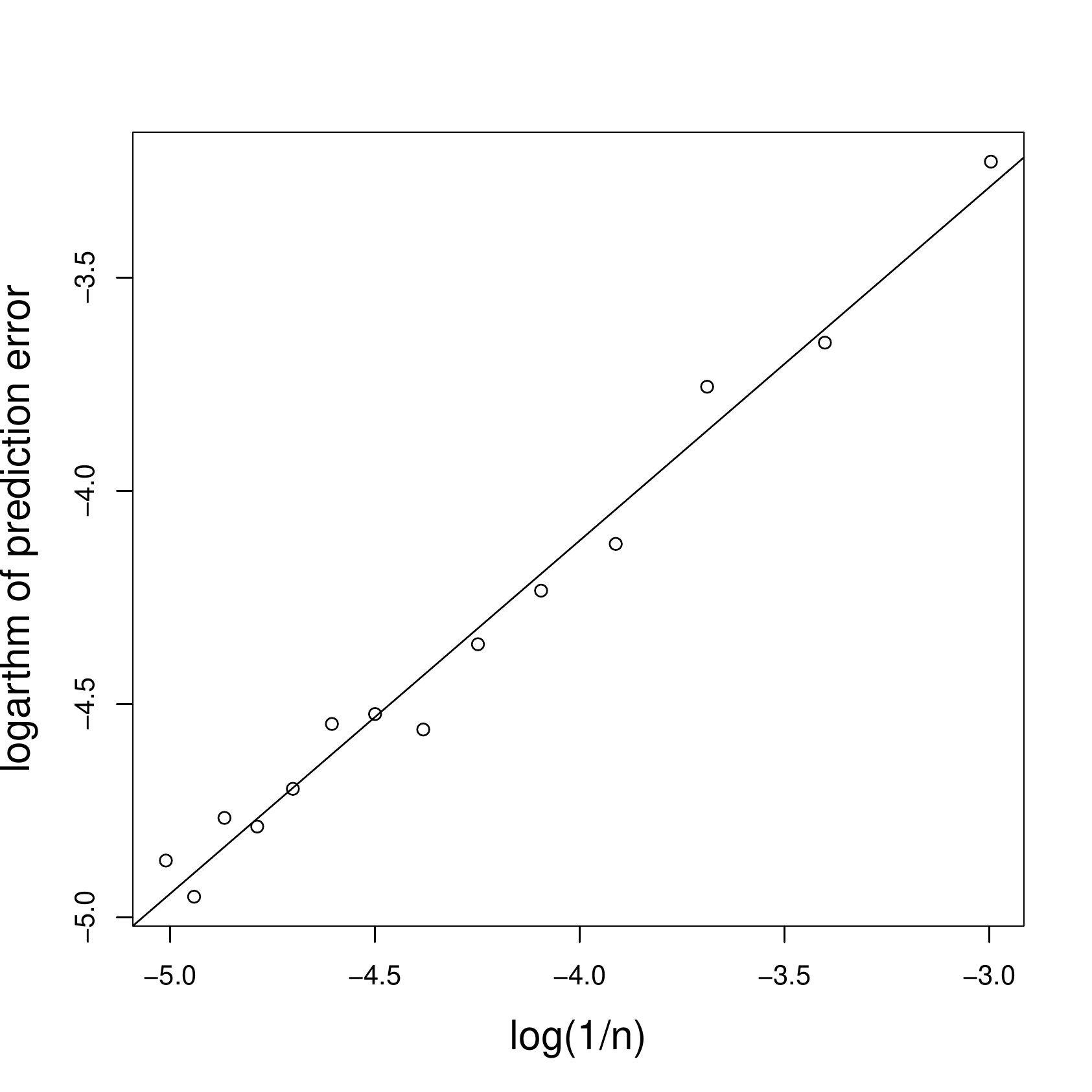}
    \end{subfigure}
   \caption{{\rm The regression line of $\log \mathcal{E}$ on $\log(1/n)$, under the four combinations of $(m_0,m)$ in Table \ref{Tab:simuResults}. Each point denotes one average prediction error for each $n$.}}
    \label{fig:nu0big}
\end{figure}

From Table \ref{Tab:simuResults}, it can be seen that the estimated slopes are close to our theoretical assertions for these cases. Figure \ref{fig:nu0big} shows the scattered points and the regression lines under the four combinations of $(m_0,m)$ in Table \ref{Tab:simuResults}.  From the $R$-squared values and Figure \ref{fig:nu0big}, we can see that the regression lines fit the scattered points well.

\section{Conclusions and discussion}\label{sec:conclu}
In this work, we provide some upper and lower error bounds for Gaussian process regression under misspecified correlation functions, when the observations are corrupted by noise. We show that the optimal convergence rate of Gaussian process regression can be achieved by using an oversmoothed Mat\'ern correlation function and a quasi-uniform sampling scheme. We also show that if the underlying truth is a deterministic function, the optimal convergence rate can still be achieved by kernel ridge regression if the kernel function is oversmoothed or not ``too undersmoothed''. Despite the difference of model assumptions and approaches in the proofs, we find an interesting connection between the constructed convergence rates of Gaussian process regression and kernel ridge regression. This connection is aligned with the connection between Gaussian process and reproducing kernel Hilbert space. The finding of the connection could serve as a bridge between Bayesian learning and frequentist learning, and may inspire new advances in these two seemingly separate fields.

There are several remaining problems. First, when the underlying truth is a Gaussian process, we consider fixed designs, which are also considered in \cite{tuo2020kriging,wang2019prediction,tuo2020improved}.
Whether the results hold for random sampling needs further study. Second, in addition to prediction, uncertainty quantification plays an important role in statistics. Since Gaussian process regression imposes a probabilistic structure on the underlying truth, it naturally induces an uncertainty quantification methodology via confidence interval. Uncertainty quantification under misspecification will be pursued in the future. 

\section*{Acknowledgements}
The authors are grateful to the AE and reviewers for their very constructive comments and suggestions. The authors also thank Rui Tuo at Texas A\&M for his helpful suggestions when writing this paper. Wang's work was supported by NSFC grant 12101149.

\bibliography{paperref}
\bibliographystyle{apalike}

\appendix
\section{Reproducing kernel Hilbert space, Sobolev space and kernel ridge regression}\label{subsec:introrkhs}
Reproducing kernel Hilbert space plays an important role in the study of kernel ridge regression and Gaussian process regression. Suppose $\Omega\subset \mathbb{R}^d$ satisfies Condition (C1). Assume that $K:\Omega \times \Omega \rightarrow \RR$ is a symmetric positive definite kernel function. Define the linear space
\begin{eqnarray}\label{FPhi}
F_{K}(\Omega)=\left\{\sum_{k=1}^n\beta_k K(\cdot,{x}_k):\beta_k\in \RR,{x}_k\in \Omega,n\in\mathbb{N}\right\},
\end{eqnarray}
and equip this space with the bilinear form
\begin{eqnarray*}
\left\langle\sum_{k=1}^n\beta_k K(\cdot,{x}_k),\sum_{j=1}^m\gamma_j K(\cdot, x'_j)\right\rangle_K:=\sum_{k=1}^n\sum_{j=1}^m\beta_k\gamma_j K({x}_k, x'_j).
\end{eqnarray*}
Then the \emph{reproducing kernel Hilbert space} $\mathcal{N}_{K}(\Omega)$ generated by the kernel function $K$ is defined as the closure of $F_{K}(\Omega)$ under the inner product $\langle\cdot,\cdot\rangle_{K}$, and the norm of  $\mathcal{N}_{K}(\Omega)$  is $\| f\|_{\mathcal{N}_{K}(\Omega)}=\sqrt{\langle f,f\rangle_{\mathcal{N}_{K}(\Omega)}}$, where $\langle\cdot,\cdot\rangle_{\mathcal{N}_{K}(\Omega)}$ is induced by $\langle \cdot,\cdot\rangle_{K}$. The following theorem gives another characterization of the reproducing kernel Hilbert space when $K$ is defined by a stationary kernel function $\Psi$, via the Fourier transform. Note that a kernel function $\Psi$ is said to be stationary if the value $\Psi(x,x')$ only depends on the difference $x-x'$. Thus, we can write $\Psi(x-x') :=\Psi(x,x')$. 

\begin{theorem}[Theorem 10.12 of \cite{wendland2004scattered}]\label{thm:NativeSpace}
	Let $\Psi$ be a positive definite kernel function which is stationary, continuous and integrable in $\RR^d$. Define 
	$$\mathcal{G}:=\{f\in L_2(\RR^d)\cap C(\RR^d):\mathcal{F}(f)/\sqrt{\mathcal{F}(\Psi)}\in L_2(\RR^d)\},$$
	with the inner product
	$$\langle f,g\rangle_{\mathcal{N}_\Psi(\RR^d)}=(2\pi)^{-d/2}\int_{\RR^d}\frac{\mathcal{F}(f)(\omega)\overline{\mathcal{F}(g)(\omega)}}{\mathcal{F}(\Psi)(\omega)}d \omega.$$ Then $\mathcal{G} = \mathcal{N}_\Psi(\RR^d)$, and both inner products coincide.
\end{theorem}

For an integer $k$, the Sobolev norm for function $g$ on $\mathbb{R}^d$ is defined by
\begin{align*}
\|g\|^2_{H^k(\mathbb{R}^d)}=\int_{\mathbb{R}^d} |\mathcal{F} (g)(\omega)|^2 (1+\|\omega\|_2^2)^{k} d \omega,
\end{align*}
and the inner product of a Sobolev space $H^k(\mathbb{R}^d)$ is defined by
$$\langle f,g\rangle_{H^k(\mathbb{R}^d)}=\int_{\RR^d}\mathcal{F}(f)(\omega)\overline{\mathcal{F}(g)(\omega)}(1+\|\omega\|_2^2)^{k}d \omega.$$
This definition can be naturally extended to Sobolev spaces with non-integer orders, which are commonly known as the \textit{fractional Sobolev spaces}, denoted by $H^m(\mathbb{R}^d)$ with a non-integer $m$. \begin{remark}\label{remarksob}
In this work, we are only interested in Sobolev spaces with $m > d/2$ because these spaces contain only continuous function according to the Sobolev embedding theorem. 
\end{remark}
A Sobolev space can also be defined on $\Omega \subset \RR^d$, denoted by $H^m(\Omega)$, with norm
\begin{eqnarray*}
	\|f\|_{H^m(\Omega)}=\inf\{\|f_e\|_{H^m(\RR^d)}:f_e\in H^m(\RR^d),f_e|_\Omega=f\},
\end{eqnarray*}
where $f_e|_\Omega$ denotes the restriction of $f_e$ to $\Omega$. It can be shown that $H^{m_0}(\mathbb{R}^d)$ coincides with the reproducing kernel Hilbert space $\mathcal{N}_\Psi(\RR^d)$ with equivalent norms, if $\Psi$ satisfies Condition (C2) (\cite{wendland2004scattered}, Corollary 10.13). By the extension theorem \citep{devore1993besov}, $\mathcal{N}_{\Psi}(\Omega)$ also coincides with $H^{m_0}(\Omega)$, and two norms are equivalent.

Given the observations $(x_k,y_k)$ with relationship \eqref{recoveringGP}, the kernel ridge regression reconstructs a function $f\in \mathcal{N}_{\Psi}(\Omega)$ by using
\begin{align}\label{KRRest}
\hat f = \operatorname*{argmin}_{g\in \mathcal{N}_{\Psi}(\Omega)}\bigg( \frac{1}{n}\sum_{k=1}^n(y_k-g(x_k))^2 + \lambda \|g\|^2_{\mathcal{N}_{\Psi}(\Omega)}\bigg),
\end{align}
where $\lambda$ is a prespecified regularization parameter. Under certain conditions and if $\Psi$ satisfies Condition (C2), the optimal order of magnitude of $\lambda$ is known in the literature \citep{geer2000empirical}, given by $\lambda = C n^{-\frac{2m_0}{2m_0+d}},$ where $C$ can be any fixed positive constant. The optimal choice of $\lambda$ leads to the optimal convergence rate under $L_2$ metric, which is $n^{-\frac{m_0}{2m_0+d}}$ \citep{stone1982optimal}.

\section{Notation}
We use $\langle \cdot, \cdot \rangle_n$ to denote the empirical inner product, which is defined by $$\langle f, g \rangle_n = \frac{1}{n}\sum_{k=1}^n f(x_k)g(x_k)$$ for two functions $f$ and $g$, and let $\|g\|_n^2 = \langle g, g \rangle_n$ be the empirical norm of function $g$. In particular, let $$\langle \epsilon, f \rangle_n = \frac{1}{n}\sum_{k=1}^n\epsilon_k f(x_k),$$ where $\epsilon = (\epsilon_1,\ldots,\epsilon_n)^T$.  
For two vectors $v$ and $w$, we use $\langle v, w \rangle=v^Tw$ to denote the inner product.

For notational simplification, let $h_n = h_{X_n,\Omega}$ and $q_n = q_{X_n}$ be the fill distance and separation radius of design $X_n$, respectively.
For the ease of treatment, in the rest of Appendix, we assume the regularization parameter $\lambda_m\asymp n^{\alpha}$ for some $\alpha\in \RR$. We use $tr(A)$ to denote the trace of a matrix $A$.

\section{Proof of Proposition \ref{prop:gpcons}}\label{subsec:pfpropgpcons}
Without loss of generality, assume $\sigma = 1$. Notice that for any $u = (u_1,...,u_n)^T \in \RR^n$ and any constant $C_1$,
\begin{align*}
\Psi(x-x) - 2\sum_{j=1}^n u_j\Psi(x - x_j) + \sum_{k=1}^n\sum_{j=1}^n u_k u_j\Psi(x_k - x_j) + C_1 \|u\|_2^2 \geq C_1 \|u\|_2^2,
\end{align*}
because $\Psi$ is positive definite. Plugging $u = (R + C_1 I_n)^{-1}r(x)$, we have 
\begin{align}\label{u2smallthanC}
C_1 r(x)^T(R + C_1 I_n)^{-2}r(x) \leq \Psi(x-x) - r(x)^T(R + C_1 I_n)^{-1}r(x).
\end{align}

If $C \leq \mu = \sigma_\epsilon^2/\sigma^2$, then direct computation shows that
\begin{align}\label{eq:prop1pf1}
& \mathbb{E}(\hat f_{G}(x) - Z(x))^2\nonumber\\
= & \Psi(x-x) - 2r(x)^T(R+ C I_n)^{-1}r(x) + r(x)^T(R+ C I_n)^{-1}(R+ \mu I_n)(R+ C I_n)^{-1}r(x)\nonumber\\
= & \Psi(x-x) - 2r(x)^T(R+ C I_n)^{-1}r(x) + r(x)^T(R+ C I_n)^{-1}(R+ C I_n)(R+ C I_n)^{-1}r(x)\nonumber\\
& + (\mu-C) r(x)^T(R+ C I_n)^{-2}r(x)\nonumber\\
= & \Psi(x-x) - r(x)^T(R+ C I_n)^{-1}r(x)  + (\mu-C) r(x)^T(R+ C I_n)^{-2}r(x)\nonumber\\
\leq & \bigg(1 + \frac{\mu-C}{C}\bigg)(\Psi(x-x) - r(x)^T(R+ C I_n)^{-1}r(x))\nonumber\\
\leq & \frac{\mu}{C}(\Psi(x-x) - r(x)^T(R+ \mu I_n)^{-1}r(x)),
\end{align}
where the first inequality is because of \eqref{u2smallthanC}, and the second inequality is because of $(R+ \mu I_n)^{-1} \preceq (R+ C I_n)^{-1}$.

If $C > \mu$, then we have 
\begin{align}\label{eq:prop1pf2}
& \mathbb{E}(\hat f_{G}(x) - Z(x))^2\nonumber\\
= & \Psi(x-x) - 2r(x)^T(R+ C I_n)^{-1}r(x) + r(x)^T(R+ C I_n)^{-1}(R+ \mu I_n)(R+ C I_n)^{-1}r(x)\nonumber\\
\leq & \Psi(x-x) - r(x)^T(R+ C I_n)^{-1}r(x),
\end{align}
where the first inequality is because of $R+ \mu I_n \preceq R+ C I_n$.

For any $u = (u_1,\ldots,u_n)^T$, the Fourier inversion theorem yields
\begin{align}\label{eq:identity1}
& \Psi(x-x) - 2\sum_{j=1}^n u_j\Psi(x - x_j) + \sum_{k=1}^n\sum_{j=1}^n u_k u_j\Psi(x_k - x_j) + C\|u\|_2^2\nonumber\\
= &  \frac{1}{(2\pi)^{d/2}}\int_{\mathbb{R}^d}\bigg|\sum_{j=1}^n u_j e^{i\langle x_j,\omega \rangle} - e^{i\langle x,\omega \rangle}\bigg|^2\mathcal{F}(\Psi)(\omega)d\omega + C\|u\|_2^2\nonumber\\
\leq & \frac{C}{\mu}\bigg(\frac{1}{(2\pi)^{d/2}}\int_{\mathbb{R}^d}\bigg|\sum_{j=1}^n u_j e^{i\langle x_j,\omega \rangle} - e^{i\langle x,\omega \rangle}\bigg|^2\mathcal{F}(\Psi)(\omega)d\omega + \mu\|u\|_2^2\bigg).
\end{align}
Let $u^{(1)}= (u_1^{(1)},...,u_n^{(1)})^T =  (R+ \mu I_n)^{-1}r(x)$. Because $u^{(2)}= (u_1^{(2)},...,u_n^{(2)})^T = (R+ C I_n)^{-1}r(x)$ is the solution to the optimization problem
\begin{align*}
\min_{u\in \RR^n} \Psi(x-x) - 2\sum_{j=1}^n u_j\Psi(x - x_j) + \sum_{k=1}^n\sum_{j=1}^n u_k u_j\Psi(x_k - x_j) + C\|u\|_2^2,
\end{align*}
we have
\begin{align}\label{eq:prop1pf3}
& \Psi(x-x) - r(x)^T(R+ C I_n)^{-1}r(x)\nonumber\\
= & \Psi(x-x) - 2\sum_{j=1}^n u_j^{(2)}\Psi(x - x_j) + \sum_{k=1}^n\sum_{j=1}^n u_k^{(2)} u_j^{(2)}\Psi(x_k - x_j) + C\|u^{(2)}\|_2^2\nonumber\\
\leq & \Psi(x-x) - 2\sum_{j=1}^n u_j^{(1)}\Psi(x - x_j) + \sum_{k=1}^n\sum_{j=1}^n u_k^{(1)} u_j^{(1)}\Psi(x_k - x_j) + C\|u^{(1)}\|_2^2\nonumber\\
\leq & \frac{C}{\mu}\bigg(\frac{1}{(2\pi)^{d/2}}\int_{\mathbb{R}^d}\bigg|\sum_{j=1}^n u_j^{(1)} e^{i\langle x_j,\omega \rangle} - e^{i\langle x,\omega \rangle}\bigg|^2\mathcal{F}(\Psi)(\omega)d\omega + \mu\|u^{(1)}\|_2^2\bigg)\nonumber\\
= & \frac{C}{\mu}(\Psi(x-x) - r(x)^T(R+ \mu I_n)^{-1}r(x)),
\end{align}
where the second inequality is by \eqref{eq:identity1}. Combining \eqref{eq:prop1pf1} and \eqref{eq:prop1pf3}, we finish the proof.

\section{Proof of Proposition \ref{LEMMAERRORWNUG}}\label{app:pflemwnug}
Let $I(x) = \Psi(x-x) - r(x)^T(R + \mu_1 I_n)^{-1}r(x)$. Consider function $g(t) = \Psi(x-t) - r(x)^T (R+\mu_1 I_n)^{-1} r(t)$. It can be seen that 
\begin{align}\label{eq:GpovergI}
I(x) = g(x) \leq \|g\|_{L_\infty(\Omega)}.
\end{align}
Direct computation shows that
\begin{align}\label{eq:lemIgNnorm}
    \|g\|_{\mathcal{N}_\Psi(\Omega)}^2 = & \Psi(x-x) - 2r(x)^T(R + \mu_1 I_n)^{-1}r(x) + r(x)^T (R+\mu_1 I_n)^{-1}R(R + \mu_1 I_n)^{-1}r(x)\nonumber\\
    \leq & \Psi(x-x) - r(x)^T(R + \mu_1 I_n)^{-1}r(x) = I(x).
\end{align}
By the Gagliardo–Nirenberg interpolation inequality for functions in Sobolev spaces \citep{leoni2017first,BrezisMironescu19}, it can be seen that
\begin{align}\label{eq:gpovergb}
\|g\|_{L_\infty(\Omega)} \lesssim \|g\|_{L_2(\Omega)}^{1-\frac{d}{2m}}\|g\|_{H^m(\Omega)}^{\frac{d}{2m}}\lesssim \|g\|_{L_2(\Omega)}^{1-\frac{d}{2m}}\|g\|_{\mathcal{N}_\Psi(\Omega)}^{\frac{d}{2m}} \lesssim \|g\|_{L_2(\Omega)}^{1-\frac{d}{2m}}I(x)^{\frac{d}{4m}}.
\end{align}
It remains to bound $\|g\|_{L_2(\Omega)}$. Let $f_1(t) = r(x)^T (R+\mu_1 I_n)^{-1} r(t)$. It can be seen from the representer theorem that
\begin{align}\label{eq:Gpoverf1}
f_1 = \argmin_{h \in \mathcal{N}_{\Psi}(\Omega)}\|\Psi(x-\cdot) - h\|_n^2+\frac{\mu_1}{n}\|h\|^2_{\mathcal{N}_{\Psi}(\Omega)}.
\end{align}
By Lemma \ref{lemmafixdesign},
\begin{align}\label{eq:Gpoverg1}
\|g\|_{L_2(\Omega)} \lesssim &  h_n^{m_0}\|g\|_{H^{m_0}(\Omega)} + \|g\|_n\lesssim h_n^{m_0}\|g\|_{\mathcal{N}_{\Psi}(\Omega)} + \|g\|_n
\lesssim h_n^{m}I(x)^{1/2} + \|g\|_n,
\end{align}
where the second inequality is by the equivalence of $\|\cdot\|_{H^{m_0}(\Omega)}$ and $\|\cdot\|_{\mathcal{N}_{\Psi}(\Omega)}$, and the last inequality is by \eqref{eq:lemIgNnorm}.

The empirical norm $\|g\|_n$ can be bounded by
\begin{align}\label{eq:Gpovergn2}
\|g\|_n^2 = & \|\Psi(x-\cdot) - f_1\|_n^2 \nonumber\\
= & \|\Psi(x-\cdot) - f_1\|_n^2 + \frac{\mu_1}{n}\|f_1\|^2_{\mathcal{N}_{\Psi}(\Omega)} - \frac{\mu_{1}}{n}\|f_1\|^2_{\mathcal{N}_{\Psi}(\Omega)}\nonumber\\
\leq & \|\Psi(x-\cdot) - \Psi(x-\cdot)\|_n^2 +  \frac{\mu_{1}}{n}\|\Psi(x-\cdot)\|^2_{\mathcal{N}_{\Psi}(\Omega)} - \frac{\mu_{1}}{n}r(x)^T (R+\mu_{1} I_n)^{-1} R (R+\mu_{1} I_n)^{-1}r(x) \nonumber\\
= & \frac{\mu_{1}}{n}\Psi(x-x) - \frac{\mu_{1}}{n}r(x)^T (R+\mu_{1} I_n)^{-1} (R+\mu_{1} I_n) (R+\mu_{1} I_n)^{-1}r(x)\nonumber\\
& + \frac{\mu_{1}^2}{n}r(x)^T (R+\mu_{1} I_n)^{-2} r(x)\nonumber\\
= & \frac{\mu_{1}}{n}\bigg(\Psi(x-x) - r(x)^T (R+\mu_1 I_n)^{-1} r(x) + \mu_1 r(x)^T (R+\mu_{1} I_n)^{-2} r(x)\bigg),
\end{align}
where the first inequality is because $f_1$ is the solution to \eqref{eq:Gpoverf1}.

Notice that for any $u = (u_1,...,u_n)^T \in \RR^n$,
\begin{align*}
\Psi(x-x) - 2\sum_{j=1}^n u_j\Psi(x-x_j) + \sum_{k=1}^n\sum_{j=1}^n u_k u_j\Psi(x_k-x_j) + \mu_1 \|u\|_2^2 \geq \mu_1 \|u\|_2^2,
\end{align*}
because $\Psi$ is positive definite. Plugging $u = (R + \mu_1 I_n)^{-1}r(x)$, we have 
\begin{align}\label{u2smallthan}
\mu_1 r(x)^T(R + \mu_1 I_n)^{-2}r(x) \leq \Psi(x-x) - r(x)^T(R + \mu_1 I_n)^{-1}r(x).
\end{align}
Therefore, \eqref{u2smallthan} and \eqref{eq:Gpovergn2} imply that
\begin{align}\label{eq:Gpovergn7}
\|g\|_n^2 
\leq & \frac{2\mu_{1}}{n}\bigg(\Psi(x-x) - r(x)^T (R+\mu_1 I_n)^{-1} r(x)\bigg) = \frac{2\mu_1}{n}I(x).
\end{align}
By \eqref{eq:GpovergI}, \eqref{eq:gpovergb}, \eqref{eq:Gpoverg1},  and \eqref{eq:Gpovergn7}, we have
\begin{align}\label{eq:GPIfinal}
I(x) \lesssim & (h_n^{m_0}I(x)^{1/2} + \|g\|_n)^{1-\frac{d}{2m_0}}I(x)^{\frac{d}{4m_0}}\nonumber\\
\lesssim & \bigg(h_n^{m_0}I(x)^{1/2} + \frac{\mu_1^{1/2}}{n^{1/2}}I(x)^{1/2}\bigg)^{1-\frac{d}{2m_0}}I(x)^{\frac{d}{4m_0}}\nonumber\\
\lesssim & \bigg(h_n^{m_0} + \frac{\mu_1^{1/2}}{n^{1/2}}\bigg)^{1-\frac{d}{2m_0}}I(x)^{1/2}\lesssim \bigg( \frac{\mu_1}{n}\bigg)^{1/2-\frac{d}{4m_0}}I(x)^{1/2},
\end{align}
where the last inequality is because $\mu_1\gtrsim n^{1-2m_0/d}$. It can be seen that \eqref{eq:GPIfinal} implies
\begin{align*}
I(x) \leq \bigg( \frac{\mu_1}{n}\bigg)^{1-\frac{d}{2m_0}}.
\end{align*}
This finishes the proof.

\section{Proof of Theorem \ref{thm:GPover}}\label{subsec:pfgpover}
We first present several lemmas used in this proof. Lemma \ref{thm211} is Lemma 24 in \cite{tuo2020kriging}. The proof of Lemma \ref{LEM:BOFKMMAX} is provided in Appendix \ref{app:pfbofKMmax}.

\begin{lemma}\label{thm211}
	Suppose $\Omega$ satisfies Condition (C1). Let $G$ be a zero-mean Gaussian process on $\Omega$ with continuous sample paths almost surely and with a finite maximum pointwise variance $\sigma^2_G = \sup_{x\in \Omega} \mathbb{E}G(x)^2<\infty$. Then for all $u>0$ and $1\leq p < \infty$, we have
	\begin{align*}
	& \mathbb{P}\left(\|G\|_{L_p(\Omega)} - \mathbb{E} \|G\|_{L_p(\Omega)} > u\right) \leq e^{-u^2/(2C_p\sigma^2_G)},\\
	& \mathbb{P}\left(\|G\|_{L_p(\Omega)} -\mathbb{E} \|G\|_{L_p(\Omega)}  <  - u\right) \leq e^{-u^2/(2C_p\sigma^2_G)},
	\end{align*}
	with $C_p={\rm Vol}(\Omega)^{2/p}$. Here ${\rm Vol}(\Omega)$ denotes the volume of $\Omega$.
\end{lemma}	
\begin{lemma}\label{LEM:BOFKMMAX}
	Suppose the design points $X=\{x_1,...,x_n\}$ and the separation radius of $X$ $q_X\lesssim 1$. Let $\Psi$ be a Mat\'ern correlation function satisfying Condition (C2) and $\Lambda_X$ be the maximum eigenvalue of matrix $(\Psi(x_j-x_k))_{jk}$. Then 
	\begin{align*}
	\Lambda_X\leq Cq_X^{-d},
	\end{align*}
	where $C$ is a constant depending on $\Psi$ and $\Omega$.
\end{lemma}

Now we begin to prove Theorem \ref{thm:GPover}. Recall that $y_j = Z(x_j) + \epsilon_j$. Let $\epsilon = (\epsilon_1,...,\epsilon_n)^T$ and $F = (Z(x_1),...,Z(x_n))^T$. Therefore, 
\begin{align*}
\hat f_{G}(x) = r_{m}(x)^T (R_{m}+\mu_{m} I_n)^{-1} F + r_{m}(x)^T (R_{m}+\mu_{m} I_n)^{-1} \epsilon.
\end{align*}
Direct computation shows that 
\begin{align}\label{eq:gpoverde}
(Z(x) - \hat f_{G}(x))^2 = & (Z(x) - r_{m}(x)^T (R_{m}+\mu_{m} I_n)^{-1} F - r_{m}(x)^T (R_{m}+\mu_{m} I_n)^{-1}\epsilon)^2\nonumber\\
\leq & 2(Z(x) - r_{m}(x)^T (R_{m}+\mu_{m} I_n)^{-1} F)^2 + 2(r_{m}(x)^T (R_{m}+\mu_{m} I_n)^{-1}\epsilon)^2\nonumber\\
= & 2I_1(x) +2I_2(x),
\end{align}
where the inequality is by the Cauchy-Schwarz inequality. Let $G(x)=Z(x) - r_{m}(x)^T (R_{m}+\mu_{m} I_n)^{-1} F$. It can be seen that $G$ is also a mean zero Gaussian process. 

By Jensen's inequality and Fubini's theorem,
\begin{align}\label{eq:gpoverEG}
& (\mathbb{E} \|G\|_{L_2(\Omega)})^2\leq \mathbb{E} \|G\|_{L_2(\Omega)}^2 = \int_\Omega \mathbb{E} (Z(x) - r_{m}(x)^T (R_{m}+\mu_{m} I_n)^{-1} F)^2 dx\nonumber\\
=  & \int_\Omega \mathbb{E}I_1(x) dx\leq {\rm Vol}(\Omega)\sup_{x\in \Omega} \mathbb{E}I_1(x),
\end{align}
where Vol$(\Omega)$ is the volumn of $\Omega$.

Next, we provide a uniform upper bound on $\mathbb{E}I_1(x)$. Direct computation gives us
\begin{align*}
\mathbb{E}I_1 = \Psi(x-x) - 2 r_{m}(x)^T (R_{m}+\mu_{m} I_n)^{-1} r(x)+ r_{m}(x)^T (R_{m}+\mu_{m} I_n)^{-1} R(R_{m}+\mu_{m} I_n)^{-1}r_{m}(x),
\end{align*}
where $r(x)$ and $R$ are as in \eqref{mean}.

By the Fourier inversion theorem, for $u=(u_1,...,u_n)^T = (R_{m}+\mu_{m} I_n)^{-1}r_{m}(x)$, we have 
\begin{align}\label{eq:eI1gpover}
& \mathbb{E}I_1(x) =  \int_{\RR^d}\left|\sum_{j=1}^nu_je^{-i\langle x_j, \omega \rangle}-e^{-i\langle x, \omega \rangle}\right|^2\mathcal{F}(\Psi)(\omega)d\omega\nonumber\\
= & \int_{\|\omega\|_2\leq \gamma }\left|\sum_{j=1}^nu_je^{-i\langle x_j, \omega \rangle}-e^{-i\langle x, \omega \rangle}\right|^2\mathcal{F}(\Psi)(\omega)d\omega + \int_{\|\omega\|_2> \gamma}\left|\sum_{j=1}^nu_je^{-i\langle x_j, \omega \rangle}-e^{-i\langle x, \omega \rangle}\right|^2\mathcal{F}(\Psi)(\omega)d\omega\nonumber\\
= & I_{11} + I_{12},
\end{align}
where $\gamma>1$ will be determined later. 

The first term $I_{11}$ can be bounded by
\begin{align}\label{eq:I11GPover}
I_{11} \lesssim &  \int_{\|\omega\|_2\leq \gamma }\left|\sum_{j=1}^nu_je^{-i\langle x_j, \omega \rangle}-e^{-i\langle x, \omega \rangle}\right|^2(1+\|\omega\|_2^2)^{-m_0}d\omega\nonumber\\
\lesssim  & \gamma^{2m-2m_0}\int_{\|\omega\|_2\leq \gamma }\left|\sum_{j=1}^nu_je^{-i\langle x_j, \omega \rangle}-e^{-i\langle x, \omega \rangle}\right|^2(1+\|\omega\|_2^2)^{-m}d\omega\nonumber\\
\lesssim & \gamma^{2m-2m_0}\int_{\|\omega\|_2\leq \gamma }\left|\sum_{j=1}^nu_je^{-i\langle x_j, \omega \rangle}-e^{-i\langle x, \omega \rangle}\right|^2\mathcal{F}(\Phi)(\omega)d\omega,
\end{align}
where the first and third inequalities are because of Conditions (C2) and (C3), respectively.

By the Cauchy-Schwarz inequality, the second term $I_{12}$ can be further split to
\begin{align}\label{eq:I12GPover}
I_{12} \leq &  2 \int_{\|\omega\|_2> \gamma}\left|\sum_{j=1}^nu_je^{-i\langle x_j, \omega \rangle}\right|^2\mathcal{F}(\Psi)(\omega)d\omega + 2\int_{\|\omega\|_2> \gamma}\left|e^{-i\langle x, \omega \rangle}\right|^2\mathcal{F}(\Psi)(\omega)d\omega\nonumber\\
= & I_{3} + I_{4}.
\end{align}
Since $\Psi$ satisfies Condition (C2), the term $I_{4}$ can be bounded by
\begin{align}\label{eq:I4GPover}
I_{4} \lesssim \int_{\|\omega\|_2> \gamma}(1+\|\omega\|_2^2)^{-m_0}d\omega\lesssim \gamma^{-2m_0+d},
\end{align}
and the term $I_{3}$ can be bounded by
\begin{align}\label{eq:I3GPover}
I_{3} & \lesssim \int_{\|\omega\|_2> \gamma}\left|\sum_{j=1}^nu_je^{-i\langle x_j, \omega \rangle}\right|^2(1+\|\omega\|_2^2)^{-m_0}d\omega\nonumber\\
& \lesssim \gamma^d\int_{\|\omega\|_2> 1}\left|\sum_{j=1}^nu_je^{-i\gamma \langle x_j, \omega \rangle}\right|^2(1+\gamma^2\|\omega\|_2^2)^{-m_0}d\omega\nonumber\\
& \lesssim  \gamma^{d-2m_0}\int_{\|\omega\|_2> 1}\left|\sum_{j=1}^nu_je^{-i\gamma \langle x_j, \omega \rangle}\right|^2(1+\|\omega\|_2^2)^{-m_0}d\omega\nonumber\\
& \lesssim  \gamma^{d-2m_0}\int_{\RR^d}\left|\sum_{j=1}^nu_je^{-i\gamma \langle x_j, \omega \rangle}\right|^2(1+\|\omega\|_2^2)^{-m_0}d\omega,
\end{align}
where the second inequality is by the change of variables, and the third inequality is by the fact that $(1+\gamma^2\|\omega\|_2^2)^{-m_0} \leq (\gamma^2(1+\|\omega\|_2^2)/2)^{-m_0}$ for $\|\omega\|_2\geq 1$.

Putting \eqref{eq:I4GPover} and \eqref{eq:I3GPover} into \eqref{eq:I12GPover}, we obtain
\begin{align}\label{eq:I12GPover2}
I_{12}& \lesssim \gamma^{d-2m_0}\int_{\RR^d}\left|\sum_{j=1}^nu_je^{-i\gamma \langle x_j, \omega \rangle}\right|^2(1+\|\omega\|_2^2)^{-m_0}d\omega + \gamma^{-2m_0+d}\nonumber\\ 
& \lesssim \gamma^{d-2m_0}\int_{\RR^d}\left|\sum_{j=1}^nu_je^{-i\gamma \langle x_j, \omega \rangle}\right|^2\mathcal{F}(\Psi)(\omega)d\omega + \gamma^{-2m_0+d}\nonumber\\ 
& = \gamma^{d-2m_0}\sum_{j,k=1}^nu_ju_k\Psi(\gamma x_j - \gamma x_k)+  \gamma^{-2m_0+d}\nonumber\\
& \leq \gamma^{d-2m_0}\Lambda_{\gamma X}\|u\|_2^2+  \gamma^{-2m_0+d},
\end{align}
where the second inequality is because $\Psi$ satisfies Condition (C2), and $\gamma X = \{\gamma x_1,...,\gamma x_n\}$. The separation distance of $\gamma X$ is $\gamma q_X$, which, together with Lemma \ref{LEM:BOFKMMAX}, implies  $\Lambda_{\gamma X}\leq C(\gamma q_X)^{-d}$ (it can be seen by the choice of $\gamma$ later, $\gamma q_X\lesssim 1$). Plugging \eqref{eq:I11GPover} and \eqref{eq:I12GPover2} into \eqref{eq:eI1gpover}, we have
\begin{align}\label{eq:eI1gpover2}
& \mathbb{E}I_1(x)\nonumber\\
\lesssim & \gamma^{2m-2m_0}\int_{\|\omega\|_2\leq \gamma }\left|\sum_{j=1}^nu_je^{-i\langle x_j, \omega \rangle}-e^{-i\langle x, \omega \rangle}\right|^2\mathcal{F}(\Phi)(\omega)d\omega + \gamma^{d-2m_0}\Lambda_{\gamma X}\|u\|_2^2+ \gamma^{-2m_0+d}\nonumber\\
\lesssim & \gamma^{2m-2m_0}\int_{\|\omega\|_2\leq \gamma }\left|\sum_{j=1}^nu_je^{-i\langle x_j, \omega \rangle}-e^{-i\langle x, \omega \rangle}\right|^2\mathcal{F}(\Phi)(\omega)d\omega + \gamma^{d-2m_0}(\gamma q_X)^{-d}\|u\|_2^2+ \gamma^{-2m_0+d}\nonumber\\
\lesssim & \gamma^{2m-2m_0}\bigg(\int_{\|\omega\|_2\leq \gamma }\left|\sum_{j=1}^nu_je^{-i\langle x_j, \omega \rangle}-e^{-i\langle x, \omega \rangle}\right|^2\mathcal{F}(\Phi)(\omega)d\omega + \gamma^{-2m}(q_X)^{-d}\|u\|_2^2\bigg)+ \gamma^{-2m_0+d}.
\end{align}
Take $\gamma = \mu_m^{-\frac{1}{2m}} (q_X)^{-d/2m}$ such that $\gamma^{-2m}(q_X)^{-d} = \mu_m$. Clearly, $\gamma q_X\lesssim 1$. By \eqref{eq:eI1gpover2}, $\mathbb{E}I_1(x)$ can be further bounded by
\begin{align}\label{eq:eI1gpover3}
\mathbb{E}I_1(x)
\lesssim  & \mu_m^{-\frac{m-m_0}{m}} q_X^{-\frac{(m-m_0)d}{m}}\bigg(\int_{\|\omega\|_2\leq \gamma }\left|\sum_{j=1}^nu_je^{-i\langle x_j, \omega \rangle}-e^{-i\langle x, \omega \rangle}\right|^2\mathcal{F}(\Phi)(\omega)d\omega + \mu_m\|u\|_2^2\bigg)\nonumber\\
&+ \mu_m^{\frac{2m_0-d}{2m}} (q_X)^{\frac{(2m_0-d)d}{2m}}\nonumber\\
\lesssim &  \mu_m^{-\frac{m-m_0}{m}} q_X^{-\frac{(m-m_0)d}{m}}\bigg(\int_{\RR^d}\left|\sum_{j=1}^nu_je^{-i\langle x_j, \omega \rangle}-e^{-i\langle x, \omega \rangle}\right|^2\mathcal{F}(\Phi)(\omega)d\omega + \mu_m\|u\|_2^2\bigg)\nonumber\\
&+ \mu_m^{\frac{2m_0-d}{2m}} (q_X)^{\frac{(2m_0-d)d}{2m}}\nonumber\\
= &  \mu_m^{-\frac{m-m_0}{m}} q_X^{-\frac{(m-m_0)d}{m}}\bigg(\Phi(x-x) -2\sum_{j=1}^n u_j\Phi(x-x_j) + \sum_{j,k=1}^n u_ju_k\Phi(x_j-x_k) + \mu_m\|u\|_2^2\bigg)
\nonumber\\
& + \mu_m^{\frac{2m_0-d}{2m}} (q_X)^{\frac{(2m_0-d)d}{2m}}\nonumber\\
=  & \mu_m^{-\frac{m-m_0}{m}} q_X^{-\frac{(m-m_0)d}{m}}(\Phi(x-x) - r_m(x)^T(R_m + \mu_mI_n)^{-1}r_m(x))+ \mu_m^{\frac{2m_0-d}{2m}} (q_X)^{\frac{(2m_0-d)d}{2m}}\nonumber\\
\lesssim & \mu_m^{-\frac{m-m_0}{m}} q_X^{-\frac{(m-m_0)d}{m}}(\mu_m/n)^{1-\frac{d}{2m}}+ \mu_m^{\frac{2m_0-d}{2m}} (q_X)^{\frac{(2m_0-d)d}{2m}},
\end{align}
where the last inequality is by applying Proposition \ref{LEMMAERRORWNUG} to the correlation function $\Phi$.

Let $T = \mu_m^{-\frac{m-m_0}{m}} q_X^{-\frac{(m-m_0)d}{m}}(\mu_m/n)^{1-\frac{d}{2m}}+ \mu_m^{\frac{2m_0-d}{2m}} (q_X)^{\frac{(2m_0-d)d}{2m}}$. Lemma \ref{thm211}, \eqref{eq:gpoverEG}, and \eqref{eq:eI1gpover3} imply that for all $t_1>0$, with probability at least $1-\exp(-Ct_1^2)$,
\begin{align*}
\|G\|_{L_2(\Omega)}\lesssim (1+t_1)T^{1/2}.
\end{align*}
which implies
\begin{align}\label{eq:gpoverI1r}
\|G\|_{L_2(\Omega)}^2\lesssim (1+t_1)^2T.
\end{align}
Next, we consider $I_2(x)$. This can be done by applying Theorem \ref{thm:krroverg}. Recall that $I_2(x)=(r_{m}(x)^T (R_{m}+\mu_{m} I_n)^{-1}\epsilon)^2$. 
Define $f_2(x) = 0$ for all $x\in \Omega$, i.e., $f_2$ is a zero function. Clearly, $f_2 \in \mathcal{N}_\Phi(\Omega)$ with $\|f_2\|_{\mathcal{N}_\Phi(\Omega)} = 0$. Let $\hat f_2(x) = r_{m}(x)^T (R_{m}+\mu_{m} I_n)^{-1}\epsilon$, then $\|\hat f_2\|_{L_2(\Omega)}^2 = \int_{x\in \Omega} I_2(x) dx$. By the representer theorem, $\hat f_2$ is the solution to 
\begin{align*}
    \min_{\hat g \in \mathcal{N}_\Phi(\Omega)}\frac{1}{n}\sum_{j=1}^n\left(\hat g(x_j)-\epsilon_j\right)^2 +\frac{\mu_m}{n}\|\hat g\|_{\mathcal{N}_\Phi(\Omega)}.
\end{align*}
Theorem \ref{thm:krroverg} tells us that for all $t>C_0$ (with appropriate changes of notation), with probability at least $1-C_1\exp(-C_2t^2)$, 
\begin{align*}
    \|f_2-\hat f_2\|_{L_2(\Omega)}^2 \leq C_3t^2 n^{-1}(\mu_m/n)^{-\frac{d}{2m}},
\end{align*}
which implies
\begin{align}\label{eq:gpoverI2r}
    \int_{x\in \Omega} I_2(x) dx = \|\hat f_2\|_{L_2(\Omega)}^2 \leq C_3t^2 n^{-1}(\mu_m/n)^{-\frac{d}{2m}}
\end{align}
since $f_2=0$. Note that \eqref{eq:gpoverde} implies $\|Z-\hat f_G\|_{L_2(\Omega)}^2 \lesssim \|G\|_{L_2(\Omega)}^2 + \|\hat f_2\|_{L_2(\Omega)}^2$. Thus, by \eqref{eq:gpoverI1r} and \eqref{eq:gpoverI2r}, we finish the proof of Theorem \ref{thm:GPover}.

\section{Proof of Theorem \ref{thm:GPunder}}\label{subsec:pfgpunder}

Note that $y_j = Z(x_j) + \epsilon_j$. Let $\epsilon = (\epsilon_1,...,\epsilon_n)^T$ and $F = (Z(x_1),...,Z(x_n))^T$. Therefore, 
\begin{align*}
\hat f_{G}(x) = r_{m}(x)^T (R_{m}+\mu_{m} I_n)^{-1} F + r_{m}(x)^T (R_{m}+\mu_{m} I_n)^{-1} \epsilon.
\end{align*}
Similar to \eqref{eq:gpoverde}, we have
\begin{align*}
(Z(x) - \hat f_{G}(x))^2 \leq & 2(Z(x) - r_{m}(x)^T (R_{m}+\mu_{m} I_n)^{-1} F)^2 + 2(r_{m}(x)^T (R_{m}+\mu_{m} I_n)^{-1}\epsilon)^2\nonumber\\
= & 2I_1(x) +2I_2(x).
\end{align*}
Let $G(x)=Z(x) - r_{m}(x)^T (R_{m}+\mu_{m} I_n)^{-1} F$, which is also a mean zero Gaussian process. 

Similar to \eqref{eq:gpoverEG}, we can obtain that 
\begin{align}\label{eq:gpunderEG}
(\mathbb{E} \|G\|_{L_2(\Omega)})^2\leq {\rm Vol}(\Omega)\sup_{x\in \Omega} \mathbb{E}I_1(x).
\end{align}
Direct computation gives us
\begin{align*}
\mathbb{E}I_1(x) = \Psi(x-x) - 2 r_{m}(x)^T (R_{m}+\mu_{m} I_n)^{-1} r(x) + r_{m}(x)^T (R_{m}+\mu_{m} I_n)^{-1} R(R_{m}+\mu_{m} I_n)^{-1}r_{m}(x).
\end{align*}
By the Fourier inversion theorem and Conditions (C2) and (C3), for $u=(u_1,...,u_n)^T = (R_{m}+\mu_{m} I_n)^{-1}r_{m}(x)$, we have 
\begin{align}\label{eq:eI1gpunder}
\mathbb{E}I_1(x) = & \int_{\RR^d}\left|\sum_{j=1}^nu_je^{-i\langle x_j, \omega \rangle}-e^{-i\langle x, \omega \rangle}\right|^2\mathcal{F}(\Psi)(\omega)d\omega\nonumber\\
\leq  & C_1\int_{\RR^d}\left|\sum_{j=1}^nu_je^{-i\langle x_j, \omega \rangle}-e^{-i\langle x, \omega \rangle}\right|^2(1+\|\omega\|_2^2)^{-m_0}d\omega\nonumber\\
\leq  & C_1\int_{\RR^d}\left|\sum_{j=1}^nu_je^{-i\langle x_j, \omega \rangle}-e^{-i\langle x, \omega \rangle}\right|^2(1+\|\omega\|_2^2)^{-m}d\omega\nonumber\\
\leq & C_2\int_{\RR^d}\left|\sum_{j=1}^nu_je^{-i\langle x_j, \omega \rangle}-e^{-i\langle x, \omega \rangle}\right|^2\mathcal{F}(\Phi)(\omega)d\omega\nonumber\\
= & C_2(\Phi(x-x) - r_m(x)^T(R_m+\mu_mI_n)^{-1}r_m(x)) =: I_3(x),
\end{align}
where the second inequality is because $m\leq m_0$.

Applying Proposition \ref{LEMMAERRORWNUG} to the correlation function $\Phi$ implies that
\begin{align}\label{eq:GpundergI51}
I_3(x) \leq C_3(\mu_m/n)^{1 - \frac{d}{2m}}.
\end{align}
Lemma \ref{thm211}, \eqref{eq:gpunderEG}, \eqref{eq:eI1gpunder}, and \eqref{eq:GpundergI51} imply that for all $t_1>0$, with probability at least $1-\exp(-Ct_1^2)$,
\begin{align}\label{eq:gpunderI1r}
\|G\|_{L_2(\Omega)}^2\leq C_4(1+t_1)^2(\mu_m/n)^{1 - \frac{d}{2m}}.
\end{align}
Similar to the proof of \eqref{eq:gpoverI2r}, for all $t_2\geq C_0$, with probability at least $1-\exp(-t_2^2)$,
\begin{align}\label{eq:gpunderI2r}
\int_{x\in \Omega}I_2(x)dx \leq C_5t_2^2\mu_m^{-\frac{d}{2m}}n^{-\left(1-\frac{d}{2m}\right)}.
\end{align}
By the fact $\|Z-\hat f_G\|_{L_2(\Omega)}^2 \lesssim \|G\|_{L_2(\Omega)}^2 + \|\hat r_{m}(x)^T (R_{m}+\mu_{m} I_n)^{-1} \epsilon\|_{L_2(\Omega)}^2$, \eqref{eq:gpunderI1r}, and \eqref{eq:gpunderI2r}, we finish the proof of Theorem \ref{thm:GPunder}.

\section{Proof of Theorem \ref{thm:GPunderlb}}\label{subsec:pfGPunderlb}

We first present several lemmas.
Lemma \ref{tracelemma} is Lemma F.7 of \cite{wang2020inference}.

\begin{lemma}\label{tracelemma}
Suppose $A,B$ and $C\in \RR^{n\times n}$ are positive definite matrices. We have
\begin{align*}
& \text{tr} ((A + B)(A + B + C)^{-1}) \geq \text{tr} (A(A + C)^{-1}),\\
\text{and } & \text{tr}((A + B)^2(A + B + C)^{-2}) \geq \text{tr} (A^2(A + C)^{-2}).
\end{align*}
\end{lemma}

Let $K$ be a stationary correlation function. Since a correlation function is positive definite, by Mercer's theorem, there exists a countable set of positive eigenvalues $\lambda_1\geq \lambda_2\geq...> 0$ and an orthonormal basis for $L_2(\Omega)$ $\{\varphi_k\}_{k\in\mathbb{N}}$ such that
\begin{align}\label{eq:AppJ1eq1}
K(x - y) = \sum_{k=1}^\infty \lambda_k \varphi_k(x)\varphi_k(y), \quad x,y\in \Omega,
\end{align}
where the summation is uniformly and absolutely convergent.

Lemma \ref{lemDecayEig} states the asymptotic rate of the eigenvalues of $K$, which is implied by the proof of Lemma 18 of \cite{tuo2020kriging}. 

\begin{lemma}\label{lemDecayEig}
Suppose Condition (C1) holds. Suppose $K$ is a stationary correlation function satisfying Condition (C2) and has an expansion as in (\ref{eq:AppJ1eq1}). 
Then, $\lambda_k\asymp k^{-2m_0/d}$.
\end{lemma}

We need the following technical assumption.

\begin{assum}\label{ass:efhasb}
    Suppose there exists a stationary correlation function $K$ satisfying Condition (C2) and a constant $A_0>0$ such that
    \begin{align}\label{Fbound}
        \left\|\frac{\mathcal{F}(K)}{\mathcal{F}(\Psi)}\right\|_{L_\infty(\RR^d)}\leq A_0,
    \end{align}
    and $K$ has an expansion as in \eqref{eq:AppJ1eq1} with eigenfunctions $\|\varphi_k\|_{L_\infty(\Omega)}\leq C$ for all $k=1,2,...$, where $C>0$ not depending on $k$.
\end{assum}

Lemma \ref{LEM:OPTRATE} states that the MSPE
$\mathbb{E}(Z(x)-\hat f_G(x))^2$ can be further bounded by the term related to $K$, and the proof is in Appendix \ref{app:pfcoropt}.

\begin{lemma}\label{LEM:OPTRATE}
    Let $\Psi$ be a correlation function satisfying Condition (C2), and $Z\sim GP(0,\sigma^2\Psi)$. Assume Assumption \ref{ass:efhasb} holds. Let $\{x_1,...,x_n\}$ be a set of design points. Then for all $x\in \Omega$,
	\begin{align*}
	    {\rm Var}[Z(x)|Y] \gtrsim K(x-x) - r_K(x)^T(R_K+\mu I_n)^{-1}r_K(x),
	\end{align*}
	where $R_K = (K(x_j-x_k))_{jk}$, $r_K(x) = (K(x-x_1),...,K(x-x_n))^T$, and $Y$ and $\mu$ are as in \eqref{mean}.
\end{lemma}

Lemma \ref{lemmafixdesign} and Lemma \ref{lemmafixdesignre} state that under fixed designs, the empirical norm is close to the $L_2$ norm. Lemma \ref{lemmafixdesign} can be found in \cite{madych1985estimate,rieger2008sampling}. The proof of Lemma \ref{lemmafixdesignre} is merely repeating the process of proving Lemma \ref{lemmafixdesign} as in \cite{madych1985estimate}, thus is omitted.

\begin{lemma}\label{lemmafixdesign}
	Suppose $g\in H^m(\Omega)$ for some $m > d/2$. Suppose Condition (C4) holds. 
	Then we have
	\begin{align*}
	\|g\|_{L_2(\Omega)} \leq C(h_n^m\|g\|_{H^m(\Omega)} + \|g\|_n)
	\end{align*}
	holds for all $n$, 	where $C$ is a positive constant not depending on $g$ and $n$. 
\end{lemma}
\begin{remark}
	Lemma \ref{lemmafixdesign} is a stronger version of Lemma 3.4 in \cite{utreras1988convergence}. In Lemma 3.4 of \cite{utreras1988convergence}, the fixed designs are assumed to be quasi-uniform. Lemma 3.4 of \cite{utreras1988convergence} is used in \cite{tuo2020improved}.
\end{remark}

\begin{lemma}\label{lemmafixdesignre}
	Suppose $g\in H^m(\Omega)$ for some $m > d/2$. Suppose Condition (C4) holds. 
    Then we have
	\begin{align*}
	\|g\|_n \leq C(h_n^m\|g\|_{H^m(\Omega)} + \|g\|_{L_2(\Omega)})
	\end{align*}
	holds for all $n$, where $C$ is a positive constant not depending on $g$ and $n$.  
\end{lemma}

By Lemma \ref{LEM:OPTRATE}, it suffices to show $$\int_{x\in \Omega} K(x-x) - r_K(x)^T(R_K+\mu I_n)^{-1}r_K(x) dx \gtrsim n^{-\frac{2m_0-d}{2m_0}},$$where $K,r_K,R_K,\mu$ are as in Lemma \ref{LEM:OPTRATE}. This is because for any linear predictor $\hat f_G$,
\begin{align*}
    \mathbb{E}(Z(x) - \hat f_G(x))^2 \geq {\rm Var}[Z(x)|Y],\forall x\in \Omega.
\end{align*}

Notice that for any $u = (u_1,...,u_n)^T \in \RR^n$,
\begin{align}\label{eq:90}
K(x-x) - 2\sum_{j=1}^n u_jK(x-x_j) + \sum_{k=1}^n\sum_{j=1}^n u_k u_jK(x_k-x_j) + \mu \|u\|_2^2 \geq \mu \|u\|_2^2,
\end{align}
since $K$ is positive definite. 
Let $u(x) = (u_1(x),...,u_n(x))^T = (R_K + \mu I_n)^{-1}r_K(x)$, \eqref{eq:90} implies 
\begin{align}\label{u2smallthanlb}
\mu r_K(x)^T(R_K + \mu I_n)^{-2}r_K(x) \leq K(x-x) - r_K(x)^T(R_K + \mu I_n)^{-1}r_K(x).
\end{align}
From \eqref{u2smallthanlb}, it can be seen that it is sufficient to provide a lower bound on $I(x) := r_K(x)^T(R_K + \mu I_n)^{-2}r_K(x)$, because $\mu$ is a constant. 

Let $p = \lfloor n^{d/(2m_0)}\rfloor$, where $\lfloor \cdot \rfloor$ is the floor function. Let $K_1=\frac{1}{\sqrt{n}}(\varphi_1(X),...,\varphi_{p}(X))$, and $K_2=\frac{1}{\sqrt{n}}(\varphi_{p+1}(X),\varphi_{p+2}(X),...)$, where $\varphi_k(X)=(\varphi_k(x_1),...,\varphi_k(x_n))^T$ for $k=1,2,...$, and $\varphi_{k}$'s are as in \eqref{eq:AppJ1eq1}. Let $\Lambda_1=\mbox{diag}(n\lambda_1,...,n\lambda_p)$ and $\Lambda_2=\mbox{diag}(n\lambda_{p+1},...)$, where $\lambda_k$'s are as in \eqref{eq:AppJ1eq1}. Therefore, $R_K=\sum_{k=1}^\infty \lambda_k\varphi_k(X)\varphi_k(X)^T=K_1\Lambda_1K_1^T+K_2\Lambda_2K_2^T$. Note that for any functions $v_1,v_2\in H^{m_0}(\Omega)$, $\|v_1v_2\|_{H^{m_0}(\Omega)}\leq C\|v_1\|_{H^{m_0}(\Omega)}\|v_2\|_{H^{m_0}(\Omega)}$ \citep{adams2003sobolev}. Because $I(x) = u(x)^Tu(x)$, we have
\begin{align}\label{eq:IHm0b}
    & \|I\|_{H^{m_0}(\Omega)} =\left\|\sum_{k=1}^n u_k^2 \right\|_{H^{m_0}(\Omega)}\leq C \sum_{k=1}^n \|u_k\|_{H^{m_0}(\Omega)}^2
    \leq C_1\sum_{k=1}^n \|u_k\|_{\mathcal{N}_K(\Omega)}^2\nonumber\\
    = & C_1\sum_{k=1}^n ((R_K + \mu I_n)^{-1})_k^TR_K ((R_K + \mu I_n)^{-1})_k = C_1{\rm tr}((R_K + \mu I_n)^{-2}R_K),
\end{align}
where the first inequality is by the triangle inequality, and $((R_K + \mu I_n)^{-1})_k$ denotes the $k$-th row of $(R_K + \mu I_n)^{-1}$.

Next we provide a lower bound on $\lambda_{\min}(K_1^TK_1)$ and an upper bound on tr$(K_2^T\Lambda_2K_2)$. For any $a = (a_1,...,a_p)^T\in \RR^p$ such that $\|a\|_2=1$, consider function $g_1 = \sum_{k=1}^p a_k \varphi_k$. Since $\varphi_i$'s are orthonormal, $\|g\|_{L_2(\Omega)}=1$. By Lemma \ref{lemDecayEig}, $\|g_1\|^2_{H^{m_0}(\Omega)}\leq C_1 \|g_1\|^2_{\mathcal{N}_K(\Omega)} \leq C_1 \|\varphi_p\|^2_{\mathcal{N}_K(\Omega)} \leq C_2/\lambda_{p}\leq C_3 p^{2m_0/d} \leq C_3 n$. Therefore, by Lemma \ref{lemmafixdesign} and Condition (C4),
\begin{align*}
    \|g\|_{L_2(\Omega)} \leq &  C_4(h_n^{m_0}\|g_1\|_{H^{m_0}(\Omega)} + \|g\|_n)\nonumber\\
    \leq & C_5(n^{-m_0/d}n^{1/2} + \|g\|_n)
\end{align*}
which implies
\begin{align*}
    \|g\|_n \geq \frac{1}{C_5}\|g\|_{L_2(\Omega)} - n^{(d-2m_0)/2d} \geq \frac{1}{2C_5},  
\end{align*}
for some $n>N_0$, since $m_0>d/2$ and $n^{(d-2m_0)/2d}$ converges to zero. Then we have
\begin{align}\label{eq:lamminb}
    \lambda_{\min}(K_1^TK_1) = \inf_{a\in \RR^p}\|g\|_n^2 \geq C_6, 
\end{align}
for some constant $C_6>0$. 

Considering ${\rm tr} (K_2^T\Lambda_2K_2)$, we have
\begin{align}\label{eq:psi2b}
    &{\rm tr} (K_2^T\Lambda_2K_2) =  \sum_{k=p+1}^\infty \lambda_k\left(\sum_{j=1}^n\varphi_k(x_j)^2\right)
    \leq  nC_7^2 \sum_{k=p+1}^\infty \lambda_k
    \leq  C_8np^{-\frac{2m_0}{d}+1} \leq C_8n^{\frac{d}{2m_0}},
\end{align}
where the first inequality is by Assumption \ref{ass:efhasb}, and the second inequality is by Lemma \ref{lemDecayEig} and the basic inequality $\sum_{k=m}^\infty k^{-2m_0/d}\lesssim m^{-2m_0/d + 1}$.

By \eqref{eq:IHm0b}, $ \|I\|_{H^{m_0}(\Omega)}$ can be further bounded by
\begin{align}\label{eq:IHm0b2}
    \|I\|_{H^{m_0}(\Omega)} \leq & C_1{\rm tr}((K_1\Lambda_1K_1^T+K_2\Lambda_2K_2^T + \mu I_n)^{-2}(K_1\Lambda_1K_1^T+K_2\Lambda_2K_2^T))\nonumber\\
    \leq & C_1({\rm tr}((K_1\Lambda_1K_1^T + \mu I_n)^{-2}K_1\Lambda_1K_1^T)) + C_1\mu^{-2}{\rm tr}(K_2\Lambda_2K_2^T).
\end{align}
Let $I_1 = {\rm tr}((K_1\Lambda_1K_1^T + \mu I_n)^{-2}K_1\Lambda_1K_1^T)$. We have
\begin{align}\label{eq:I1gplb23}
    I_1 = & \sum_{i = 1}^{p} \frac{\lambda_i(K_1\Lambda_1K_1^T)}{(\lambda_i(K_1\Lambda_1K_1^T) + \mu)^2}\leq \mu^{-2} p \leq \mu^{-2} n^{\frac{d}{2m_0}},
\end{align}
where $\lambda_i(K_1\Lambda_1K_1^T)$ denote the $i$-th eigenvalue of $K_1\Lambda_1K_1^T$. 

Combining \eqref{eq:psi2b}, \eqref{eq:IHm0b2}, and \eqref{eq:I1gplb23}, we find that $\|I\|_{H^{m_0}(\Omega)} \leq C_9 n^{\frac{d}{2m_0}}$.

Together with Lemma \ref{lemmafixdesignre}, we have
\begin{align}\label{eq:InandIL2}
    \|I\|_{n}\leq C_{10}(h_n^{m_0}\|I\|_{H^{m_0}(\Omega)} + \|I\|_{L_2(\Omega)}) \leq C_{11}(n^{\frac{d}{2m_0}-\frac{m_0}{d}} + \|I\|_{L_2(\Omega)}).
\end{align}
Let $I_n=\text{tr} (R_K^2 (R_K+ \mu I_n)^{-2})$. Note that $\lambda_i(K_1\Lambda_1K_1^T) = \lambda_i(K_1^TK_1\Lambda_1)$ for $i=1,...,p$,  because if $v_i$ is eigenvector corresponding to $i$-th eigenvalue of $K_1\Lambda_1K_1^T$, then 
\begin{align*}
K_1\Lambda_1K_1^Tv_i = \lambda_i v_i \Rightarrow  K_1^TK_1\Lambda_1K_1^Tv_i = \lambda_i K_1^Tv_i.
\end{align*} 
By Lemma \ref{tracelemma} and $R_K=K_1\Lambda_1K_1^T+K_2\Lambda_2K_2^T$, it can be shown that
\begin{align*}
I_n \geq & \text{tr}((K_1\Lambda_1K_1^T)^2(K_1\Lambda_1K_1^T + \mu I)^{-2}) \nonumber\\
= & \sum_{i = 1}^{p} \bigg(\frac{\lambda_i(K_1\Lambda_1K_1^T)}{\lambda_i(K_1\Lambda_1K_1^T) + \mu}\bigg)^2\nonumber\\
= & \sum_{i = 1}^{p} \bigg(\frac{\lambda_i(\Lambda_1K_1^TK_1)}{\lambda_i(\Lambda_1K_1^TK_1) + \mu}\bigg)^2,
\end{align*}
which implies
\begin{align}\label{thm34lbfangcha2}
I_n \geq & \text{tr}((\Lambda_1K_1^TK_1)^2(\Lambda_1K_1^TK_1 + \mu I)^{-2}) \geq \text{tr}(\Lambda_1^2(\Lambda_1 + \mu (K_1^TK_1)^{-1})^{-2}).
\end{align}
Combining \eqref{eq:lamminb} with \eqref{thm34lbfangcha2}, we conclude that 
\begin{align}\label{thm34lbfangcha3}
I_n  \geq \text{tr}(\Lambda_1^2(\Lambda_1 + C_{12}I_p)^{-2})
\geq \sum_{k=1}^p \frac{\lambda_k^2}{(\lambda_k+C_{12}/n)^2} \geq C_{13}p \geq C_{14} n^{d/(2m_0)},
\end{align}
where the third inequality is by Lemma \ref{lemDecayEig}. 
It follows the Cauchy-Schwarz inequality that
\begin{align*}
    \|I\|_n = \sqrt{\frac{1}{n}\sum_{k=1}^n (u(x_k)^Tu(x_k))^2}\geq \frac{1}{n}\sum_{k=1}^nu(x_k)^Tu(x_k) = I_n/n \geq C_{14} n^{d/(2m_0)-1}.
\end{align*}
By \eqref{eq:InandIL2}, we have
\begin{align*}
    \|I\|_{L_2(\Omega)} \geq \frac{1}{C_{10}}\|I\|_n - n^{d/(2m_0)-\frac{m_0}{d}} \geq C_{15} n^{d/(2m_0)-1} - n^{d/(2m_0)-\frac{m_0}{d}} \gtrsim n^{d/(2m_0)-1},
\end{align*}
for some $n>N_1$ such that $C_{15} n^{d/(2m_0)-1} > 2n^{d/(2m_0)-\frac{m_0}{d}}$, which can be done since $m_0>d$. Thus, for $n>\max(N_0,N_1)$, we have $\mathbb{E}\|Z - \hat f_G\|_{L_2(\Omega)}^2\gtrsim n^{d/(2m_0)-1}$. But for $n\leq\max(N_0,N_1)$, taking $C_{16}=\inf_{n\leq\max(N_0,N_1)} n^{d/(2m_0)-1}/\mathbb{E}\|Z - \hat f_G\|_{L_2(\Omega)}^2$ (which is clearly larger than zero), we can see that $\mathbb{E}\|Z - \hat f_G\|_{L_2(\Omega)}^2\geq C_{16}n^{d/(2m_0)-1}$ for all $n\leq\max(N_0,N_1)$. This finishes the proof.

\section{Proof of Lemma \ref{LEMMA1}}\label{pfLEMMA1}

In this section, we set $m_g:=m_0(g)$ for notational simplicity. Since $g\notin H^{m_g}(\mathbb{R}^d)$, we have
\begin{eqnarray}
& \int_{\mathbb{R}^d}  |\mathcal{F}(g)(\omega)|^2(1+\|\omega\|^2)^{m_g}d\omega = \infty,\label{eq:LEMMA1fou1}\\
{\rm and } & \int_{\mathbb{R}^d} |\mathcal{F}(g)(\omega)|^2(1+\|\omega\|^2)^{m_g-\delta}d\omega < \infty, \forall \delta>0.\label{eq:LEMMA1fou2}
\end{eqnarray}

Using the hyperspherical coordinate transformation, we can represent $\omega$ by a radial coordinate $r$, and $d-1$ angular coordinates $\phi_1,\phi_2,...,\phi_d$. Let $\phi = (\phi_1,\phi_2,...,\phi_d)^T$, and the Jacobian of the transformation be $J$. We can rewrite the left-hand side in (\ref{eq:LEMMA1fou1}) as
\begin{eqnarray*}
\int_0^\infty \int_{[0,2\pi]^{d-1}}|\mathcal{F}(g)(r,\phi)|^2(1+r^2)^{m_g} |\mbox{det}(J)|d\phi dr.
\end{eqnarray*}
Let $g_1(r) = (1+r^2)^{m_g}\int_{[0,2\pi]^{d-1}}|\mathcal{F}(g)(r,\phi)|^2 |\mbox{det}(J)|d\phi.$ Therefore, \eqref{eq:LEMMA1fou1} is equal to $\int_0^{\infty} g_1(r) dr$, which is infinite. It suffices to find an increasing function $Q(r)$ satisfying
\begin{eqnarray}\label{firstcondition_h1}
\int_0^{\infty} \frac{g_1(r)}{Q(r)} dr \leq C_0
\end{eqnarray}
and 
\begin{align}\label{condition_h1}
\lim_{r\rightarrow +\infty} \frac{\log Q(r)}{\log r} = 0,
\end{align}
where $C_0$ is a constant. This is because by \eqref{condition_h1}, we naturally have
\begin{align}\label{anotherconditionneededinLEMMA1}
\int_{\mathbb{R}^d} \frac{|\mathcal{F}(g)(\omega)|^2(1+\|\omega\|_2^2)^{m_g+\delta_1}}{Q(\|\omega\|_2)}d\omega = \infty
\end{align}
for any $\delta_1 > 0$, and more specifically, if \eqref{anotherconditionneededinLEMMA1} is false, then there exists $\delta_1 > 0$ such that 
\begin{align}\label{LEMMA1contraeq1}
\int_{\mathbb{R}^d} \frac{|\mathcal{F}(g)(\omega)|^2(1+\|\omega\|_2^2)^{m_g+\delta_1}}{Q(\|\omega\|_2)}d\omega < \infty.
\end{align}
By \eqref{condition_h1}, there exists a constant $C$ such that for $r > C$, $\frac{\log Q(r)}{\log r} < \delta_1/4$, which is the same as $Q(r) < r^{\delta_1/4}$. This implies that there exists a constant $C_0$ such that $Q(r) < C_0(1 + r^2)^{\delta_1/2}$ for all $r\geq 0$. Therefore, \eqref{LEMMA1contraeq1} yields
\begin{align*}
\infty > & \int_{\mathbb{R}^d} \frac{|\mathcal{F}(g)(\omega)|^2(1+\|\omega\|^2)^{m_g+\delta_1}}{Q(\|\omega\|)}d\omega\\
> &  C_0 \int_{\mathbb{R}^d} |\mathcal{F}(g)(\omega)|^2(1+\|\omega\|^2)^{m_g+\delta_1/2}d\omega = \infty,
\end{align*}
which leads to a contradiction.

We construct $Q(r)$ by the following recurrence way. Let $\alpha_i = 2^{-i}$ for $i\in \mathbb{N}_+$, $\alpha_0=1$, $x_0=0$, and $x_1 = 1$. Let $Q(r) = 1$ for $0 \leq r < x_1$. Since \eqref{eq:LEMMA1fou2} implies that $\int_{0}^\infty g_1(r) r^{-\alpha} dr < \infty$ for any $\alpha > 0$, there exists $x_2 > x_1^{x_1}$ such that $\int_{x_2}^{\infty} g_1(r) r^{-\alpha_1}dr < 1$. Let $Q(r) = x_1^{-\alpha_i} r^{\alpha_i} $ for $r\in (x_1,x_2]$. Suppose we have specified $Q(r)$ for $r\in (x_{i-1}, x_i]$. Clearly there exists an $x_{i+1} > x_i^{x_i}$ such that $\int_{x_{i+1}}^\infty g_1(r) r^{-\alpha_i}dr < 2^{-i}$. Take $Q(r) = Q(x_i) x_i^{-\alpha_i} r^{\alpha_i} $ for $r\in (x_i, x_{i+1}]$. It can be seen that
\begin{align*}
Q(r) = \bigg[\prod_{j=1}^i x_j^{\alpha_{j-1} - \alpha_j}\bigg] r^{\alpha_i}
\end{align*}
and $Q(r)$ is an increasing function. To show that $Q(r)$ satisfies \eqref{firstcondition_h1}, we note that
\begin{align*}
\int_0^{\infty} \frac{g_1(r)}{Q(r)} dr = \sum_{i=1}^\infty \int_{x_{i-1}}^{x_i} \frac{g_1(r)}{Q(r)}dr \leq \sum_{i=1}^\infty \int_{x_{i-1}}^{\infty} \frac{g_1(r)}{Q(r)}dr \leq \int_0^1 g_1(r) dr + \sum_{i=2}^\infty 2^{-i}\leq C_1,
\end{align*}
where the second inequality is because of the choice of $x_i$'s and $Q(r)\geq r^{\alpha_i}$ for $r\in (x_i,x_{i+1}]$.

Next, we show $Q(r)$ satisfies \eqref{condition_h1}. For any $\delta > 0$, there exists an integer $N$ such that for all $r > x_N$, $r\in (x_M,x_{M+1}]$ for some $M>0$,
\begin{align*}
\frac{\log Q(r)}{\log r} = & \frac{\big(\sum_{i=1}^{M-1} (\alpha_{i-1} - \alpha_i) \log x_i\big) +(\alpha_{M-1} - \alpha_M)\log x_M + \alpha_M \log r }{\log r} \\
= & \frac{\big(\sum_{i=1}^{M-1} (\alpha_{i-1} - \alpha_i) \log x_i\big)}{\log r} +\frac{(\alpha_{M-1} - \alpha_M)\log x_M + \alpha_M \log r }{\log r}\\
< & \frac{\big(\sum_{i=1}^{M-1} (\alpha_{i-1} - \alpha_i) \log x_{M-1}\big)}{\log r} +\frac{\alpha_{M-1}\log x_M}{\log r} + \alpha_M\\
< & \frac{1}{x_{M-1}} + \alpha_{M-1} + \alpha_{M}< \frac{1}{x_{N-1}} + 2\alpha_{N-1} < \delta,
\end{align*}
where the second inequality is because $r \geq x_M \geq x_{M-1}^{x_{M-1}}$, and the last inequality is because $x_{N-1} \rightarrow\infty$ and $\alpha_{N-1}\rightarrow 0$.

\section{Proof of Theorem \ref{thm:krrover}}\label{subsec:pfkrrover}
In this section, we set $m_0:=m_0(f)$ for notational simplicity. 
We prove more general results of Theorem \ref{thm:krrover}, as follows. Note that in Theorem \ref{thm:krroverg}, $H^{m_0}(\Omega)$ coincides with $\mathcal{N}_\Psi(\Omega)$.

\begin{theorem}\label{thm:krroverg}
    Suppose the conditions of Theorem \ref{thm:krrover} hold. Suppose $\lambda_m = o(1)$ if $f\in \mathcal{N}_{\Psi}(\Omega)$, and $\lambda_m = o(Q(n)^{-m/m_0})$ if $f\notin \mathcal{N}_{\Psi}(\Omega)$, where $Q(n)$ is as in Lemma \ref{LEMMA1} with $g=f$. Furthermore, suppose $\lambda_m \gtrsim n^{-\frac{2m^2}{d(2m-m_0)}}$. If $f$ has smoothness $m_0$ and $f\in \mathcal{N}_{\Psi}(\Omega)$, for all $t>C_0$ and $n$, with probability at least $1 - C_1\exp(-C_2t^2)$, 
\begin{align}\label{eq:pfpvercaseallX}
	& \|f-\hat f_m\|_{L_2(\Omega)}^2 \leq  CT,
	\mbox{ and }\|\hat f_m\|_{\mathcal{N}_{\Phi}(\Omega)}^2 \leq  C\lambda_m^{-1}T,
	\end{align}
	where
	\begin{align*}
	T = &\max\{t^{\frac{4m}{2m+d}}n^{-\frac{2m}{2m+d}}\lambda_m^{\frac{d(m_0-m)}{m(2m+d)}}\|f\|_{\mathcal{N}_{\Psi}(\Omega)}^{\frac{2d}{2m+d}}, \lambda_m^{\frac{m_0}{m}}\|f\|_{\mathcal{N}_{\Psi}(\Omega)}^2 + 4tn^{-\frac{1}{2}}\lambda_m^{\frac{2m_0-d}{4m}}\|f\|_{\mathcal{N}_{\Psi}(\Omega)},\nonumber\\
	& t^{\frac{4m}{4m-d}}n^{-\frac{2m}{4m-d}} \lambda_m^{\frac{4mm_0-2m_0d-2md}{2m(4m-d)}}\|f\|_{\mathcal{N}_{\Psi}(\Omega)}^{\frac{2(2m-d)}{4m-d}}, t^2n^{-1}\lambda_m^{-\frac{d}{2m}}\}.
	\end{align*}
	If $f$ has smoothness $m_0$ but $f\notin \mathcal{N}_{\Psi}(\Omega)$, for all $t>C_0$, with probability at least $1 - C_1\exp(-C_2t^2)$, 
	we can obtain that 
	\begin{align}\label{eq:pfpvercaseallXX}
	& \|f-\hat f_m\|_{L_2(\Omega)}^2 \leq  CT,
	\mbox{ and }\|\hat f_m\|_{\mathcal{N}_{\Phi}(\Omega)}^2 \leq  C\lambda_m^{-1}T,
	\end{align}
	where
	\begin{align*}
	T = &\max\{t^{\frac{4m}{2m+d}}n^{-\frac{2m}{2m+d}}\lambda_m^{\frac{d(m_0-m)}{m(2m+d)}}Q(n)^{\frac{d}{2m+d}}, \lambda_m^{\frac{m_0}{m}}Q(n) + 4tn^{-\frac{1}{2}}\lambda_m^{\frac{2m_0-d}{4m}}Q(n)^{1/2},\nonumber\\
	& t^{\frac{4m}{4m-d}}n^{-\frac{2m}{4m-d}} \lambda_m^{\frac{4mm_0-2m_0d-2md}{2m(4m-d)}}Q(n)^{\frac{2m-d}{4m-d}}, t^2n^{-1}\lambda_m^{-\frac{d}{2m}}\}.
	\end{align*}
\end{theorem}

Note that Theorem \ref{thm:krrover} can be obtained by taking $\lambda_m \asymp n^{-\frac{2m}{2m_0+d}}$ in Theorem \ref{thm:krroverg}. 

Now we begin to prove Theorem \ref{thm:krroverg}. The following lemmas are used. Lemma \ref{lemmaefsmall} is Lemma A.1 in \cite{tuo2020improved}, which states that the inner product $\langle \epsilon, g \rangle_n$ is small; also see Lemma 8.4 of \cite{geer2000empirical}.
\begin{lemma}\label{lemmaefsmall}
	Suppose Condition (C5) holds. Let $K$ be a kernel function, which is stationary, positive definite and intergrable on $\RR^d$. Suppose there exist constants $c_2 \geq c_1>0$ and $m> d/2$ such that, for all $\omega\in\mathbb{R}^d$,
	$$ c_1(1+\|\omega\|_2^2)^{-m} \leq  \mathcal{F}(K)(\omega)\leq c_2(1+\|\omega\|_2^2)^{-m}. $$
	Then for all $t > C$, with probability at least $1 - C_1\exp(-C_2 t^2)$,
	\begin{align*}
	\sup_{g\in \mathcal{N}_{K}(\Omega)}\frac{|\langle \epsilon, g \rangle_n|}{\|g\|_n^{1 - \frac{d}{2m}}\|g\|_{\mathcal{N}_{K}(\Omega)}^{\frac{d}{2m}}} \leq tn^{-\frac{1}{2}}.
	\end{align*}
\end{lemma}

Lemma \ref{LEM:FSTAR} states the solution to the expectation version of \eqref{KRRest1} obtained by replacing $\|\cdot\|_n$ with $\|\cdot\|_{L_2(\Omega)}$, denoted by $f^*$, can approximate $f$ well. The proof of Lemma \ref{LEM:FSTAR} can be found in Appendix \ref{app:pflemmafstar}.
\begin{lemma}\label{LEM:FSTAR}
	Suppose the conditions in Theorem \ref{thm:krroverg} hold. Let $f^*$ be the solution to the optimization problem 
	\begin{align}\label{eq:overoptOmega}
	\min_{\tilde f \in \mathcal{N}_{\Phi}(\Omega)}\|f-\tilde f\|_{L_2(\Omega)}^2 + \lambda_m \|\tilde f\|^2_{\mathcal{N}_{\Phi}(\Omega)}.
	\end{align}
	If $f\in H^{m_0}(\Omega)$, then 
	\begin{align*}
	\|f-f^*\|_{L_2(\Omega)}^2 + \lambda_m \|f^*\|^2_{\mathcal{N}_{\Phi}(\Omega)} \leq C\lambda_m^{\frac{m_0}{m}}\|f\|_{\mathcal{N}_{\Psi}(\Omega)}^2,
	\end{align*}
	where $C$ is a constant only depending on $\Omega$, $\Phi$, and $\Psi$ including $m$ and $m_0$. In particular, we have
	\begin{align}\label{eq:fstarbound}
	\|f-f^*\|_{L_2}^2 \leq C\lambda_m^{\frac{m_0}{m}}\|f\|_{\mathcal{N}_{\Psi}(\Omega)}^2, \mbox{ and } \|f^*\|_{\mathcal{N}_{\Phi}(\Omega)}^2 \leq C\lambda_m^{\frac{m_0-m}{m}}\|f\|_{\mathcal{N}_{\Psi}(\Omega)}^2. 
	\end{align}
	If $f\notin H^{m_0}(\Omega)$, then there exists an increasing $Q:\RR_+ \mapsto\RR_+$ such that
	\begin{align*}
	\|f-f^*\|_{L_2(\Omega)}^2 + \lambda_m \|f^*\|^2_{\mathcal{N}_{\Phi}(\Omega)} \leq C\lambda_m^{\frac{m_0}{m}}Q(n),
	\end{align*}
	where $C$ is a constant only depending on $\Omega$, $\Phi$, and $\Psi$ including $m$ and $m_0$. 
	In particular, we have
	\begin{align}\label{eq:fstarboundX}
	\|f-f^*\|_{L_2}^2 \leq C\lambda_m^{\frac{m_0}{m}}Q(n), \mbox{ and } \|f^*\|_{\mathcal{N}_{\Phi}(\Omega)}^2 \leq C\lambda_m^{\frac{m_0-m}{m}}Q(n).
	\end{align}
\end{lemma}

\textit{Proof of Theorem \ref{thm:krroverg}.} The proof of Theorem \ref{thm:krroverg} consists of two parts, according to $f$ lies in $\mathcal{N}_\Psi(\Omega)$ or not.

\noindent\textbf{Case 1: $f$ has smoothness $m_0$ and $f\in \mathcal{N}_\Psi(\Omega)$.}

We first consider the case that $f$ has smoothness $m_0$ and $f\in \mathcal{N}_\Psi(\Omega)$. Let $f^*$ be as in Lemma \ref{LEM:FSTAR}. Because $\hat f_m$ is the solution to \eqref{KRRest1}, we have 
\begin{align}\label{eq:basicineq1}
\|y-\hat f_m\|_n^2 + \lambda_m\|\hat f_m\|_{\mathcal{N}_{\Phi}(\Omega)}^2 \leq \|y-f^*\|_n^2 + \lambda_m\|f^*\|_{\mathcal{N}_{\Phi}(\Omega)}^2,
\end{align}
where $y=(y_1,...,y_n)^T$. By rearrangement, \eqref{eq:basicineq1} yields the \textit{basic inequality}
\begin{align}\label{eq:basicineq2}
\|f-\hat f_m\|_n^2 + \lambda_m\|\hat f_m\|_{\mathcal{N}_{\Phi}(\Omega)}^2 \leq \|f-f^*\|_n^2 + \lambda_m\|f^*\|_{\mathcal{N}_{\Phi}(\Omega)}^2 + 2\langle \epsilon,\hat f_m - f^*\rangle_n.
\end{align}
We apply Lemma \ref{lemmaefsmall} to $\langle \epsilon,\hat f_m - f^*\rangle_n$ and obtain that with probability at least $1 - C_1\exp(-C_2 t^2)$, 
\begin{align}\label{eq:efdsmall1}
\langle \epsilon,\hat f_m - f^*\rangle_n \leq tn^{-\frac{1}{2}} \|\hat f_m - f^*\|_n^{1-\frac{d}{2m}}\|\hat f_m - f^*\|_{\mathcal{N}_{\Phi}(\Omega)}^{\frac{d}{2m}}.
\end{align}
Plugging \eqref{eq:efdsmall1} into \eqref{eq:basicineq2}, the inequality
\begin{align}\label{eq:overineq1}
& \|f-\hat f_m\|_n^2 + \lambda_m\|\hat f\|_{\mathcal{N}_{\Phi}(\Omega)}^2\nonumber\\
\leq & \|f-f^*\|_n^2 + \lambda_m\|f^*\|_{\mathcal{N}_{\Phi}(\Omega)}^2 + 2tn^{-\frac{1}{2}} \|\hat f_m - f^*\|_n^{1-\frac{d}{2m}}\|\hat f_m - f^*\|_{\mathcal{N}_{\Phi}(\Omega)}^{\frac{d}{2m}}\nonumber\\
\leq & \|f-f^*\|_n^2 + \lambda_m\|f^*\|_{\mathcal{N}_{\Phi}(\Omega)}^2 + 2tn^{-\frac{1}{2}} \|f - \hat f_m \|_n^{1-\frac{d}{2m}}\|\hat f_m - f^*\|_{\mathcal{N}_{\Phi}(\Omega)}^{\frac{d}{2m}}\nonumber\\
& + 2tn^{-\frac{1}{2}} \|f - f^*\|_n^{1-\frac{d}{2m}}\|\hat f_m - f^*\|_{\mathcal{N}_{\Phi}(\Omega)}^{\frac{d}{2m}}
\end{align}
holds with probability at least $1 - C_1\exp(-C_2 t^2)$. The last inequality in \eqref{eq:overineq1} is because of the triangle inequality and the basic inequality $(a+b)^q \leq a^q + b^q$ for any $a,b\geq 0$ and $q\in [0,1]$.

By the Gagliardo–Nirenberg interpolation inequality,
\begin{align}\label{eq:GNinter}
\|f^*\|_{\mathcal{N}_{\Psi}(\Omega)} \lesssim \|f^*\|_{H^{m_0}(\Omega)}\lesssim \|f^*\|_{L_2(\Omega)}^{\frac{m-m_0}{m}}\|f^*\|_{H^m(\Omega)}^{\frac{m_0}{m}} \lesssim \|f^*\|_{L_2(\Omega)}^{\frac{m-m_0}{m}}\|f^*\|_{\mathcal{N}_{\Phi}(\Omega)}^{\frac{m_0}{m}},
\end{align}
where the first and last inequalities are because $\|\cdot\|_{\mathcal{N}_{\Psi}(\Omega)}$ and $\|\cdot\|_{\mathcal{N}_{\Phi}(\Omega)}$ are equivalent to $\|\cdot\|_{H^{m_0}(\Omega)}$ and $\|\cdot\|_{H^m(\Omega)}$, respectively. It can be seen from \eqref{eq:GNinter} that $f - f^*\in \mathcal{N}_{\Psi}(\Omega)$, which is equivalent to $f - f^*\in H^{m_0}(\Omega)$. By Lemma \ref{lemmafixdesignre}, we have
\begin{align}\label{eq:diffofff1overfix}
\|f-f^*\|_n \lesssim & h_n^{m_0}\|f-f^*\|_{N_{\Psi}(\Omega)} + \|f-f^*\|_{L_2(\Omega)}\nonumber\\
\lesssim & h_n^{m_0}\|f\|_{N_{\Psi}(\Omega)} +h_n^{m_0}\|f^*\|_{N_{\Psi}(\Omega)} + \|f-f^*\|_{L_2(\Omega)}\nonumber\\
\lesssim & h_n^{m_0}\|f\|_{N_{\Psi}(\Omega)} +h_n^{m_0}\|f^*\|_{L_2(\Omega)}^{\frac{m-m_0}{m}}\|f^*\|_{\mathcal{N}_{\Phi}(\Omega)}^{\frac{m_0}{m}} + \|f-f^*\|_{L_2(\Omega)},
\end{align}
where the second inequality is because of the triangle inequality, and the third inequality is because of \eqref{eq:GNinter}.

The reproducing property \citep{wendland2004scattered} implies that for any $x\in \Omega$,
\begin{align}\label{eq:reproP}
    |f(x)| = |\langle f,\Psi(x-\cdot)\rangle_{\mathcal{N}_{\Psi}(\Omega)}| \leq \|f\|_{\mathcal{N}_{\Psi}(\Omega)}\Psi(x-x),
\end{align}
which implies $\|f\|_{L_2(\Omega)} \lesssim \|f\|_{L_\infty(\Omega)}\lesssim \|f\|_{\mathcal{N}_{\Psi}(\Omega)}$. By Lemma \ref{LEM:FSTAR}, we obtain
\begin{align*}
\|f^*\|_{L_2(\Omega)} \leq \|f\|_{L_2(\Omega)} + \|f-f^*\|_{L_2(\Omega)} \lesssim \|f\|_{L_2(\Omega)}\lesssim \|f\|_{\mathcal{N}_{\Psi}(\Omega)},
\end{align*}
which, together with \eqref{eq:fstarbound} and \eqref{eq:diffofff1overfix}, yields \begin{align}\label{eq:1}
\|f-f^*\|_n \lesssim & h_n^{m_0}\|f\|_{N_{\Psi}(\Omega)} +h_n^{m_0}(C\lambda_m^{\frac{m_0-m}{2m}})^{\frac{m_0}{m}}\|f\|_{\mathcal{N}_{\Psi}(\Omega)} + \lambda_m^{\frac{m_0}{2m}}\|f\|_{\mathcal{N}_{\Psi}(\Omega)}\nonumber\\
\lesssim & (h_n^{m_0} +h_n^{m_0}\lambda_m^{\frac{(m_0-m)m_0}{2m^2}}+ \lambda_m^{\frac{m_0}{2m}})\|f\|_{\mathcal{N}_{\Psi}(\Omega)}\nonumber\\
\lesssim & (h_n^{m_0}\lambda_m^{\frac{(m_0-m)m_0}{2m^2}}+ \lambda_m^{\frac{m_0}{2m}})\|f\|_{\mathcal{N}_{\Psi}(\Omega)},
\end{align}
where the second inequality is because $\lambda_m\leq C_3$ and $m\geq m_0$. By Condition (C4) and the condition $\lambda_m \gtrsim n^{-\frac{2m^2}{d(2m-m_0)}}$, it can be checked that $h_n^{m_0}\lambda_m^{\frac{(m_0-m)m_0}{2m^2}}\lesssim \lambda_m^{\frac{m_0}{2m}}$. Therefore, \eqref{eq:1} can be further bounded by
\begin{align}\label{eq:ffstarnnorm}
\|f-f^*\|_n \lesssim  (h_n^{m_0}\lambda_m^{\frac{(m_0-m)m_0}{2m^2}}+ \lambda_m^{\frac{m_0}{2m}})\|f\|_{\mathcal{N}_{\Psi}(\Omega)} \lesssim \lambda_m^{\frac{m_0}{2m}}\|f\|_{\mathcal{N}_{\Psi}(\Omega)}.
\end{align}

Plugging \eqref{eq:ffstarnnorm} into \eqref{eq:overineq1}, we have that with probability at least $1 - C_1\exp(-C_2t^2)$, 
\begin{align}\label{eq:overineq2}
& \|f-\hat f_m\|_n^2 + \lambda_m\|\hat f_m\|_{\mathcal{N}_{\Phi}(\Omega)}^2\nonumber\\
\lesssim & \lambda_m^{\frac{m_0}{m}}\|f\|_{\mathcal{N}_{\Psi}(\Omega)}^2 + \lambda_m\|f^*\|_{\mathcal{N}_{\Phi}(\Omega)}^2 + 2tn^{-\frac{1}{2}} \|f - \hat f_m \|_n^{1-\frac{d}{2m}}\|\hat f_m - f^*\|_{\mathcal{N}_{\Phi}(\Omega)}^{\frac{d}{2m}}\nonumber\\
& + 2tn^{-\frac{1}{2}} (\lambda_m^{\frac{m_0}{2m}}\|f\|_{\mathcal{N}_{\Psi}(\Omega)})^{1-\frac{d}{2m}}\|\hat f_m - f^*\|_{\mathcal{N}_{\Phi}(\Omega)}^{\frac{d}{2m}}\nonumber\\
\lesssim & \lambda_m^{\frac{m_0}{m}}\|f\|_{\mathcal{N}_{\Psi}(\Omega)}^2 + C\lambda_m\lambda_m^{\frac{m_0-m}{m}}\|f\|_{\mathcal{N}_{\Psi}(\Omega)}^2 + 2tn^{-\frac{1}{2}} \|f - \hat f_m \|_n^{1-\frac{d}{2m}}\|\hat f_m - f^*\|_{\mathcal{N}_{\Phi}(\Omega)}^{\frac{d}{2m}}\nonumber\\
& + 2tn^{-\frac{1}{2}} (C\lambda_m^{\frac{m_0}{2m}}\|f\|_{\mathcal{N}_{\Psi}(\Omega)})^{1-\frac{d}{2m}}\|\hat f_m - f^*\|_{\mathcal{N}_{\Phi}(\Omega)}^{\frac{d}{2m}}\nonumber\\
\lesssim & \lambda_m^{\frac{m_0}{m}}\|f\|_{\mathcal{N}_{\Psi}(\Omega)}^2 + 2tn^{-\frac{1}{2}} \|f - \hat f_m \|_n^{1-\frac{d}{2m}}\|\hat f_m\|_{\mathcal{N}_{\Phi}(\Omega)}^{\frac{d}{2m}}\nonumber\\
&  + 2tn^{-\frac{1}{2}} \|f - \hat f_m \|_n^{1-\frac{d}{2m}}\| f^*\|_{\mathcal{N}_{\Phi}(\Omega)}^{\frac{d}{2m}}+ 2tn^{-\frac{1}{2}} (\lambda_m^{\frac{m_0}{2m}}\|f\|_{\mathcal{N}_{\Psi}(\Omega)})^{1-\frac{d}{2m}}\|\hat f_m\|_{\mathcal{N}_{\Phi}(\Omega)}^{\frac{d}{2m}}\nonumber\\
& + 2tn^{-\frac{1}{2}} (\lambda_m^{\frac{m_0}{2m}}\|f\|_{\mathcal{N}_{\Psi}(\Omega)})^{1-\frac{d}{2m}}\| f^*\|_{\mathcal{N}_{\Phi}(\Omega)}^{\frac{d}{2m}},
\end{align}
where the second inequality is because of Lemma \ref{LEM:FSTAR}, and the third inequality is because of the triangle inequality $\|\hat f_m - f^*\|_{\mathcal{N}_{\Phi}(\Omega)}\leq \|\hat f_m\|_{\mathcal{N}_{\Phi}(\Omega)} +\|f^*\|_{\mathcal{N}_{\Phi}(\Omega)}$ and the basic inequality $(a+b)^q \leq a^q + b^q$ for any $a,b\geq 0$ and $q\in [0,1]$.

Next, we consider two subcases.

\noindent\textbf{Case 1.1: $\|\hat f_m\|_{\mathcal{N}_{\Phi}(\Omega)}\leq \|f^*\|_{\mathcal{N}_{\Phi}(\Omega)}$.} Then \eqref{eq:overineq2} implies
\begin{align*}
\|f-\hat f_m\|_n^2 + \lambda_m\|\hat f_m\|_{\mathcal{N}_{\Phi}(\Omega)}^2
\lesssim & \lambda_m^{\frac{m_0}{m}}\|f\|_{\mathcal{N}_{\Psi}(\Omega)}^2 + 4tn^{-\frac{1}{2}} \|f - \hat f_m \|_n^{1-\frac{d}{2m}}\| f^*\|_{\mathcal{N}_{\Phi}(\Omega)}^{\frac{d}{2m}}\nonumber\\
& + 4tn^{-\frac{1}{2}} (\lambda_m^{\frac{m_0}{2m}}\|f\|_{\mathcal{N}_{\Psi}(\Omega)})^{1-\frac{d}{2m}}\| f^*\|_{\mathcal{N}_{\Phi}(\Omega)}^{\frac{d}{2m}},
\end{align*}
which implies either
\begin{align}\label{eq:pfovercase111}
\|f-\hat f_m\|_n^2 + \lambda_m\|\hat f_m\|_{\mathcal{N}_{\Phi}(\Omega)}^2
\lesssim \lambda_m^{\frac{m_0}{m}}\|f\|_{\mathcal{N}_{\Psi}(\Omega)}^2 + tn^{-\frac{1}{2}} (\lambda_m^{\frac{m_0}{2m}}\|f\|_{\mathcal{N}_{\Psi}(\Omega)})^{1-\frac{d}{2m}}\| f^*\|_{\mathcal{N}_{\Phi}(\Omega)}^{\frac{d}{2m}},
\end{align}
or 
\begin{align}\label{eq:pfovercase112}
& \|f-\hat f_m\|_n^2 + \lambda_m\|\hat f_m\|_{\mathcal{N}_{\Phi}(\Omega)}^2   \lesssim tn^{-\frac{1}{2}} \|f - \hat f_m \|_n^{1-\frac{d}{2m}}\| f^*\|_{\mathcal{N}_{\Phi}(\Omega)}^{\frac{d}{2m}}\nonumber\\
\lesssim & tn^{-\frac{1}{2}} \|f - \hat f_m \|_n^{1-\frac{d}{2m}}(\lambda_m^{\frac{m_0-m}{2m}}\|f\|_{\mathcal{N}_{\Psi}(\Omega)})^{\frac{d}{2m}}.
\end{align}
By Lemma \ref{LEM:FSTAR}, \eqref{eq:pfovercase111} implies
\begin{align}\label{eq:pfovercase113} 
& \|f-\hat f_m\|_n^2 + \lambda_m\|\hat f_m\|_{\mathcal{N}_{\Phi}(\Omega)}^2\nonumber\\
\lesssim & \lambda_m^{\frac{m_0}{m}}\|f\|_{\mathcal{N}_{\Psi}(\Omega)}^2 + tn^{-\frac{1}{2}} (\lambda_m^{\frac{m_0}{2m}}\|f\|_{\mathcal{N}_{\Psi}(\Omega)})^{1-\frac{d}{2m}}(\lambda_m^{\frac{m_0-m}{2m}}\|f\|_{\mathcal{N}_{\Psi}(\Omega)})^{\frac{d}{2m}}\nonumber\\
= & \lambda_m^{\frac{m_0}{m}}\|f\|_{\mathcal{N}_{\Psi}(\Omega)}^2 + tn^{-\frac{1}{2}}\lambda_m^{\frac{2m_0-d}{4m}}\|f\|_{\mathcal{N}_{\Psi}(\Omega)}.
\end{align}
Solving \eqref{eq:pfovercase112} leads to
\begin{align}\label{eq:pfovercase114} 
\|f-\hat f_m\|_n^2 \lesssim & t^{\frac{4m}{2m+d}}n^{-\frac{2m}{2m+d}}\lambda_m^{\frac{d(m_0-m)}{m(2m+d)}}\|f\|_{\mathcal{N}_{\Psi}(\Omega)}^{\frac{2d}{2m+d}},\nonumber\\
\mbox{ and } \|\hat f_m\|_{\mathcal{N}_{\Phi}(\Omega)}^2  \lesssim & t^{\frac{4m}{2m+d}}n^{-\frac{2m}{2m+d}}\lambda_m^{\frac{d(m_0-m)(2m-d)}{4m^2(2m+d)}-1} \lambda_m^{\frac{d(m_0-m)}{4m^2}}\|f\|_{\mathcal{N}_{\Psi}(\Omega)}^{\frac{2d}{2m+d}},\nonumber\\
= & t^{\frac{4m}{2m+d}}n^{-\frac{2m}{2m+d}}\lambda_m^{\frac{d(m_0-m)}{m(2m+d)}-1}\|f\|_{\mathcal{N}_{\Psi}(\Omega)}^{\frac{2d}{2m+d}},.
\end{align}
It can be seen that \eqref{eq:pfovercase113} yields
\begin{align}\label{eq:pfovercase115} 
& \|f-\hat f_m\|_n^2\lesssim \lambda_m^{\frac{m_0}{m}}\|f\|_{\mathcal{N}_{\Psi}(\Omega)}^2 + tn^{-\frac{1}{2}}\lambda_m^{\frac{2m_0-d}{4m}}\|f\|_{\mathcal{N}_{\Psi}(\Omega)},\nonumber\\
\mbox{ and } &  \|\hat f_m\|_{\mathcal{N}_{\Phi}(\Omega)}^2
\lesssim \lambda_m^{\frac{m_0-m}{m}}\|f\|_{\mathcal{N}_{\Psi}(\Omega)}^2 + tn^{-\frac{1}{2}}\lambda_m^{\frac{2m_0-d -4m}{4m}}\|f\|_{\mathcal{N}_{\Psi}(\Omega)}.
\end{align}

\noindent\textbf{Case 1.2: $\|\hat f_m\|_{\mathcal{N}_{\Phi}(\Omega)}> \|f^*\|_{\mathcal{N}_{\Phi}(\Omega)}$.} Then \eqref{eq:overineq2} implies
\begin{align*}
\|f-\hat f_m\|_n^2 + \lambda_m\|\hat f_m\|_{\mathcal{N}_{\Phi}(\Omega)}^2
\lesssim & \lambda_m^{\frac{m_0}{m}}\|f\|_{\mathcal{N}_{\Psi}(\Omega)}^2 + 4tn^{-\frac{1}{2}} \|f - \hat f_m \|_n^{1-\frac{d}{2m}}\|\hat f_m\|_{\mathcal{N}_{\Phi}(\Omega)}^{\frac{d}{2m}}\nonumber\\
& + 4tn^{-\frac{1}{2}} (\lambda_m^{\frac{m_0}{2m}}\|f\|_{\mathcal{N}_{\Psi}(\Omega)})^{1-\frac{d}{2m}}\|\hat f_m\|_{\mathcal{N}_{\Phi}(\Omega)}^{\frac{d}{2m}},
\end{align*}
which implies either
\begin{align}\label{eq:pfovercase121}
& \|f-\hat f_m\|_n^2 + \lambda_m\|\hat f_m\|_{\mathcal{N}_{\Phi}(\Omega)}^2\nonumber\\
\lesssim & \lambda_m^{\frac{m_0}{m}}\|f\|_{\mathcal{N}_{\Psi}(\Omega)}^2 + 4tn^{-\frac{1}{2}} (\lambda_m^{\frac{m_0}{2m}}\|f\|_{\mathcal{N}_{\Psi}(\Omega)})^{1-\frac{d}{2m}}\|\hat f_m\|_{\mathcal{N}_{\Phi}(\Omega)}^{\frac{d}{2m}},
\end{align}
or 
\begin{align}\label{eq:pfovercase122}
& \|f-\hat f_m\|_n^2 + \lambda_m\|\hat f_m\|_{\mathcal{N}_{\Phi}(\Omega)}^2 \lesssim tn^{-\frac{1}{2}} \|f - \hat f_m \|_n^{1-\frac{d}{2m}}\|\hat f_m\|_{\mathcal{N}_{\Phi}(\Omega)}^{\frac{d}{2m}}.
\end{align}
It can be seen that \eqref{eq:pfovercase121} implies either
\begin{align}\label{eq:pfovercase1211}
& \|f-\hat f_m\|_n^2 + \lambda_m\|\hat f_m\|_{\mathcal{N}_{\Phi}(\Omega)}^2  \lesssim \lambda_m^{\frac{m_0}{m}}\|f\|_{\mathcal{N}_{\Psi}(\Omega)}^2
\end{align}
or
\begin{align}\label{eq:pfovercase1212}
& \|f-\hat f_m\|_n^2 + \lambda_m\|\hat f_m\|_{\mathcal{N}_{\Phi}(\Omega)}^2  \lesssim tn^{-\frac{1}{2}} (\lambda_m^{\frac{m_0}{2m}}\|f\|_{\mathcal{N}_{\Psi}(\Omega)})^{1-\frac{d}{2m}}\|\hat f_m\|_{\mathcal{N}_{\Phi}(\Omega)}^{\frac{d}{2m}}.
\end{align}
Solving \eqref{eq:pfovercase1211} leads to
\begin{align}\label{eq:pfovercase1213}
\|f-\hat f_m\|_n^2 \lesssim \lambda_m^{\frac{m_0}{m}}\|f\|_{\mathcal{N}_{\Psi}(\Omega)}^2, \mbox{ and } \|\hat f_m\|_{\mathcal{N}_{\Phi}(\Omega)}^2 \lesssim \lambda_m^{\frac{m_0-m}{m}}\|f\|_{\mathcal{N}_{\Psi}(\Omega)}^2.
\end{align}
Solving \eqref{eq:pfovercase1212} yields
\begin{align}\label{eq:pfovercase1214}
\|f-\hat f_m\|_n^2 \lesssim \lambda_m^{-\frac{2mm_0-m_0d-md}{(4m-d)m}}t^{\frac{4m}{4m-d}}n^{-\frac{2m}{4m-d}}\|f\|_{\mathcal{N}_{\Psi}(\Omega)}^{\frac{2(2m-d)}{4m-d}},\nonumber\\
\mbox{ and } \|\hat f_m\|_{\mathcal{N}_{\Phi}(\Omega)} \lesssim \lambda_m^{-\frac{2mm_0-m_0d-md}{(4m-d)m}-1}t^{\frac{4m}{4m-d}}n^{-\frac{2m}{4m-d}}\|f\|_{\mathcal{N}_{\Psi}(\Omega)}^{\frac{2(2m-d)}{4m-d}}.
\end{align}
Solving \eqref{eq:pfovercase122} yields
\begin{align}\label{eq:pfovercase1215}
\|f-\hat f_m\|_n \lesssim tn^{-\frac{1}{2}}\lambda_m^{-\frac{d}{4m}} \mbox{ and }\|\hat f_m\|_{\mathcal{N}_{\Phi}(\Omega)} \lesssim tn^{-\frac{1}{2}}\lambda_m^{-\frac{2m+d}{4m}}.
\end{align}
Combining all the cases listed in \eqref{eq:pfovercase114}, \eqref{eq:pfovercase115}, \eqref{eq:pfovercase1213}, \eqref{eq:pfovercase1214} and \eqref{eq:pfovercase1215}, we have
\begin{align}\label{eq:pfpvercaseall}
\|f-\hat f_m\|_n^2 \lesssim &  T, \mbox{ and }\|\hat f_m\|_{\mathcal{N}_{\Phi}(\Omega)}^2 \lesssim \lambda_m^{-1}T,
\end{align}
where
\begin{align*}
	T = &\max\{t^{\frac{4m}{2m+d}}n^{-\frac{2m}{2m+d}}\lambda_m^{\frac{d(m_0-m)}{m(2m+d)}}\|f\|_{\mathcal{N}_{\Psi}(\Omega)}^{\frac{2d}{2m+d}}, \lambda_m^{\frac{m_0}{m}}\|f\|_{\mathcal{N}_{\Psi}(\Omega)}^2 + tn^{-\frac{1}{2}}\lambda_m^{\frac{2m_0-d}{4m}}\|f\|_{\mathcal{N}_{\Psi}(\Omega)},\nonumber\\
	& t^{\frac{4m}{4m-d}}n^{-\frac{2m}{4m-d}} \lambda_m^{\frac{2mm_0-m_0d-md}{m(4m-d)}}\|f\|_{\mathcal{N}_{\Psi}(\Omega)}^{\frac{2(2m-d)}{4m-d}}, t^2n^{-1}\lambda_m^{-\frac{d}{2m}}\}.
\end{align*}
It remains to bound the difference between $\|f-\hat f_m\|_n$ and $\|f-\hat f_m\|_{L_2(\Omega)}$, which can be done by applying Lemma \ref{lemmafixdesign}. To this end, note that $f^*-\hat f_m\in \mathcal{N}_{\Phi}(\Omega)$, which is equivalent to $f^*-\hat f_m\in H^m(\Omega)$. By the triangle inequality and Lemma \ref{LEM:FSTAR},
\begin{align}\label{eq:pfpvercaseL2fix}
\|f-\hat f_m\|_{L_2(\Omega)}\leq & \|f-f^*\|_{L_2(\Omega)} + \|f^*-\hat f_m\|_{L_2(\Omega)}\nonumber\\
\lesssim & \lambda_m^{\frac{m_0}{2m}}\|f\|_{\mathcal{N}_{\Psi}(\Omega)} + h_n^{m} \|f^*-\hat f_m\|_{H^{m}(\Omega)} + \|f^*-\hat f_m\|_n\nonumber\\
\lesssim & \lambda_m^{\frac{m_0}{2m}}\|f\|_{\mathcal{N}_{\Psi}(\Omega)} + h_n^{m} \|f^*-\hat f_m\|_{\mathcal{N}_{\Phi}(\Omega)} + \|f^*-\hat f_m\|_n\nonumber\\
\lesssim & \lambda_m^{\frac{m_0}{2m}}\|f\|_{\mathcal{N}_{\Psi}(\Omega)} + h_n^{m} \|f^*\|_{\mathcal{N}_{\Phi}(\Omega)} + h_n^{m} \|\hat f_m\|_{\mathcal{N}_{\Phi}(\Omega)} + \|f-f^*\|_n +\|f - \hat f_m\|_n\nonumber\\
\lesssim & \lambda_m^{\frac{m_0}{2m}}\|f\|_{\mathcal{N}_{\Psi}(\Omega)} + h_n^{m}\lambda_m^{\frac{m_0-m}{2m}}\|f\|_{\mathcal{N}_{\Psi}(\Omega)} + h_n^{m}\|\hat f_m\|_{\mathcal{N}_{\Phi}(\Omega)} +\|f - \hat f_m\|_n\nonumber\\
\lesssim & \lambda_m^{\frac{m_0}{2m}}\|f\|_{\mathcal{N}_{\Psi}(\Omega)} +h_n^{m}\lambda_m^{\frac{m_0-m}{2m}}\|f\|_{\mathcal{N}_{\Psi}(\Omega)} + h_n^{m}\lambda_m^{-1/2}T^{1/2} + T^{1/2}\nonumber\\
\lesssim & T^{1/2}.
\end{align}
The second inequality is by Lemma \ref{LEM:FSTAR}; the third inequality is by the equivalence of $H^m(\Omega)$ and $\mathcal{N}_{\Phi}(\Omega)$; the fourth inequality is by the triangle inequality; the fifth inequality is by Lemma \ref{LEM:FSTAR} and \eqref{eq:ffstarnnorm}; the sixth inequality is by \eqref{eq:pfpvercaseall}; the last inequality is because of Condition (C4) and the condition $\lambda_m \gtrsim n^{-\frac{2m^2}{d(2m-m_0)}} \geq n^{-\frac{2m}{d}}$. Then the desired results of Case 1 follow from \eqref{eq:pfpvercaseall} and \eqref{eq:pfpvercaseL2fix}.

\noindent\textbf{Case 2: $f$ has smoothness $m_0$ and $f\notin \mathcal{N}_\Psi(\Omega)$.}

Let $f^*$ be as in Lemma \ref{LEM:FSTAR}. Let $\delta = \frac{m_0^2(m_0-d/2)}{2(m_0+m)(m_0-d/2)+m_0^2}>0$, we have $f\in H^{m_0-\delta}$ and $m_0-\delta > d/2$. Similar to the proof in Case 1, we can change \eqref{eq:diffofff1overfix} to 
\begin{align}\label{eq:diffofff1overfixX}
\|f-f^*\|_n \lesssim & h_n^{m_0-\delta}\|f-f^*\|_{H^{m_0-\delta}(\Omega)} + \|f-f^*\|_{L_2(\Omega)}\nonumber\\
\lesssim & h_n^{m_0-\delta}\|f\|_{H^{m_0-\delta}(\Omega)} +h_n^{m_0-\delta}\|f^*\|_{H^{m_0-\delta}(\Omega)} + \|f-f^*\|_{L_2(\Omega)}\nonumber\\
\lesssim & h_n^{m_0-\delta}\|f\|_{H^{m_0-\delta}(\Omega)} +h_n^{m_0-\delta}\|f^*\|_{L_2(\Omega)}^{\frac{m-m_0+\delta}{m}}\|f^*\|_{\mathcal{N}_{\Phi}(\Omega)}^{\frac{m_0-\delta}{m}} + \|f-f^*\|_{L_2(\Omega)},
\end{align}
where the second inequality is by the triangle inequality, and last inequality is because of the Gagliardo–Nirenberg interpolation inequality. The Sobolev embedding theorem suggests that $\|f\|_{L_2(\Omega)}\leq C_3\|f\|_{H^{m_0-\delta}(\Omega)}$. Therefore, Lemma \ref{LEM:FSTAR} gives us
\begin{align}\label{eq:2}
\|f^*\|_{L_2(\Omega)} \lesssim \|f\|_{L_2(\Omega)} + \|f-f^*\|_{L_2(\Omega)} \lesssim \|f\|_{L_2(\Omega)}\lesssim \|f\|_{H^{m_0-\delta}(\Omega)},
\end{align}
where the second inequality is because $\lambda_m = o(Q(n)^{-m/m_0})$ yields $\|f-f^*\|_{L_2(\Omega)}\rightarrow 0$, which implies 
$\|f-f^*\|_{L_2(\Omega)} \leq (C-1)\|f\|_{L_2(\Omega)}$ for some $C>0$.

By \eqref{eq:fstarboundX}, \eqref{eq:diffofff1overfixX}, and \eqref{eq:2}, we have
\begin{align}\label{eq:ffstarnnormX}
\|f-f^*\|_n \lesssim &h_n^{m_0-\delta}\|f\|_{H^{m_0-\delta}(\Omega)} +h_n^{m_0-\delta}\|f^*\|_{L_2(\Omega)}^{\frac{m-m_0+\delta}{m}}\|f^*\|_{\mathcal{N}_{\Phi}(\Omega)}^{\frac{m_0-\delta}{m}} + \|f-f^*\|_{L_2(\Omega)}\nonumber\\
\lesssim & h_n^{m_0-\delta}\|f\|_{H^{m_0-\delta}(\Omega)} +h_n^{m_0-\delta}(\|f\|_{H^{m_0-\delta}(\Omega)})^{\frac{m-m_0+\delta}{m}}(\lambda_m^{\frac{m_0-m}{2m}}Q(n)^{1/2})^{\frac{m_0-\delta}{m}} + \lambda_m^{\frac{m_0}{2m}}Q(n)^{1/2}\nonumber\\
\lesssim & (h_n^{m_0-\delta}\lambda_m^{\frac{(m_0-m)(m_0-\delta)}{2m^2}} + \lambda_m^{\frac{m_0}{2m}})\max\{Q(n)^{1/2}, \|f\|_{H^{m_0-\delta}(\Omega)}\},
\end{align}
where the last inequality is because $\lambda_m = o(Q(n)^{-m/m_0})$ and $m_0\leq m$. Since $Q$ is an increasing function and satisfies
\begin{align*}
\lim_{r\rightarrow +\infty} \frac{\log Q(r)}{\log r} = 0,
\end{align*}
there exists a constant $C_4$ such that $Q(r)\leq C_4(1 + r^{2})^\delta$ for all $r\geq 0$. Therefore, by the extension theorem, there exists an extension of $f_e$ such that $f=f_e|_\Omega$, and 
\begin{align*}
& \|f\|_{H^{m_0-\delta}(\Omega)}^2 \lesssim \|f_e\|_{H^{m_0-\delta}(\RR^d)}^2 = \int_{\RR^d} (1+\|\omega\|_2^2)^{m_0-\delta} |\mathcal{F}(f_e)(\omega)|^2d\omega\\
\leq & C_5\int_{\RR^d} \frac{(1+\|\omega\|_2^2)^{m_0} |\mathcal{F}(f_e)(\omega)|^2}{Q(\|\omega\|)}d\omega \leq C_5,
\end{align*}
which implies $\max\{Q(n)^{1/2}, \|f\|_{H^{m_0-\delta}(\Omega)}\}\lesssim Q(n)^{1/2}$. Noting that $h_n\lesssim n^{-1/d}$ because $\mathcal{X}$ satisfies Condition (C4), we have $h_n^{m_0-\delta}\lambda_m^{\frac{(m_0-m)(m_0-\delta)}{2m^2}} \lesssim \lambda_m^{\frac{m_0}{2m}}$. Therefore, \eqref{eq:ffstarnnormX} can be further bounded by
\begin{align}\label{eq:ffstarnnormX1}
\|f-f^*\|_n \lesssim \lambda_m^{\frac{m_0}{2m}}Q(n)^{1/2}.
\end{align}

Repeating the proof in the case of $f\in \mathcal{N}_{\Psi}(\Omega)$, we can obtain that 
\begin{align*}
\|f-\hat f_m\|_n^2 \leq &  C_6T, \mbox{ and }\|\hat f_m\|_{\mathcal{N}_{\Phi}(\Omega)}^2 \leq C_6\lambda_m^{-1}T,
\end{align*}
where
\begin{align*}
	T = &\max\{t^{\frac{4m}{2m+d}}n^{-\frac{2m}{2m+d}}\lambda_m^{\frac{d(m_0-m)}{m(2m+d)}}Q(n)^{\frac{d}{2m+d}}, \lambda_m^{\frac{m_0}{m}}Q(n) + tn^{-\frac{1}{2}}\lambda_m^{\frac{2m_0-d}{4m}}Q(n)^{1/2},\nonumber\\
	& t^{\frac{4m}{4m-d}}n^{-\frac{2m}{4m-d}} \lambda_m^{\frac{4mm_0-2m_0d-2md}{2m(4m-d)}}Q(n)^{\frac{2m-d}{4m-d}}, t^2n^{-1}\lambda_m^{-\frac{d}{2m}}\}.
\end{align*}
Similar to the proof of \eqref{eq:pfpvercaseL2fix}, we have
\begin{align*}
\|f-\hat f_m\|_{L_2(\Omega)}\leq & \|f-f^*\|_{L_2(\Omega)} + \|f^*-\hat f_m\|_{L_2(\Omega)}\nonumber\\
\lesssim & \lambda_m^{\frac{m_0}{2m}}Q(n)^{1/2} + h_n^{m} \|f^*-\hat f_m\|_{H^{m}(\Omega)} + \|f^*-\hat f_m\|_n\nonumber\\
\lesssim & \lambda_m^{\frac{m_0}{2m}}Q(n)^{1/2} +h_n^{m}\lambda_m^{\frac{m_0-m}{2m}}Q(n)^{1/2} + h_n^{m}\lambda_m^{-1}T + T)\nonumber\\
\lesssim & T,
\end{align*}
where the second inequality is by Lemma \ref{LEM:FSTAR}. This finishes the proof of the case $f\notin \mathcal{N}_{\Psi}(\Omega)$, thus finishes the proof of Theorem \ref{thm:krroverg}.

\section{Proof of Theorem \ref{thm:krrunder}}\label{subsec:pfkrrunder}

In this section, we set $m_0:=m_0(f)$ for notational simplicity. We show a generalized version of Theorem \ref{thm:krrunder} as follows. Recall that  $H^{m_0}(\Omega)$ coincides with $\mathcal{N}_\Psi(\Omega)$.

\begin{theorem}\label{thm:krrunderg}
Suppose conditions in Theorem \ref{thm:krrunder} hold.
	Suppose $\lambda_m = o(1)$ if $f\in \mathcal{N}_{\Psi}$, and $\lambda_m = o(Q(n)^{-2m/m_0})$ if $f\notin \mathcal{N}_{\Psi}$.
Furthermore, suppose $\lambda_m \gtrsim n^{-\frac{2m}{d}}$. If $f\in \mathcal{N}_{\Psi}(\Omega)$, let 
	\begin{align*}
	T =\max\{ & n^{-\frac{2m_0}{2m_0+d}}\|f\|_{\mathcal{N}_{\Psi}(\Omega)}^2 + \lambda_{m}n^{-\frac{2(m_0-m)}{2m_0+d}}\|f\|_{\mathcal{N}_{\Psi}(\Omega)}^2, \nonumber\\
    & t^{\frac{4m}{4m-d}}n^{-{\frac{8mm_0+2md-2m_0d}{(4m-d)(2m_0+d)}}}\lambda_m^{-\frac{d}{4m-d}} \|f\|_{\mathcal{N}_{\Psi}(\Omega)}^{\frac{2m-d}{4m-d}}, t^2n^{-1}\lambda_m^{-\frac{d}{2m}} \}.
	\end{align*}
	For all $t\geq C_0$ and $n$, with probability at least $1 - C_1\exp(-C_2t^2)$,
	\begin{align*}
	\| f - \hat f_m\|_{L_2(\Omega)}\lesssim & T^{1/2}.
	\end{align*}
	If $f\notin \mathcal{N}_{\Psi}(\Omega)$, for all $t\geq C_0$ and $n$, with probability at least $1 - C_1\exp(-C_2t^2)$,
	\begin{align*}
	\| f - \hat f_m\|_{L_2(\Omega)} \lesssim T_1^{1/2},
	\end{align*}
where
\begin{align*}
    T_1 =  \max\{ & n^{-\frac{2m_0}{2m_0+d}}Q(n) + \lambda_{m}n^{-\frac{2(m_0-m)}{2m_0+d}}Q(n), \nonumber\\
    & t^{\frac{4m}{4m-d}}n^{-{\frac{8mm_0+2md-2m_0d}{(4m-d)(2m_0+d)}}}\lambda_m^{-\frac{d}{4m-d}} Q(n)^{\frac{2m-d}{8m-2d}}, t^2n^{-1}\lambda_m^{-\frac{d}{2m}} \}.
\end{align*}
\end{theorem}

We first show that Theorem \ref{thm:krrunderg} implies Theorem \ref{thm:krrunder}. The results of Theorem \ref{thm:krrunder} under the case of $m_0/2\leq m < m_0$ can be obtained by setting $\lambda_m\asymp n^{-\frac{2m}{2m_0+d}}$.

If $d/2<m<m_0/2$, we have that $f\in H^{2m}(\Omega)$. Therefore, replacing $m_0$ by $2m$ and setting $\lambda_m\asymp n^{-\frac{2m}{4m+d}}$ in the case of $m_0/2\leq m < m_0$, we obtain 
\begin{align*}
\| f - \hat f_m\|_{L_2(\Omega)} \lesssim n^{-\frac{2m}{4m + d}}.
\end{align*}
Therefore, we conclude that Theorem \ref{thm:krrunderg} implies Theorem \ref{thm:krrunder}.

Now we begin to prove Theorem \ref{thm:krroverg}. We need the following lemmas.
Lemma \ref{lem:handhe} is Proposition 2.1 of \cite{tuo2015theoretical}. 
\begin{lemma}\label{lem:handhe}
	Each $h\in \mathcal{N}_{\Phi}(\Omega)$ has an extension $h\in \mathcal{N}_{\Phi}(\RR^d)$ which defines an isometric map from $\mathcal{N}_{\Phi}(\Omega)$ to $\mathcal{N}_{\Phi}(\RR^d)$. In other words,$h_e|_\Omega\in \mathcal{N}_{\Phi}(\Omega)$, and $\langle h_e,h_e'\rangle_{\mathcal{N}_{\Phi}(\RR^d)} = \langle h,h'\rangle_{\mathcal{N}_{\Phi}(\Omega)}$ for all $h,h'\in \mathcal{N}_{\Phi}(\Omega)$, where $h_e|_\Omega$ denotes the restriction of $h_e$ on the region $\Omega$.
\end{lemma}
\begin{remark}\label{Reextmap}
    As shown in \cite{tuo2020improved}, the map is extended by the map from $F_\Phi(\Omega)$ defined in \eqref{FPhi} to $F_\Phi(\RR^d)$ given by
\begin{align*}
    \sum_{k=1}^n \beta_k\Phi(x -x_k), x\in \Omega \mapsto \sum_{k=1}^n \beta_k\Phi(x -x_k), x\in \RR^d.
\end{align*}
\end{remark}

Lemma \ref{lem:Thasv} is implied by the proof of Theorem 2.2 of \cite{tuo2020improved}.
\begin{lemma}\label{lem:Thasv}
	Let $\Phi$ satisfy Condition (C3) and $f\in \mathcal{N}_{\Phi}(\Omega)$. Let $f_e$ be an extended function by the map in Lemma \ref{lem:handhe}. Suppose $f_e\in H^{2m}(\RR^d)$, then the integral equation
	\begin{align*}
	f(x) = \int_\Omega \Phi(x-y)v(y)dy,\quad x\in \Omega,
	\end{align*}
	has a solution $v = h_f|_\Omega\in L_2(\Omega)$, where $h_f= \mathcal{F}^{-1}(\mathcal{F}(f_e)/\mathcal{F}(\Phi))$.
\end{lemma}

The proof has three steps. In Step 1, we establish an \textit{improved basic inequality}. In Step 2, we prove the results under the scenario that $f\in \mathcal{N}_\Psi(\Omega)$. In Step 3, we prove the results under the scenario that $f\notin \mathcal{N}_\Psi(\Omega)$.

\noindent\textbf{Step 1: Establish the improved basic inequality.}

Let 
\begin{align*}
\Psi_{2m}(x)=\frac{1}{\Gamma(2m-d/2)2^{2m-d/2 -1}}\| x\|_2^{2m-d/2} K_{2m-d/2}(\| x\|_2),
\end{align*}
i.e., $\Psi_{2m}$ be the Mat\'ern kernel function as in \eqref{materngai} with $\nu = 2m-d/2$ and $\phi = (2\sqrt{2m-d/2})^{-1}$. Therefore, \eqref{maternspec} implies that there exist constants $C_2\geq C_1 > 0$ such that
\begin{align*}
C_1(1+\|\omega\|_2^2)^{-2m} \leq  \mathcal{F}(\Psi_{2m})(\omega)\leq C_2(1+\|\omega\|_2^2)^{-2m}, 
\end{align*}
and $\mathcal{N}_{\Psi_{2m}}(\Omega)$ coincides with the Sobolev space $H^{2m}(\Omega)$. 

Let $f^*$ be the solution to the optimization problem 
\begin{align}\label{eq:underoptOmega}
\min_{\tilde f \in \mathcal{N}_{\Psi_{2m}}(\Omega)}\|f-\tilde f\|_{L_2(\Omega)}^2 + \lambda_{2m} \|\tilde f\|^2_{\mathcal{N}_{\Psi_{2m}}(\Omega)},
\end{align}
and $f^*_n$ be the solution to the optimization problem 
\begin{align}\label{eq:underoptOmegan}
\min_{\tilde f \in \mathcal{N}_{\Psi_{2m}}(\Omega)}\|f-\tilde f\|_{n}^2 + \lambda_{2m} \|\tilde f\|^2_{\mathcal{N}_{\Psi_{2m}}(\Omega)},
\end{align}
where $\lambda_{2m}= n^{-\frac{4m}{2m_0+d}}$ is a regularization parameter.

Because $\hat f_m$ is the solution to \eqref{KRRest1}, we have
\begin{align}\label{eq:basicineq1under}
\|y-\hat f_m\|_n^2 + \lambda_m\|\hat f\|_{\mathcal{N}_{\Phi}(\Omega)}^2 \leq \|y-f_n^*\|_n^2 + \lambda_m\|f_n^*\|_{\mathcal{N}_{\Phi}(\Omega)}^2.
\end{align}
By rearrangement, \eqref{eq:basicineq1under} yields
\begin{align}\label{eq:basicineq2under}
\|f-\hat f_m\|_n^2 + \lambda_m\|\hat f_m\|_{\mathcal{N}_{\Phi}(\Omega)}^2 \leq \|f-f_n^*\|_n^2 + \lambda_m\|f_n^*\|_{\mathcal{N}_{\Phi}(\Omega)}^2 + 2\langle \epsilon,\hat f_m - f_n^*\rangle_n.
\end{align}
Notice that 
\begin{align}\label{eq:idunder}
\|f_n^*\|_{\mathcal{N}_{\Phi}(\Omega)}^2 - \|\hat f_m\|_{\mathcal{N}_{\Phi}(\Omega)}^2 = 2\langle f_n^* - \hat f_m, f_n^*\rangle_{\mathcal{N}_{\Phi}(\Omega)} - \|f_n^* - \hat f_m\|_{\mathcal{N}_{\Phi}(\Omega)}^2.
\end{align}
Plugging \eqref{eq:idunder} into \eqref{eq:basicineq2under} gives us
\begin{align}\label{eq:basicineq3under}
& \|f-\hat f_m\|_n^2 + \lambda_m\|f_n^* - \hat f_m\|_{\mathcal{N}_{\Phi}(\Omega)}^2\leq \|f-f_n^*\|_n^2 + 2\lambda_m\langle f_n^* - \hat f_m, f_n^*\rangle_{\mathcal{N}_{\Phi}(\Omega)} + 2\langle \epsilon,\hat f_m - f_n^*\rangle_n.
\end{align}
Next, we consider the term $\langle f_n^* - \hat f_m, f_n^*\rangle_{\mathcal{N}_{\Phi}(\Omega)}$ in \eqref{eq:basicineq3under}. By the representer theorem, the solution of \eqref{eq:underoptOmegan} can be expressed by
\begin{align*}
f_n^*(x) = r_{2m}(x)^T (R_{2m}+n\lambda_{2m} I_n)^{-1} F,
\end{align*}
where $r_{2m}(x)=(\Psi_{2m}(x-x_1),\ldots,\Psi_{2m}(x-x_n))^T, R_{2m}=(\Psi_{2m}(x_j-x_k))_{j k}$, $I_n$ is an identity matrix, and $F = (f(x_1),\ldots,f(x_n))^T$. We can extend $f_n^*(x)$ from $\mathcal{N}_{\Psi_{2m}}(\Omega)$ to $\mathcal{N}_{\Psi_{2m}}(\RR^d)$ by the map
\begin{align*}
f_n^*(x) = r_{2m}(x)^T (R_{2m}+n\lambda_{2m} I_n)^{-1} F, x\in \Omega \mapsto f_{n,e}^*(x) = r_{2m}(x)^T (R_{2m}+n\lambda_{2m} I_n)^{-1} F, x\in \RR^d.
\end{align*}
Clearly, $\|f_n^*\|_{\mathcal{N}_{\Psi_{2m}}(\Omega)} = \|f_{n,e}^*\|_{\mathcal{N}_{\Psi_{2m}}(\RR^d)}$. As explained in Remark \ref{Reextmap}, $f_{n,e}$ is the extension as in Lemma \ref{lem:handhe}. By the equivalence of $\mathcal{N}_{\Psi_{2m}}(\RR^d)$ and $H^{2m}(\RR^d)$, we can apply Lemma \ref{lem:Thasv} to $f_n^*$, and obtain
\begin{align*}
f_n^*(x) = \int_\Omega \Phi(x-y)v(y)dy,\quad x\in \Omega,
\end{align*}
with $v = h_f|_\Omega\in L_2(\Omega)$, where $h_f= \mathcal{F}^{-1}(\mathcal{F}(f_{n,e}^*)/\mathcal{F}(\Phi))$. Proposition 10.28 of \cite{wendland2004scattered} shows that for any function $g\in \mathcal{N}_{\Phi}(\Omega)$, $\langle g, f_n^*\rangle_{\mathcal{N}_{\Phi}(\Omega)} = \langle g, v\rangle_{L_2(\Omega)}$. Together with \eqref{eq:basicineq3under}, we have
\begin{align}\label{eq:basicinequ1}
& \|f-\hat f_m\|_n^2 + \lambda_m\|f_n^* - \hat f_m\|_{\mathcal{N}_{\Phi}(\Omega)}^2\leq \|f-f_n^*\|_n^2 + 2\lambda_m\langle f_n^* - \hat f_m, v\rangle_{L_2(\Omega)} + 2\langle \epsilon,\hat f_m - f_n^*\rangle_n\nonumber\\
\leq & \|f-f_n^*\|_n^2 + 2\lambda_m\| f_n^* - \hat f_m\|_{L_2(\Omega)} \|v\|_{L_2(\Omega)} + 2\langle \epsilon,\hat f_m - f_n^*\rangle_n.
\end{align}
We call \eqref{eq:basicinequ1} the \textit{improved basic inequality}, which improves the basic inequality \eqref{eq:basicineq2} in the oversmoothed case.

\noindent\textbf{Step 2: Prove the results under the case $f\in \mathcal{N}_\Psi(\Omega)$.}

Applying Lemma \ref{lemmaefsmall} to $\langle \epsilon,\hat f_m - f^*\rangle_n$, we obtain that for all $t\geq C_0$, with probability at least $1 - C_1\exp(-C_2 t^2)$, 
\begin{align}\label{eq:efdsmall2}
\langle \epsilon,\hat f_m - f_n^*\rangle_n \leq tn^{-\frac{1}{2}} \|\hat f_m - f_n^*\|_n^{1-\frac{d}{2m}}\|\hat f_m - f_n^*\|_{\mathcal{N}_{\Phi}(\Omega)}^{\frac{d}{2m}}.
\end{align}
Plugging \eqref{eq:efdsmall2} into \eqref{eq:basicinequ1}, we obtain
\begin{align}\label{eq:impbasicunder}
& \|f-\hat f_m\|_n^2 + \lambda_m\|f_n^* - \hat f_m\|_{\mathcal{N}_{\Phi}(\Omega)}^2\nonumber\\
\leq & \|f-f_n^*\|_n^2 + 2\lambda_m\| f_n^* - \hat f_m\|_{L_2(\Omega)} \|v\|_{L_2(\Omega)} + Ctn^{-\frac{1}{2}} \|\hat f_m - f_n^*\|_n^{1-\frac{d}{2m}}\|\hat f_m - f_n^*\|_{\mathcal{N}_{\Phi}(\Omega)}^{\frac{d}{2m}}\nonumber\\
\leq & \|f-f_n^*\|_n^2 + 2\lambda_m\| f - f_n^*\|_{L_2(\Omega)} \|v\|_{L_2(\Omega)} + 2\lambda_m\| f - \hat f_m\|_{L_2(\Omega)} \|v\|_{L_2(\Omega)}\nonumber\\
& + Ctn^{-\frac{1}{2}} \|f - f_n^*\|_n^{1-\frac{d}{2m}}\|\hat f_m - f_n^*\|_{\mathcal{N}_{\Phi}(\Omega)}^{\frac{d}{2m}} + Ctn^{-\frac{1}{2}} \|f - \hat f_m\|_n^{1-\frac{d}{2m}}\|\hat f_m - f_n^*\|_{\mathcal{N}_{\Phi}(\Omega)}^{\frac{d}{2m}},
\end{align}
where the second inequality is because of the triangle inequality $\|\hat f_m - f_n^*\|_n\leq \|f - \hat f_m\|_n +\|f - f_n^*\|_n$ and the basic inequality $(a+b)^q \leq a^q + b^q$ for any $a,b\geq 0$ and $q\in [0,1]$.

Since $m_0/2\leq m$, we have $2m \geq m_0$. Recall that $f^*$ is the solution to \eqref{eq:underoptOmega}. Lemma \ref{LEM:FSTAR} gives us that 
\begin{align}\label{eq:unfstarbound}
\|f-f^*\|_{L_2}^2 \lesssim \lambda_{2m}^{\frac{m_0}{2m}}\|f\|_{\mathcal{N}_{\Psi}(\Omega)}^2, \mbox{ and } \|f^*\|_{\mathcal{N}_{\Psi_{2m}}(\Omega)}^2 \lesssim \lambda_{2m}^{\frac{m_0-2m}{2m}}\|f\|_{\mathcal{N}_{\Psi}(\Omega)}^2. 
\end{align}
As shown in \eqref{eq:ffstarnnorm}, we have
\begin{align}\label{eq:3}
&  \|f-f^*\|_n^2 \lesssim (h_n^{m_0}\lambda_{2m}^{\frac{(m_0-2m)m_0}{8m^2}}+ \lambda_{2m}^{\frac{m_0}{4m}})^2\|f\|_{\mathcal{N}_{\Psi}(\Omega)}^2\nonumber\\ 
\lesssim & (h_n^{2m_0}\lambda_{2m}^{\frac{(m_0-2m)m_0}{4m^2}}+ \lambda_{2m}^{\frac{m_0}{2m}})\|f\|_{\mathcal{N}_{\Psi}(\Omega)}^2
\lesssim \lambda_{2m}^{\frac{m_0}{2m}}\|f\|_{\mathcal{N}_{\Psi}(\Omega)}^2,
\end{align}
where the second inequality is because of the Cauchy-Schwarz inequality, and the last inequality is because $h_n \lesssim n^{-1/d}$ and $\lambda_{2m} = n^{-\frac{4m}{2m_0+d}}$ (thus $h_n^{2m_0}\lambda_{2m}^{\frac{(m_0-2m)m_0}{4m^2}}\lesssim \lambda_{2m}^{\frac{m_0}{2m}}$).

Because $f_n^*$ is the solution to \eqref{eq:underoptOmegan}, we have
\begin{align}\label{eq:underffstar}
\|f- f^*_n\|_{n}^2 + \lambda_{2m} \|f^*_n\|^2_{\mathcal{N}_{\Psi_{2m}}(\Omega)} \leq & \|f- f^*\|_{n}^2 + \lambda_{2m} \|f^*\|^2_{\mathcal{N}_{\Psi_{2m}}(\Omega)}.
\end{align}

Combining \eqref{eq:unfstarbound}, \eqref{eq:3}, and \eqref{eq:underffstar}  yields
\begin{align*}
\|f- f^*_n\|_{n}^2 + \lambda_{2m} \|f^*_n\|^2_{\mathcal{N}_{\Psi_{2m}}(\Omega)} \lesssim \lambda_{2m}^{\frac{m_0}{2m}}\|f\|_{\mathcal{N}_{\Psi}(\Omega)}^2,
\end{align*}
which implies
\begin{align}\label{eq:underffstar1}
& \|f- f^*_n\|_{n}^2\lesssim \lambda_{2m}^{\frac{m_0}{2m}}\|f\|_{\mathcal{N}_{\Psi}(\Omega)}^2,
\mbox{ and } \|f^*_n\|^2_{\mathcal{N}_{\Psi_{2m}}(\Omega)} \lesssim \lambda_{2m}^{\frac{m_0-2m}{2m}}\|f\|_{\mathcal{N}_{\Psi}(\Omega)}^2.
\end{align}
Similar to \eqref{eq:GNinter}, by the Gagliardo–Nirenberg interpolation inequality, we can show that
\begin{align}\label{eq:GNinterunder}
\|f^*_n\|_{\mathcal{N}_{\Psi}(\Omega)} \lesssim \|f^*_n\|_{H^{m_0}(\Omega)} \lesssim  \|f^*_n\|_{L_2(\Omega)}^{\frac{2m-m_0}{2m}}\|f^*_n\|_{H^{2m}(\Omega)}^{\frac{m_0}{2m}} \lesssim \|f^*_n\|_{L_2(\Omega)}^{\frac{2m-m_0}{2m}}\|f^*_n\|_{\mathcal{N}_{\Psi_{2m}}(\Omega)}^{\frac{m_0}{2m}}.
\end{align}
Applying Lemma \ref{lemmafixdesign} to $\|f- f^*_n\|_{L_2(\Omega)}$, it can be seen that
\begin{align}\label{eq:ffstarnb2}
\|f- f^*_n\|_{L_2(\Omega)} \lesssim & h_n^{m_0}\|f- f^*_n\|_{H^{m_0}(\Omega)} + \|f- f^*_n\|_n\nonumber\\
\lesssim & h_n^{m_0}\|f- f^*_n\|_{N_{\Psi}(\Omega)} + \|f- f^*_n\|_n\nonumber\\
\lesssim & h_n^{m_0}\|f\|_{N_{\Psi}(\Omega)} + h_n^{m_0}\|f^*_n\|_{N_{\Psi}(\Omega)} + \|f- f^*_n\|_n\nonumber\\
\lesssim & h_n^{m_0}\|f\|_{N_{\Psi}(\Omega)} + h_n^{m_0}\|f^*_n\|_{L_2(\Omega)}^{\frac{2m-m_0}{2m}}\|f^*_n\|_{\mathcal{N}_{\Psi_{2m}}(\Omega)}^{\frac{m_0}{2m}} + \lambda_{2m}^{\frac{m_0}{4m}}\|f\|_{\mathcal{N}_{\Psi}(\Omega)}\nonumber\\
\lesssim &  h_n^{m_0}\|f^*_n\|_{L_2(\Omega)}^{\frac{2m-m_0}{2m}}\|f^*_n\|_{\mathcal{N}_{\Psi_{2m}}(\Omega)}^{\frac{m_0}{2m}} + \lambda_{2m}^{\frac{m_0}{4m}}\|f\|_{\mathcal{N}_{\Psi}(\Omega)},
\end{align}
where the second inequality is because $\|\cdot\|_{H^{m_0}(\Omega)}$ is equivalent to $\|\cdot\|_{N_{\Psi}(\Omega)}$, the third inequality is because of the triangle inequality, the fourth inequality is because of \eqref{eq:underffstar1} and \eqref{eq:GNinterunder}, and the last inequality is because $h_n^{m_0}\lesssim \lambda_{2m}^{\frac{m_0}{4m}}$.

Applying Lemma \ref{lemmafixdesign} to $\|f^*_n\|_{L_2(\Omega)}$ leads to
\begin{align}\label{eq:fstarnBun}
\|f^*_n\|_{L_2(\Omega)} \lesssim & h_n^{2m}\| f^*_n\|_{H^{2m}(\Omega)} + \|f^*_n\|_n
\lesssim  h_n^{2m}\| f^*_n\|_{N_{\Psi_{2m}}(\Omega)} + \|f^*_n\|_n\nonumber\\
\lesssim & h_n^{2m}\|f^*_n\|_{N_{\Psi_{2m}}(\Omega)} + \|f\|_n + \|f-f^*\|_n\nonumber\\
\lesssim & h_n^{2m} \lambda_{2m}^{\frac{m_0}{4m}-1/2}\|f\|_{\mathcal{N}_{\Psi}(\Omega)} + \|f\|_n + \lambda_{2m}^{\frac{m_0}{4m}}\|f\|_{\mathcal{N}_{\Psi}(\Omega)} 
\lesssim \|f\|_{\mathcal{N}_{\Psi}(\Omega)},
\end{align}
where the second inequality is because $\|\cdot\|_{H^{2m}(\Omega)}$ is equivalent to $\|\cdot\|_{N_{\Psi_{2m}}(\Omega)}$, the third inequality is because of the triangle inequality, the fourth inequality is because of \eqref{eq:3} and \eqref{eq:underffstar1}, and the last inequality is because of $h_n\lesssim n^{-1/d}$ and \eqref{eq:reproP}.

Plugging \eqref{eq:underffstar1} and \eqref{eq:fstarnBun} into \eqref{eq:ffstarnb2}, we have
\begin{align}\label{eq:ffstarBNunder2}
\|f- f^*_n\|_{L_2(\Omega)}
\lesssim & h_n^{m_0}\lambda_{2m}^{\frac{m_0(m_0-2m)}{8m^2}}\|f\|_{\mathcal{N}_{\Psi}(\Omega)} + \lambda_{2m}^{\frac{m_0}{4m}}\|f\|_{\mathcal{N}_{\Psi}(\Omega)}\lesssim \lambda_{2m}^{\frac{m_0}{4m}}\|f\|_{\mathcal{N}_{\Psi}(\Omega)},
\end{align}
where the last inequality is because $h_n^{m_0}\lambda_{2m}^{\frac{m_0(m_0-2m)}{8m^2}}\lesssim \lambda_{2m}^{\frac{m_0}{4m}}$, which can be checked by noting the fact $h_n\lesssim n^{-1/d}$ and $\lambda_{2m} = n^{-\frac{4m}{2m_0+d}}$.

The Plancherel theorem \citep{bracewell1986fourier} implies that
\begin{align}\label{eq:vboundu}
& \|v\|_{L_2(\Omega)}^2\lesssim \|h_f\|_{L_2(\RR^d)}^2 = \|\mathcal{F}(f_{n,e}^*)/\mathcal{F}(\Phi)\|_{L_2(\RR^d)}^2 \lesssim \|f_{n,e}^*\|_{\mathcal{N}_{\Psi_{2m}}(\RR^d)}^2 \lesssim \lambda_{2m}^{\frac{m_0-2m}{2m}}\|f\|_{\mathcal{N}_{\Psi}(\Omega)}^2.
\end{align}
So far, we have provided upper bounds of $\|f-f_n^*\|_n$, $\| f - f_n^*\|_{L_2(\Omega)}$, and $\|v\|_{L_2(\Omega)}$. It remains to solve \eqref{eq:impbasicunder}, which can be divided to several cases. Note that \eqref{eq:impbasicunder} implies that either
\begin{align}\label{eq:underc1}
&  \|f-\hat f_m\|_n^2 + \lambda_m\|f_n^* - \hat f_m\|_{\mathcal{N}_{\Phi}(\Omega)}^2
\leq  2(\|f-f_n^*\|_n^2 + 2\lambda_m\| f - f_n^*\|_{L_2(\Omega)} \|v\|_{L_2(\Omega)}),
\end{align}
or 
\begin{align}\label{eq:underc2}
&  \|f-\hat f_m\|_n^2 + \lambda_m\|f_n^* - \hat f_m\|_{\mathcal{N}_{\Phi}(\Omega)}^2\nonumber\\
\leq & 2(2\lambda_m\| f - \hat f_m\|_{L_2(\Omega)} \|v\|_{L_2(\Omega)}+ Ctn^{-\frac{1}{2}} \|f - f_n^*\|_n^{1-\frac{d}{2m}}\|\hat f_m - f_n^*\|_{\mathcal{N}_{\Phi}(\Omega)}^{\frac{d}{2m}}\nonumber\\
&  + Ctn^{-\frac{1}{2}} \|f - \hat f_m\|_n^{1-\frac{d}{2m}}\|f_n^* - \hat f_m\|_{\mathcal{N}_{\Phi}(\Omega)}^{\frac{d}{2m}}).
\end{align}
Plugging \eqref{eq:underffstar1}, \eqref{eq:ffstarBNunder2}, and \eqref{eq:vboundu} into \eqref{eq:underc1} leads to
\begin{align*}
&  \|f-\hat f_m\|_n^2 + \lambda_m\|f_n^* - \hat f_m\|_{\mathcal{N}_{\Phi}(\Omega)}^2\lesssim \lambda_{2m}^{\frac{m_0}{2m}}\|f\|_{\mathcal{N}_{\Psi}(\Omega)}^2 + \lambda_{m}\lambda_{2m}^{\frac{m_0-m}{2m}}\|f\|_{\mathcal{N}_{\Psi}(\Omega)}^2,
\end{align*}
which implies
\begin{align}\label{eq:underr1}
&  \|f-\hat f_m\|_n^2\lesssim \lambda_{2m}^{\frac{m_0}{2m}}\|f\|_{\mathcal{N}_{\Psi}(\Omega)}^2 + \lambda_{m}\lambda_{2m}^{\frac{m_0-m}{2m}}\|f\|_{\mathcal{N}_{\Psi}(\Omega)}^2,\nonumber\\
\mbox{and } &  \|f_n^* - \hat f_m\|_{\mathcal{N}_{\Phi}(\Omega)}^2 \lesssim \lambda_m^{-1}\lambda_{2m}^{\frac{m_0}{2m}}\|f\|_{\mathcal{N}_{\Psi}(\Omega)}^2 + \lambda_{2m}^{\frac{m_0-m}{2m}}\|f\|_{\mathcal{N}_{\Psi}(\Omega)}^2.
\end{align}
Solving \eqref{eq:underc2} is more complicated. By Lemma \ref{lemmafixdesign}, it holds that
\begin{align}\label{eq:u1}
\| f - \hat f_m\|_{L_2(\Omega)} \lesssim & h_n^m\| f - \hat f_m\|_{H^m(\Omega)} + \| f - \hat f_m\|_{n}\nonumber\\
\lesssim & h_n^m\| f - \hat f_m\|_{\mathcal{N}_{\Phi}(\Omega)} + \| f - \hat f_m\|_{n}\nonumber\\
\lesssim & h_n^m\| f_n^* - \hat f_m\|_{\mathcal{N}_{\Phi}(\Omega)} +h_n^m\|f -  f_n^*\|_{\mathcal{N}_{\Phi}(\Omega)} + \| f - \hat f_m\|_{n},
\end{align}
where the second inequality is by the equivalence of $\|\cdot\|_{H^m(\Omega)}$ and $\| \cdot\|_{\mathcal{N}_{\Phi}(\Omega)}$, and the third inequality is by the triangle inequality. 

By the Gagliardo–Nirenberg interpolation inequality, we have
\begin{align}\label{eq:GNinteru1}
& \|f -  f_n^*\|_{\mathcal{N}_{\Phi}(\Omega)} \lesssim \|f -  f_n^*\|_{H^{m}(\Omega)}\lesssim \|f -  f_n^*\|_{L_2(\Omega)}^{\frac{m_0-m}{m_0}}\|f -  f_n^*\|_{H^{m_0}(\Omega)}^{\frac{m}{m_0}},
\end{align}
where the first inequality is because of the equivalence of $\|\cdot\|_{H^m(\Omega)}$ and $\| \cdot\|_{\mathcal{N}_{\Phi}(\Omega)}$. Using the Gagliardo–Nirenberg interpolation inequality again, we find that
\begin{align}\label{eq:GNinteru2}
& \|f -  f_n^*\|_{H^{m_0}(\Omega)} \leq \|f\|_{H^{m_0}(\Omega)} +  \|f_n^*\|_{H^{m_0}(\Omega)} \lesssim \|f\|_{H^{m_0}(\Omega)} +  \|f_n^*\|_{L_2(\Omega)}^{\frac{2m - m_0}{2m}}\|f_n^*\|_{H^{2m}(\Omega)}^{\frac{m_0}{2m}}\nonumber\\
\lesssim & \|f\|_{\mathcal{N}_{\Psi}(\Omega)} +  \|f_n^*\|_{L_2(\Omega)}^{\frac{2m - m_0}{2m}}\|f_n^*\|_{\mathcal{N}_{\Psi_{2m}}(\Omega)}^{\frac{m_0}{2m}}
\lesssim  \|f\|_{\mathcal{N}_{\Psi}(\Omega)} +  \lambda_{2m}^{\frac{(m_0-2m)m_0}{8m^2}}\|f\|_{\mathcal{N}_{\Psi}(\Omega)}\lesssim  \lambda_{2m}^{\frac{(m_0-2m)m_0}{8m^2}}\|f\|_{\mathcal{N}_{\Psi}(\Omega)},
\end{align}
where the first inequality is by the triangle inequality, the third inequality is by the equivalence of $\|\cdot\|_{H^{m_0}(\Omega)}$ and $\| \cdot\|_{\mathcal{N}_{\Psi}(\Omega)}$ and the equivalence of $\|\cdot\|_{H^{2m}(\Omega)}$ and $\| \cdot\|_{\mathcal{N}_{\Psi_{2m}}(\Omega)}$, the fourth inequality is by \eqref{eq:underffstar1} and \eqref{eq:fstarnBun}, and the last inequality is because $\lambda_{2m}\lesssim 1$ and $m_0\leq 2m$. Combining \eqref{eq:ffstarBNunder2}, \eqref{eq:u1}, \eqref{eq:GNinteru1}, and \eqref{eq:GNinteru2} leads to
\begin{align}\label{eq:5}
    \| f - \hat f_m\|_{L_2(\Omega)}
\lesssim & h_n^m\| f_n^* - \hat f_m\|_{\mathcal{N}_{\Phi}(\Omega)} +h_n^m\|f -  f_n^*\|_{L_2(\Omega)}^{\frac{m_0-m}{m_0}}\|f -  f_n^*\|_{H^{m_0}(\Omega)}^{\frac{m}{m_0}} + \| f - \hat f_m\|_{n}\nonumber\\
\lesssim & h_n^m\| f_n^* - \hat f_m\|_{\mathcal{N}_{\Phi}(\Omega)} +h_n^m\lambda_{2m}^{\frac{3m_0-4m}{8m}}\|f\|_{\mathcal{N}_{\Psi}(\Omega)} + \| f - \hat f_m\|_{n}.
\end{align}
Plugging \eqref{eq:5} into \eqref{eq:underc2} leads to
\begin{align}\label{eq:undercc2}
&  \|f-\hat f_m\|_n^2 + \lambda_m\|f_n^* - \hat f_m\|_{\mathcal{N}_{\Phi}(\Omega)}^2\nonumber\\
\lesssim & h_n^m\lambda_m\|v\|_{L_2(\Omega)} \| f_n^* - \hat f_m\|_{\mathcal{N}_{\Phi}(\Omega)} +h_n^m\lambda_m\|v\|_{L_2(\Omega)} \lambda_{2m}^{\frac{3m_0-4m}{8m}}\|f\|_{\mathcal{N}_{\Psi}(\Omega)} + \lambda_m\|v\|_{L_2(\Omega)} \| f - \hat f_m\|_{n}\nonumber\\
& + tn^{-\frac{1}{2}} \|f - f_n^*\|_n^{1-\frac{d}{2m}}\|\hat f_m - f_n^*\|_{\mathcal{N}_{\Phi}(\Omega)}^{\frac{d}{2m}} + tn^{-\frac{1}{2}} \|f - \hat f_m\|_n^{1-\frac{d}{2m}}\|f_n^* - \hat f_m\|_{\mathcal{N}_{\Phi}(\Omega)}^{\frac{d}{2m}}\nonumber\\
= & I_1 + I_2 + I_3 + I_4 + I_5.
\end{align}
Depending on which $I_k$ is equal to $\max\{I_1,I_2,...,I_5\}$, \eqref{eq:undercc2} implies one of the following case is true:

\noindent\textbf{Case 1.1: $I_1=\max\{I_1,I_2,...,I_5\}$.} Under this case, we have
\begin{align}\label{eq:undercc21}
	\|f-\hat f_m\|_n^2 + \lambda_m\|f_n^* - \hat f_m\|_{\mathcal{N}_{\Phi}(\Omega)}^2\lesssim h_n^m\lambda_m\|v\|_{L_2(\Omega)} \| f_n^* - \hat f_m\|_{\mathcal{N}_{\Phi}(\Omega)}.
\end{align}
For the conciseness of this proof, we only provide details on solving \eqref{eq:undercc21} in Case 1.1. The inequalities in other cases can be solved similarly. 

Note that \eqref{eq:undercc21} implies that
\begin{align*}
	\lambda_m\|f_n^* - \hat f_m\|_{\mathcal{N}_{\Phi}(\Omega)}^2 \lesssim h_n^m\lambda_m\|v\|_{L_2(\Omega)} \| f_n^* - \hat f_m\|_{\mathcal{N}_{\Phi}(\Omega)},
\end{align*}
which implies 
\begin{align}\label{eq:undercr111}
	\|f_n^* - \hat f_m\|_{\mathcal{N}_{\Phi}(\Omega)}^2 \lesssim h_n^{2m}\|v\|_{L_2(\Omega)}^2 \lesssim  h_n^{2m}\lambda_{2m}^{\frac{m_0-2m}{2m}}\|f\|_{\mathcal{N}_{\Psi}(\Omega)}^2,
\end{align}
where the last inequality is because of \eqref{eq:vboundu}. Together with \eqref{eq:undercc21}, we have 
\begin{align}\label{eq:undercr112}
	\|f-\hat f_m\|_n^2 \lesssim &  h_n^m\lambda_m\|v\|_{L_2(\Omega)} \| f_n^* - \hat f_m\|_{\mathcal{N}_{\Phi}(\Omega)}
	\lesssim  h_n^{2m}\lambda_m\lambda_{2m}^{\frac{m_0-2m}{2m}}\|f\|_{\mathcal{N}_{\Psi}(\Omega)}^2.
\end{align}
\noindent\textbf{Case 1.2: $I_2=\max\{I_1,I_2,...,I_5\}$.} Under this case, we have
\begin{align}\label{eq:undercc22}
&  \|f-\hat f_m\|_n^2 + \lambda_m\|f_n^* - \hat f_m\|_{\mathcal{N}_{\Phi}(\Omega)}^2
\lesssim h_n^m\lambda_m\|v\|_{L_2(\Omega)} \lambda_{2m}^{\frac{3m_0-4m}{8m}}\|f\|_{\mathcal{N}_{\Psi}(\Omega)}.
\end{align}
By \eqref{eq:vboundu}, \eqref{eq:undercc22} yields
\begin{align}\label{eq:underr22}
	&  \|f-\hat f_m\|_n^2   \lesssim h_n^m\lambda_m\lambda_{2m}^{\frac{5m_0-8m}{8m}}\|f\|_{\mathcal{N}_{\Psi}(\Omega)}^2,\nonumber\\
	&  \|f_n^* - \hat f_m\|_{\mathcal{N}_{\Phi}(\Omega)}^2\lesssim h_n^m\lambda_{2m}^{\frac{5m_0-8m}{8m}}\|f\|_{\mathcal{N}_{\Psi}(\Omega)}^2.
\end{align}

\noindent\textbf{Case 1.3: $I_3=\max\{I_1,I_2,...,I_5\}$.} Under this case, we have
\begin{align}\label{eq:undercc23}
\|f-\hat f_m\|_n^2 + \lambda_m\|f_n^* - \hat f_m\|_{\mathcal{N}_{\Phi}(\Omega)}^2
\lesssim  \lambda_m\|v\|_{L_2(\Omega)} \| f - \hat f_m\|_{n}.
\end{align}
Solving \eqref{eq:undercc23} leads to
\begin{align}\label{eq:underr23}
	&  \|f-\hat f_m\|_n^2   \lesssim \lambda_m^2\lambda_{2m}^{\frac{m_0-2m}{2m}}\|f\|_{\mathcal{N}_{\Psi}(\Omega)}^2,\nonumber\\
	&  \|f_n^* - \hat f_m\|_{\mathcal{N}_{\Phi}(\Omega)}^2\lesssim\lambda_m\lambda_{2m}^{\frac{m_0-2m}{2m}}\|f\|_{\mathcal{N}_{\Psi}(\Omega)}^2,
\end{align}
where we use \eqref{eq:vboundu} to bound $\|v\|_{L_2(\Omega)}$.

\noindent\textbf{Case 1.4: $I_4=\max\{I_1,I_2,...,I_5\}$.} Under this case, we have
\begin{align}\label{eq:undercc24}
\|f-\hat f_m\|_n^2 + \lambda_m\|f_n^* - \hat f_m\|_{\mathcal{N}_{\Phi}(\Omega)}^2
\lesssim tn^{-\frac{1}{2}} \|f - f_n^*\|_n^{1-\frac{d}{2m}}\|\hat f_m - f_n^*\|_{\mathcal{N}_{\Phi}(\Omega)}^{\frac{d}{2m}}.
\end{align}
Solving \eqref{eq:undercc24} leads to
\begin{align}\label{eq:underr24}
		&  \|f-\hat f_m\|_n^2	\lesssim t^{\frac{4m}{4m-d}}n^{-{\frac{2m}{4m-d}}}\lambda_m^{-\frac{d}{4m-d}} (\lambda_{2m}^{\frac{m_0}{2m}}\|f\|_{\mathcal{N}_{\Psi}(\Omega)}^2)^{\frac{2m-d}{4m-d}},\nonumber\\
	&  \|f_n^* - \hat f_m\|_{\mathcal{N}_{\Phi}(\Omega)}^2
	\lesssim t^{\frac{4m}{4m-d}}n^{-\frac{2m}{4m-d}}\lambda_m^{-\frac{4m}{4m-d}} (\lambda_{2m}^{\frac{m_0}{2m}}\|f\|_{\mathcal{N}_{\Psi}(\Omega)}^2)^{\frac{2m-d}{4m-d}},
\end{align}
where we use \eqref{eq:underffstar1} to bound $\|f- f^*_n\|_{n}$.

\noindent\textbf{Case 1.5: $I_5=\max\{I_1,I_2,...,I_5\}$.} Under this case, we have
\begin{align}\label{eq:undercc25}
&  \|f-\hat f_m\|_n^2 + \lambda_m\|f_n^* - \hat f_m\|_{\mathcal{N}_{\Phi}(\Omega)}^2\lesssim tn^{-\frac{1}{2}} \|f - \hat f_m\|_n^{1-\frac{d}{2m}}\|f_n^* - \hat f_m\|_{\mathcal{N}_{\Phi}(\Omega)}^{\frac{d}{2m}}.
\end{align}
Solving \eqref{eq:undercc25} leads to
\begin{align}\label{eq:underr25}
	\|f-\hat f_m\|_n^2 \lesssim t^2n^{-1}\lambda_m^{-\frac{d}{2m}},
	\|f_n^* - \hat f_m\|_{\mathcal{N}_{\Phi}(\Omega)}^2 \lesssim t^2n^{-1}\lambda_m^{-\frac{2m+d}{2m}}.
\end{align}
Combining \eqref{eq:underr1}, \eqref{eq:undercr111}, \eqref{eq:undercr112}, \eqref{eq:underr22}, \eqref{eq:underr23}, \eqref{eq:underr24}, and \eqref{eq:underr25} we have
\begin{align}\label{eq:underall}
\|f-\hat f_m\|_n^2 \lesssim T,
\|f_n^* - \hat f_m\|_{\mathcal{N}_{\Phi}(\Omega)}^2 \lesssim \lambda_m^{-1}T,
\end{align}
where
\begin{align*}
T = \max\{ & \lambda_{2m}^{\frac{m_0}{2m}}\|f\|_{\mathcal{N}_{\Psi}(\Omega)}^2 + \lambda_{m}\lambda_{2m}^{\frac{m_0-m}{2m}}\|f\|_{\mathcal{N}_{\Psi}(\Omega)}^2,  h_n^{2m}\lambda_m\lambda_{2m}^{\frac{m_0-2m}{2m}}\|f\|_{\mathcal{N}_{\Psi}(\Omega)}^2,\nonumber\\
& h_n^m\lambda_m\lambda_{2m}^{\frac{5m_0-8m}{8m}}\|f\|_{\mathcal{N}_{\Psi}(\Omega)}^2, \lambda_m^2\lambda_{2m}^{\frac{m_0-2m}{2m}}\|f\|_{\mathcal{N}_{\Psi}(\Omega)}^2,\nonumber\\
& t^{\frac{4m}{4m-d}}n^{-{\frac{2m}{4m-d}}}\lambda_m^{-\frac{d}{4m-d}} (\lambda_{2m}^{\frac{m_0}{2m}}\|f\|_{\mathcal{N}_{\Psi}(\Omega)}^2)^{\frac{2m-d}{4m-d}}, t^2n^{-1}\lambda_m^{-\frac{d}{2m}} \}\nonumber\\
= \max\{ & \lambda_{2m}^{\frac{m_0}{2m}}\|f\|_{\mathcal{N}_{\Psi}(\Omega)}^2 + \lambda_{m}\lambda_{2m}^{\frac{m_0-m}{2m}}\|f\|_{\mathcal{N}_{\Psi}(\Omega)}^2, \nonumber\\
& t^{\frac{4m}{4m-d}}n^{-{\frac{2m}{4m-d}}}\lambda_m^{-\frac{d}{4m-d}} (\lambda_{2m}^{\frac{m_0}{2m}}\|f\|_{\mathcal{N}_{\Psi}(\Omega)}^2)^{\frac{2m-d}{4m-d}}, t^2n^{-1}\lambda_m^{-\frac{d}{2m}} \},
\end{align*}
because $\lambda_m^2\lambda_{2m}^{\frac{m_0-2m}{2m}} \lesssim \lambda_m\lambda_{2m}^{\frac{m_0-2m}{2m}}$, $h_n \lesssim n^{-1/d}$ and $\lambda_{2m} = n^{-\frac{4m}{2m_0+d}}$.

We finish Step 2 by bounding the difference between $\|f-\hat f_m\|_n^2$ and $\|f-\hat f_m\|_{L_2(\Omega)}^2$. 
By \eqref{eq:ffstarBNunder2}, \eqref{eq:GNinteru1}, and \eqref{eq:GNinteru2}, it can be seen that 
\begin{align}\label{eq:6}
    \|f -  f_n^*\|_{\mathcal{N}_{\Phi}(\Omega)} \lesssim 
\|f -  f_n^*\|_{L_2(\Omega)}^{\frac{m_0-m}{m_0}}\|f -  f_n^*\|_{H^{m_0}(\Omega)}^{\frac{m}{m_0}} \lesssim \lambda_{2m}^{\frac{3m_0-4m}{8m}}\|f\|_{\mathcal{N}_{\Psi}(\Omega)}.
\end{align}
By \eqref{eq:u1}, \eqref{eq:underall}, and \eqref{eq:6}, we have
\begin{align*}
\| f - \hat f_m\|_{L_2(\Omega)} 
\lesssim & h_n^m\| f_n^* - \hat f_m\|_{\mathcal{N}_{\Phi}(\Omega)} +h_n^m\|f -  f_n^*\|_{\mathcal{N}_{\Phi}(\Omega)} + \| f - \hat f_m\|_{n}\nonumber\\
\lesssim & h_n^m\lambda_m^{-1/2}T^{1/2} +h_n^m\lambda_{2m}^{\frac{3m_0-4m}{8m}}\|f\|_{\mathcal{N}_{\Psi}(\Omega)} + T^{1/2}\nonumber\\
\lesssim & T^{1/2} + n^{-\frac{4m_0m + 3m_0d-2md}{2d(2m_0+d)}}\lesssim T^{1/2},
\end{align*}
where the second inequality is by $\lambda_m \gtrsim n^{-\frac{2m}{d}}$, the third inequality is by $h_n \lesssim n^{-1/d}$ and $\lambda_{2m} = n^{-\frac{4m}{2m_0+d}}$, and the last inequality is because the optimal convergence rate of $T$ is $n^{-\frac{2m_0}{2m_0+d}}$, which is always larger than $n^{-\frac{4m_0m + 3m_0d-2md}{d(2m_0+d)}}$. Therefore, we finish the proof of the case $f\in \mathcal{N}_{\Psi(\Omega)}$.

\noindent\textbf{Step 3: Proof of the case $f\notin \mathcal{N}_{\Psi(\Omega)}$.}

Let $f^*$ be as in \eqref{eq:underoptOmega}. We still set $\lambda_{2m} = n^{-\frac{4m}{2m_0+d}}$. If $f\notin \mathcal{N}_{\Psi(\Omega)}$, Lemma \ref{LEM:FSTAR} gives us that 
\begin{align}\label{eq:unfstarboundX}
\|f-f^*\|_{L_2}^2 \leq C\lambda_{2m}^{\frac{m_0}{2m}}Q(n), \mbox{ and } \|f^*\|_{\mathcal{N}_{\Psi_{2m}}(\Omega)}^2 \leq C\lambda_{2m}^{\frac{m_0-2m}{2m}}Q(n),
\end{align}
where $Q: \RR_+\mapsto \RR_+$ satisfies
$$
\lim_{r\rightarrow +\infty} \frac{\log Q(r)}{\log r} = 0.
$$
Let $\delta = \frac{m_0^2(m_0-d/2)}{2(m_0+2m)(m_0-d/2)+m_0^2}>0$. Thus, $m_0-\delta > d/2$. As shown in \eqref{eq:ffstarnnormX1} (note that we replace $m$ in \eqref{eq:ffstarnnormX1} by $2m$ in \eqref{eq:ffstarnnormXXXX}), we have
\begin{align}\label{eq:ffstarnnormXXXX}
\|f-f^*\|_n \lesssim \lambda_{2m}^{\frac{m_0}{4m}}Q(n)^{1/2}.
\end{align}
Note that \eqref{eq:underffstar} is also valid under the case of $f\notin \mathcal{N}_{\Psi(\Omega)}$. Combining \eqref{eq:ffstarnnormXXXX} with \eqref{eq:underffstar} and \eqref{eq:unfstarboundX}, yields
\begin{align}\label{eq:underffstar1X}
\|f- f^*_n\|_{n}^2 + \lambda_{2m} \|f^*_n\|^2_{\mathcal{N}_{\Psi_{2m}}(\Omega)} \lesssim \lambda_{2m}^{\frac{m_0}{2m}}Q(n).
\end{align}
By the Gagliardo–Nirenberg interpolation inequality, we have
\begin{align}\label{eq:GNinterunderX}
\|f^*_n\|_{H^{m_0-\delta}(\Omega)} \lesssim  \|f^*_n\|_{L_2(\Omega)}^{\frac{2m-m_0+\delta}{2m}}\|f^*_n\|_{H^{2m}(\Omega)}^{\frac{m_0-\delta}{2m}} \lesssim \|f^*_n\|_{L_2(\Omega)}^{\frac{2m-m_0+\delta}{2m}}\|f^*_n\|_{\mathcal{N}_{\Psi_{2m}}(\Omega)}^{\frac{m_0-\delta}{2m}},
\end{align}
where the last inequality is because of the equivalence of $\|\cdot\|_{H^{2m}(\Omega)}$ and $\| \cdot\|_{\mathcal{N}_{\Psi_{2m}}(\Omega)}$. By Lemma \ref{lemmafixdesign}, it can be seen that 
\begin{align}\label{eq:ffstarnb2X}
\|f- f^*_n\|_{L_2(\Omega)} \lesssim & h_n^{m_0-\delta}\|f- f^*_n\|_{H^{m_0-\delta}(\Omega)} + \|f- f^*_n\|_n\nonumber\\
\lesssim & h_n^{m_0-\delta}\|f\|_{H^{m_0-\delta}(\Omega)} + h_n^{m_0}\|f^*_n\|_{H^{m_0-\delta}(\Omega)} + \|f- f^*_n\|_n\nonumber\\
\lesssim & h_n^{m_0-\delta}\|f\|_{H^{m_0-\delta}(\Omega)} + h_n^{m_0-\delta}\|f^*_n\|_{L_2(\Omega)}^{\frac{2m-m_0+\delta}{2m}}\|f^*_n\|_{\mathcal{N}_{\Psi_{2m}}(\Omega)}^{\frac{m_0-\delta}{2m}} + \lambda_{2m}^{\frac{m_0}{4m}}Q(n)^{1/2},
\end{align}
where the second inequality is because of the triangle inequality, and the last inequality is because of \eqref{eq:underffstar1X} and\eqref{eq:GNinterunderX}.

Since $f_n^*\in N_{\Psi_{2m}}(\Omega)$, $f_n^*\in H^{2m}(\Omega)$. Applying Lemma \ref{lemmafixdesign} to $f_n^*$ leads to
\begin{align}\label{eq:fstarnBunX}
\|f^*_n\|_{L_2(\Omega)} \lesssim & h_n^{2m}\| f^*_n\|_{H^{2m}(\Omega)} + \|f^*_n\|_n\nonumber\\
\lesssim & h_n^{2m}\| f^*_n\|_{N_{\Psi_{2m}}(\Omega)} + \|f\|_n + \|f-f^*_n\|_n\nonumber\\
\lesssim & h_n^{2m}\lambda_{2m}^{\frac{m_0-2m}{2m}}Q(n) + \|f\|_n + \lambda_{2m}^{\frac{m_0}{4m}}Q(n)^{1/2}\nonumber\\ 
\lesssim & 1,
\end{align}
where the second inequality is by the equivalence of $\|\cdot\|_{H^{2m}(\Omega)}$ and $\|\cdot\|_{N_{\Psi_{2m}}(\Omega)}$ and the triangle inequality, the third inequality is by \eqref{eq:underffstar1X}, and the last inequality is because of \eqref{eq:reproP} and $h_n^{2m}\lambda_{2m}^{\frac{m_0-2m}{2m}}Q(n) \lesssim n^{-\frac{2m(2m_0-d)+2m_0d}{(2m_0+d)d}}Q(n) \lesssim 1$.

Similar to the proof of Theorem \ref{thm:krroverg}, we have $\|f\|_{H^{m_0-\delta}(\Omega)}\lesssim Q(n)^{1/2}$. Plugging \eqref{eq:fstarnBunX} into \eqref{eq:ffstarnb2X}, we have
\begin{align}\label{eq:ffstarBNunder2X}
\|f- f^*_n\|_{L_2(\Omega)} \lesssim & h_n^{m_0-\delta}\|f\|_{H^{m_0-\delta}(\Omega)} + h_n^{m_0-\delta}\|f^*_n\|_{L_2(\Omega)}^{\frac{m-m_0+\delta}{2m}}\|f^*_n\|_{\mathcal{N}_{\Psi_{2m}}(\Omega)}^{\frac{m_0-\delta}{2m}} + \lambda_{2m}^{\frac{m_0}{4m}}Q(n)^{1/2}\nonumber\\
\lesssim & h_n^{m_0-\delta}\|f\|_{H^{m_0-\delta}(\Omega)} + h_n^{m_0-\delta}(\lambda_{2m}^{\frac{m_0-2m}{4m}}Q(n)^{1/2})^{\frac{m_0-\delta}{2m}} + \lambda_{2m}^{\frac{m_0}{4m}}Q(n)^{1/2}\nonumber\\
\lesssim & \lambda_{2m}^{\frac{m_0}{4m}}Q(n)^{1/2},
\end{align}
where the second inequality is by \eqref{eq:underffstar1X}, and the last inequality is because $h_n^{m_0-\delta}\lambda_{2m}^{\frac{(m_0-2m)(m_0-\delta)}{8m^2}}\lesssim \lambda_{2m}^{\frac{m_0}{4m}}$ and $Q(n) \gtrsim 1$.

The Plancherel theorem \citep{bracewell1986fourier} implies that
\begin{align*}
& \|v\|_{L_2(\Omega)}^2\lesssim \|h_f\|_{L_2(\RR^d)}^2 = \|\mathcal{F}(f_{n,e}^*)/\mathcal{F}(\Phi)\|_{L_2(\RR^d)}^2 \lesssim \|f_{n,e}^*\|_{\mathcal{N}_{\Psi_{2m}}(\RR^d)}^2 \lesssim \lambda_{2m}^{\frac{m_0-2m}{2m}}Q(n).
\end{align*}
Note that \eqref{eq:impbasicunder} also holds for $f\notin \mathcal{N}_{\Psi(\Omega)}$. Repeating the process of obtaining \eqref{eq:underr1}, \eqref{eq:undercr111}, \eqref{eq:undercr112}, \eqref{eq:underr22}, \eqref{eq:underr23}, \eqref{eq:underr24}, and \eqref{eq:underr25}, we have
\begin{align}\label{eq:underallX}
\|f-\hat f_m\|_n^2 \lesssim T_1,
\|f_n^* - \hat f_m\|_{\mathcal{N}_{\Phi}(\Omega)}^2 \lesssim \lambda_m^{-1}T_1,
\end{align}
where
\begin{align*}
    T_1 =  \max\{ & \lambda_{2m}^{\frac{m_0}{2m}}Q(n) + \lambda_{m}\lambda_{2m}^{\frac{m_0-m}{2m}}Q(n), \nonumber\\
& t^{\frac{4m}{4m-d}}n^{-{\frac{2m}{4m-d}}}\lambda_m^{-\frac{d}{4m-d}} (\lambda_{2m}^{\frac{m_0}{2m}}Q(n))^{\frac{2m-d}{4m-d}}, t^2n^{-1}\lambda_m^{-\frac{d}{2m}} \}.
\end{align*}
It remains to bound the difference between $\|f-\hat f_m\|_n^2$ and $\|f-\hat f_m\|_{L_2(\Omega)}^2$. Using the Gagliardo–Nirenberg interpolation inequality, for $\delta_1 = \min\{m_0-m,\frac{m_0^2d+2mm_0(2m_0-d)}{8mm_0+4md}\} > 0$, we find that
\begin{align}\label{eq:GNinteru2X}
\|f -  f_n^*\|_{H^{m_0-\delta_1}(\Omega)}&  \leq \|f\|_{H^{m_0-\delta_1 }(\Omega)} +  \|f_n^*\|_{H^{m_0-\delta_1 }(\Omega)}\nonumber\\
& \lesssim \|f\|_{H^{m_0-\delta_1 }(\Omega)} +  \|f_n^*\|_{L_2(\Omega)}^{\frac{2m - m_0+\delta_1 }{2m}}\|f_n^*\|_{H^{2m}(\Omega)}^{\frac{m_0-\delta_1 }{2m}}\nonumber\\
&  \lesssim  \|f\|_{H^{m_0-\delta_1 }(\Omega)} +  \|f_n^*\|_{L_2(\Omega)}^{\frac{2m - m_0+\delta_1}{2m}}\|f_n^*\|_{\mathcal{N}_{\Psi_{2m}}(\Omega)}^{\frac{m_0-\delta_1}{2m}}\nonumber\\
&  \lesssim  \|f\|_{H^{m_0-\delta_1 }(\Omega)} +  (\lambda_{2m}^{\frac{m_0-2m}{4m}}Q(n)^{1/2})^{\frac{m_0-\delta_1}{2m}}\nonumber\\
&  \lesssim  (\lambda_{2m}^{\frac{m_0-2m}{4m}}Q(n)^{1/2})^{\frac{m_0-\delta_1}{2m}},
\end{align}
where the first inequality is by the triangle inequality, the third inequality is by the equivalence of $\|\cdot\|_{H^{2m}(\Omega)}$ and $\|\cdot\|_{\mathcal{N}_{\Psi_{2m}}}$, the fourth inequality is by \eqref{eq:underffstar1X} and \eqref{eq:fstarnBunX}, and the last inequality is by $\|f\|_{H^{m_0-\delta_1 }(\Omega)}\lesssim Q(n)^{1/2}$ (which can be shown similarly as showing $\|f\|_{H^{m_0-\delta}(\Omega)}\lesssim Q(n)^{1/2}$), and $Q(n)^{\frac{2m-m_0}{4m}}\lesssim \lambda_{2m}^{\frac{m_0(m_0-2m)}{8m^2}}$.

By the Gagliardo–Nirenberg interpolation inequality, we have
\begin{align}\label{eq:GNinteru1X}
& \|f -  f_n^*\|_{\mathcal{N}_{\Phi}(\Omega)} \lesssim \|f -  f_n^*\|_{H^{m}(\Omega)}\lesssim
\|f -  f_n^*\|_{L_2(\Omega)}^{\frac{m_0-\delta_1-m}{m_0-\delta_1}}\|f -  f_n^*\|_{H^{m_0-\delta_1}(\Omega)}^{\frac{m}{m_0-\delta_1}}\nonumber\\
\lesssim &  (\lambda_{2m}^{\frac{m_0}{4m}}Q(n)^{1/2})^{\frac{m_0-\delta_1-m}{m_0-\delta_1}}((\lambda_{2m}^{\frac{m_0-2m}{4m}}Q(n)^{1/2})^{\frac{m_0-\delta_1}{2m}})^{\frac{m}{m_0-\delta_1}}\nonumber\\
= & \lambda_{2m}^{\frac{3m_0^2 -3m_0\delta_1 - 4mm_0 + 2m\delta_1}{8m(m_0-\delta_1)}}Q(n)^{\frac{3m_0-2m-3\delta_1}{4(m_0-\delta_1)}},
\end{align}
where the first inequality is because of the equivalence of $\|\cdot\|_{H^m(\Omega)}$ and $\| \cdot\|_{\mathcal{N}_{\Phi}(\Omega)}$, and third inequality is by \eqref{eq:ffstarBNunder2X} and \eqref{eq:GNinteru2X}.

By Lemma \ref{lemmafixdesign}, we have
\begin{align*}
\| f - \hat f_m\|_{L_2(\Omega)} \lesssim & h_n^m\| f - \hat f_m\|_{H^m(\Omega)} + \| f - \hat f_m\|_{n}\nonumber\\
\lesssim & h_n^m\| f - \hat f_m\|_{\mathcal{N}_{\Phi}(\Omega)} + \| f - \hat f_m\|_{n}\nonumber\\
\lesssim & h_n^m\| f_n^* - \hat f_m\|_{\mathcal{N}_{\Phi}(\Omega)} +h_n^m\|f -  f_n^*\|_{\mathcal{N}_{\Phi}(\Omega)} + \| f - \hat f_m\|_{n}\nonumber\\
\lesssim & h_n^m\lambda_m^{-1/2}T_1^{1/2} +h_n^m\lambda_{2m}^{\frac{3m_0^2 -3m_0\delta_1 - 4mm_0 + 2m\delta_1}{8m(m_0-\delta_1)}}Q(n)^{\frac{3m_0-2m-3\delta_1}{4(m_0-\delta_1)}} + T_1^{1/2}\nonumber\\
\lesssim & T_1^{1/2},
\end{align*}
where the second inequality is by the equivalence of $\|\cdot\|_{H^{m}(\Omega)}$ and $\|\cdot\|_{\mathcal{N}_{\Phi}(\Omega)}$, the third inequality is by the triangle inequality. The fourth inequality is by \eqref{eq:underall}, \eqref{eq:underallX}, \eqref{eq:GNinteru2X} and \eqref{eq:GNinteru1X}, and the last inequality is by the facts that $\lambda_m \gtrsim n^{-\frac{2m}{d}}$, $h_n \lesssim n^{-1/d}$, $\lambda_{2m} = n^{-\frac{4m}{2m_0+d}}$, and that since $T_1$ is always larger than $n^{-\frac{2m_0}{2m_0+d}}$, it can be checked that $T_1^{1/2}$ is always larger than $h_n^m\lambda_{2m}^{\frac{3m_0^2 -3m_0\delta_1 - 4mm_0 + 2m\delta_1}{8m(m_0-\delta_1)}}Q(n)^{\frac{3m_0-2m-3\delta_1}{4(m_0-\delta_1)}}$, which is also because of our choice of $\delta_1$. Therefore, we finish the proof.

\section{Proof of Proposition \ref{prop:kgrela1}}\label{subsec:pfofpropkg}
Let $u = (u_1,...,u_n)^T =  (R_m + \mu_m I_n)^{-1}r_m(x)$. By the representer theorem, 
\begin{align*}
\hat f_m(x) = u^T Y = u^T F + u^T \epsilon,
\end{align*}
where $F= (f(x_1),...,f(x_n))^T$ and $\epsilon = (\epsilon_1,...,\epsilon_n)^T$. Therefore,
\begin{align}\label{eq:pfpropkgf1}
\mathbb{E}(f(x) - \hat f_m(x))^2 = & (f(x) - u^TF)^2 + \sigma_\epsilon^2 u^Tu,
\end{align}
where $\sigma_\epsilon^2$ is the variance of $\epsilon_j$. The Fourier inverse theorem implies
\begin{align}\label{eq:pfpropkgf}
& (f(x) - \sum_{j=1}^n u_j f(x_j))^2\nonumber\\
= & \bigg|\frac{1}{(2\pi)^d} \int_{\RR^d} \bigg(\sum_{j=1}^n u_j e^{-i\langle x_j, \omega \rangle} - e^{-i\langle x, \omega \rangle}\bigg)\mathcal{F}(f)(\omega) d\omega  \bigg| ^2\nonumber\\
\leq & \frac{1}{(2\pi)^d}\int_{\RR^d} \bigg| \sum_{j=1}^n u_j e^{-i\langle x_j, \omega \rangle} - e^{-i\langle x, \omega \rangle}\bigg|^2 \mathcal{F}(\Psi)(\omega) d\omega  \frac{1}{(2\pi)^d}\int_{\RR^d} \frac{|\mathcal{F}(f)(\omega)|^2}{\mathcal{F}(\Psi)(\omega)} d\omega \nonumber\\
= & \bigg(\Psi(x-x) - 2\sum_{j=1}^n u_j\Psi(x-x_j) + \sum_{j=1}^n\sum_{j=1}^n u_k u_j\Psi(x_k-x_j)\bigg)\|f\|_{\mathcal{N}_\Psi(\Omega)}^2,
\end{align}
where the inequality is by the Cauchy-Schwarz inequality.

Now consider $\mathbb{E}(Z(x) - \hat f_G(x))^2$. Direct computation shows that
\begin{align}\label{eq:pfpropkgz}
\mathbb{E}(Z(x) - \hat f_G(x))^2 = \bigg(\Psi(x-x) - 2\sum_{j=1}^n u_j\Psi(x-x_j) + \sum_{j=1}^n\sum_{j=1}^n u_k u_j\Psi(x_k-x_j)\bigg) + \sigma_\epsilon^2u^Tu.
\end{align}
Plugging \eqref{eq:pfpropkgf} into \eqref{eq:pfpropkgf1}, we have
\begin{align*}
\mathbb{E}(f(x) - \hat f_m(x))^2 \leq & \bigg(\Psi(x-x) - 2\sum_{j=1}^n u_j\Psi(x-x_j) + \sum_{j=1}^n\sum_{j=1}^n u_k u_j\Psi(x_k-x_j)\bigg)\|f\|_{\mathcal{N}_\Psi(\Omega)}^2\nonumber\\
& + \sigma_\epsilon^2u^Tu\nonumber\\
\leq & C\bigg(\Psi(x-x) - 2\sum_{j=1}^n u_j\Psi(x-x_j) + \sum_{j=1}^n\sum_{j=1}^n u_k u_j\Psi(x_k-x_j) + \sigma_\epsilon^2u^Tu\bigg)\nonumber\\
= & C\mathbb{E}(Z(x) - \hat f_G(x))^2,
\end{align*}
where $C= \max\{1,\|f\|_{\mathcal{N}_\Psi(\Omega)}^2\}$, and the last equality is by \eqref{eq:pfpropkgz}. This finishes the proof.

\section{Proof of auxiliary lemmas}

\subsection{Proof of Lemma \ref{LEM:BOFKMMAX}}\label{app:pfbofKMmax}
We first state some lemmas used in this proof. Lemma \ref{LEM:BOFKMMAXX} can be obtained by repeating the arguments used to establish Theorem 2.2 of \cite{narcowich1994condition}; also see Equation (4) in \cite{narcowich2006sobolev}. Lemma \ref{lem:bofMatern} is Lemma 3.2 of \cite{narcowich2006sobolev}. Lemma \ref{LEM:BOFKMGX} can be obtained by elementary mathematical analysis, and the proof is in Appendix \ref{app:pflembofkMgx}. 

\begin{lemma}\label{LEM:BOFKMMAXX}
		Let $\Psi$ be a Mat\'ern correlation function satisfying Condition (C2). Suppose the design points $X=\{x_1,...,x_n\}$. Let $\Lambda_X$ be the maximum eigenvalue of matrix $(\Psi(x_j-x_k))_{jk}$. Then 
	\begin{align*}
	\Lambda_X\leq \Psi(0)+\sum_{k=1}^\infty 3d(k+2)^{d-1}\Psi(kq_X),
	\end{align*}
	where $q_X$ is the separation radius of $X$.
\end{lemma}

\begin{lemma}\label{lem:bofMatern}
	Let $\Psi$ be a Mat\'ern correlation function satisfying Condition (C2). Then we have that $\Psi$ is positive definite, decreasing on $[0,\infty)$, and satisfies the bound 
	\begin{align*}
	\Psi(x)\leq \sqrt{2\pi c_{m_0}}r^{\nu - 1/2}e^{-r+\nu^2/2r},\quad r=\|x\|_2>0,
	\end{align*}
	where $\nu=m_0-d/2$.
\end{lemma}

\begin{lemma}\label{LEM:BOFKMGX}
	Define function $g$ as
	\begin{align*}
	g(x) = x^{d/2 + m_0 - 3/2} e^{-x},\quad x\geq 1.
	\end{align*}
	We have $g(x) \leq Ce^{-x/2}$ for all $x\geq 1$, where $C = (d+2m_0)^{d/2+m_0}$.
\end{lemma}

Now we begin to prove Lemma \ref{LEM:BOFKMMAX}. Note that $k+2\leq 3k$ for $k\geq 1$, which, together with Lemma \ref{LEM:BOFKMMAXX}, leads to
\begin{align}\label{eq:pflemeigenM}
\Lambda_X\leq & \Psi(0)+\sum_{k=1}^\infty 3d(k+2)^{d-1}\Psi(kq_X) \nonumber\\
\leq & \Psi(0)+3^dd\sum_{k=1}^\infty k^{d-1}\Psi(kq_X)\nonumber\\
= & \Psi(0) + 3^dd\sum_{k=1}^{\lfloor 1/q_X\rfloor} k^{d-1}\Psi(kq_X) + 3^dd\sum_{k=\lfloor 1/q_X\rfloor + 1}^\infty k^{d-1}\Psi(kq_X)\nonumber\\
= & \Psi(0) + I_1 + I_2,
\end{align}
where $\lfloor a \rfloor$ is the integer part of $a\in \RR_+$. Since $\Psi(x)$ is a decreasing function, the first term $I_1$ can be bounded by 
\begin{align}\label{eq:pflemeigenMI1}
I_1 \leq 3^dd\sum_{k=1}^{\lfloor 1/q_X\rfloor} k^{d-1}\Psi(0) \leq C_13^dd\Psi(0)(\lfloor 1/q_X\rfloor)^d\leq C_13^dd\Psi(0)(1/q_X)^d,
\end{align}
where we utilize the basic inequality $\sum_{k=1}^m k^{d-1}\leq C_1m^d$.

Using Lemma \ref{lem:bofMatern}, the second term $I_2$ can be bounded by 
\begin{align*}
I_2 \leq 3^dd\sum_{k=\lfloor 1/q_X\rfloor+1}^\infty k^{d-1}\Psi(kq_X)\leq C_2\sum_{k=\lfloor 1/q_X\rfloor+1}^\infty k^{d-1}(kq_X)^{m_0- d/2 - 1/2}e^{-kq_X+(m_0- d/2)^2/(2kq_X)}.
\end{align*}
Clearly, $kq_X \geq 1$ for $k\geq \lfloor 1/q_X\rfloor+1$, which implies $e^{(m_0- d/2)^2/(2kq_X)}\leq e^{(m_0- d/2)^2/2}$. Therefore, $I_2$ can be further bounded by
\begin{align}\label{eq:pflemeigenMI22}
I_2 \leq & C_2e^{(m_0- d/2)^2/2}q_X^{-d+1}\sum_{k=\lfloor 1/q_X\rfloor+1}^\infty (kq_X)^{d-1}(kq_X)^{m_0- d/2 - 1/2}e^{-kq_X}\nonumber\\
= & C_2e^{(m_0- d/2)^2/2}q_X^{-d+1}\sum_{k=\lfloor 1/q_X\rfloor+1}^\infty (kq_X)^{m_0 + d/2 - 3/2}e^{-kq_X}\nonumber\\
= & C_2e^{(m_0- d/2)^2/2}q_X^{-d+1}I_3,
\end{align}
where
\begin{align*}
I_3 = \sum_{k=\lfloor 1/q_X\rfloor+1}^\infty (kq_X)^{d/2 + m_0 - 3/2}e^{-kq_X}.
\end{align*}
Consider function $
g(x) = x^{d/2 + m_0 - 3/2} e^{-x}, x\geq 1.$
By Lemma \ref{LEM:BOFKMGX}, we have $g(x) \leq C_3e^{-x/2}$, where $C_3 = (d+2m_0)^{d/2+m_0}$. This implies that $I_3$ can be bounded by 
\begin{align}\label{eq:pflemeigenMI3}
I_3 & \leq C_3\sum_{k=\lfloor 1/q_X\rfloor+1}^\infty e^{-kq_X/2}\leq C_3\sum_{k=0}^\infty e^{-kq_X/2} \nonumber\\
& = \frac{C_3}{1-e^{-q_X/2}} \leq C_3\frac{q_X/2 + 1}{q_X/2} \lesssim q_X^{-1},
\end{align}
where third inequality is because of the basic inequality $1-e^{-x}>\frac{x}{x+1}$ for $x>0$, and the last equality is because $q_X\lesssim 1$. Plugging \eqref{eq:pflemeigenMI3} into \eqref{eq:pflemeigenMI22}, we have $I_2\leq C_4 q_X^{-d}$. Together with \eqref{eq:pflemeigenMI1} and \eqref{eq:pflemeigenM}, we can see the desired result holds.

\subsection{Proof of Lemma \ref{LEM:OPTRATE}}\label{app:pfcoropt}
Let $u=(u_1,...,u_n)^T = (R + \mu I_n)^{-1}r(x)$, where $R$, $r$ and $\mu$ are as in \eqref{mean}. 
The Fourier inversion theorem implies that 
\begin{align*}
    & K(x-x) - r_K(x)^T(R_K + \mu I_n)^{-1}r_K(x)\nonumber\\
    \leq & K(x-x) - 2\sum_{j=1}^nu_jK(x-x_j) + \sum_{j=1}^n\sum_{k=1}^n u_ju_kK(x_j-x_k) + \mu\|u\|_2^2\nonumber\\
    = & \frac{1}{(2\pi)^d}\int_{\RR^d} \bigg| \sum_{j=1}^n u_j e^{-i\langle x_j, \omega \rangle} - e^{-i\langle x, \omega \rangle}\bigg|^2 \mathcal{F}(K)(\omega) d\omega + \mu\|u\|_2^2 \nonumber\\
    \leq & A_0\frac{1}{(2\pi)^d}\int_{\RR^d} \bigg| \sum_{j=1}^n u_j e^{-i\langle x_j, \omega \rangle} - e^{-i\langle x, \omega \rangle}\bigg|^2 \mathcal{F}(\Psi)(\omega) d\omega + \mu\|u\|_2^2\nonumber\\
    \leq & \max\{A_0,1\}\left(\frac{1}{(2\pi)^d}\int_{\RR^d} \bigg| \sum_{j=1}^n u_j e^{-i\langle x_j, \omega \rangle} - e^{-i\langle x, \omega \rangle}\bigg|^2 \mathcal{F}(\Psi)(\omega) d\omega + \mu\|u\|_2^2\right)\nonumber\\
    = & \max\{A_0,1\}(\Psi(x-x) - r(x)^T(R + \mu I_n)^{-1}r(x)),
\end{align*}
where $r$ and $R$ are as in \eqref{mean}. This finishes the proof.

\subsection{Proof of Lemma \ref{LEM:FSTAR}}\label{app:pflemmafstar}
In this section, we set $m_0:=m_0(f)$ for notational simplicity. Since $\Omega$ has Lipschitz boundary, there exists an extension operator from $L_2(\Omega)$ to $L_2(\mathbb{R}^d)$, such that the smoothness of each function is maintained \citep{devore1993besov,rychkov1999restrictions}. Therefore, there exist constants $0 < C_1 \leq C_2$ such that for any functions $h_1 \in H^m(\Omega)$ and $h_2 \in H^{m_0}(\Omega)$, there exist $h_{1,e} \in H^m(\RR^d)$ and $h_{2,e} \in H^{m_0}(\RR^d)$ satisfying
\begin{align}
C_1\|h_{1,e}\|_{H^m(\RR^d)} \leq \|h_1\|_{H^m(\Omega)} \leq C_2\|h_{1,e}\|_{H^m(\RR^d)},\label{eq:Hmext}\\
C_1\|h_{2,e}\|_{H^{m_0}(\RR^d)} \leq \|h_2\|_{H^{m_0}(\Omega)} \leq C_2\|h_{2,e}\|_{H^{m_0}(\RR^d)},\label{eq:Hm0ext}
\end{align}
and $h_{1,e}(x) = h_1(x)$ and $h_{2,e}(x) = h_2(x)$ for any $x\in \Omega$. Let $f_1$ be the solution to the optimization problem
\begin{align}\label{eq:overoptR}
\min_{\tilde f \in H^m(\RR^d)}\|f_e-\tilde f\|_{L_2(\RR^d)}^2 + \lambda_m \|\tilde f\|^2_{H^{m}(\RR^d)},
\end{align}
where $f_e$ is the extension of $f$ satisfying \eqref{eq:Hm0ext} with $h_2 = f$.

Because $\Phi$ satisfies Condition (C3), $\mathcal{N}_{\Phi}(\Omega)$ coincides with $H^m(\Omega)$. Since $f^*$ is the solution to \eqref{eq:overoptOmega}, by \eqref{eq:Hmext}, we have
\begin{align}\label{eq:lem1pf1}
& \|f-f^*\|_{L_2(\Omega)}^2 + \lambda_m \|f^*\|^2_{\mathcal{N}_{\Phi}(\Omega)} \leq \|f-f_{1,r}\|_{L_2(\Omega)}^2 + \lambda_m \|f_{1,r}\|^2_{\mathcal{N}_{\Phi}(\Omega)}\nonumber\\
\leq & C_3\left(\|f-f_{1,r}\|_{L_2(\Omega)}^2 + \lambda_m \|f_{1,r}\|^2_{H^m(\Omega)}\right)\leq C_4\left(\|f_e-f_1\|_{L_2(\RR^d)}^2 + \lambda_m \|f_1\|^2_{H^m(\RR^d)}\right),
\end{align}
where $f_{1,r}$ is the restriction of $f_1$ onto $\Omega$.

By the Fourier transform and the Plancherel theorem \citep{bracewell1986fourier}, we have 
\begin{align}\label{eq:lem1pff1}
&\|f_e-f_1\|_{L_2(\RR^d)}^2+\lambda_m\|f_1\|^2_{H^m(\RR^d)} \nonumber\\
= & \int_{\RR^d} |\mathcal{F}(f_e)(\omega)-\mathcal{F}(f_1)(\omega)|^2 d\omega + \lambda_m\int_{\RR^d} |\mathcal{F}(f_1)(\omega)|^2(1+\|\omega\|_2^2)^{m}d\omega\nonumber\\
= & \int_{\RR^d}\big( |\mathcal{F}(f_e)(\omega)-\mathcal{F}(f_1)(\omega)|^2 + \lambda_m |\mathcal{F}(f_1)(\omega)|^2(1+\|\omega\|_2^2)^{m}\big)d\omega\nonumber\\
= & \int_{\RR^d} \frac{\lambda_m(1+\|\omega\|_2^2)^{m}}{1+\lambda_m(1+\|\omega\|_2^2)^{m}}|\mathcal{F}(f_e)(\omega)|^2d\omega\nonumber\\
= & \int_{\Omega_1} \frac{\lambda_m(1+\|\omega\|_2^2)^{m}}{1+\lambda_m(1+\|\omega\|_2^2)^{m}}|\mathcal{F}(f_e)(\omega)|^2d\omega + \int_{\Omega_1^C} \frac{\lambda_m(1+\|\omega\|_2^2)^{m}}{1+\lambda_m(1+\|\omega\|_2^2)^{m}}|\mathcal{F}(f_e)(\omega)|^2d\omega\nonumber\\
\leq & \int_{\Omega_1} \lambda_m(1+\|\omega\|_2^2)^{m}|\mathcal{F}(f_e)(\omega)|^2d\omega + \int_{\Omega_1^C} |\mathcal{F}(f_e)(\omega)|^2d\omega,
\end{align}
where $\Omega_1=\{\omega: \lambda_m(1+\|\omega\|_2^2)^{m}<1\}$, $\Omega_1^C = \RR^d\setminus \Omega$, and the third equality follows that $f_1$ is the solution to \eqref{eq:overoptR}. 

On the set $\Omega_1$, since $\lambda_m(1+\|\omega\|_2^2)^{m}<1$ and $m\geq m_0$, it can be verified that $\lambda_m(1+\|\omega\|_2^2)^{m}\leq \lambda_m^{\frac{m_0}{m}}(1+\|\omega\|_2^2)^{m_0}$. On the other hand, since $\lambda_m(1+\|\omega\|_2^2)^{m}\geq 1$ on the set $\Omega_1^C$, we have $ \lambda_m^{\frac{m_0}{m}}(1+\|\omega\|_2^2)^{m_0}\geq 1$. Together with \eqref{eq:lem1pff1}, we obtain 
\begin{align}\label{eq:lem1pff12}
&\|f_e-f_1\|_{L_2(\RR^d)}^2+\lambda_m\|f_1\|^2_{H^m(\RR^d)} \nonumber\\
\leq & \int_{\Omega_1} \lambda_m(1+\|\omega\|_2^2)^{m}|\mathcal{F}(f_e)(\omega)|^2d\omega + \int_{\Omega_1^C} |\mathcal{F}(f_e)(\omega)|^2d\omega\nonumber\\
\leq & \lambda_m^{\frac{m_0}{m}}\int_{\Omega_1} (1+\|\omega\|_2^2)^{m_0}|\mathcal{F}(f_e)(\omega)|^2d\omega + \lambda_m^{\frac{m_0}{m}}\int_{\Omega_1^C} (1+\|\omega\|_2^2)^{m_0}|\mathcal{F}(f_e)(\omega)|^2d\omega\nonumber\\
= & \lambda_m^{\frac{m_0}{m}}\|f_e\|_{H^{m_0}(\RR^d)}^2 \leq C_5\lambda_m^{\frac{m_0}{m}}\|f\|_{H^{m_0}(\Omega)}^2\leq C_6\lambda_m^{\frac{m_0}{m}}\|f\|_{\mathcal{N}_{\Psi}(\Omega)}^2,
\end{align}
where the third inequality follows from \eqref{eq:Hm0ext} and the last inequality is because of the equivalence of $\|\cdot\|_{H^{m_0}(\Omega)}$ and $\|\cdot\|_{\mathcal{N}_{\Psi}(\Omega)}$.

Combining \eqref{eq:lem1pf1} and \eqref{eq:lem1pff12} yields 
\begin{align*}
\|f-f^*\|_{L_2(\Omega)}^2 + \lambda_m \|f^*\|^2_{\mathcal{N}_{\Phi}(\Omega)} \leq C_6\lambda_m^{\frac{m_0}{m}}\|f\|_{\mathcal{N}_{\Psi}(\Omega)}^2,
\end{align*}
which implies
\begin{align*}
\|f^*-f\|_{L_2(\Omega)}^2 \leq C_6\lambda_m^{\frac{m_0}{m}}\|f\|_{\mathcal{N}_{\Psi}(\Omega)}^2
\mbox{ and }
\|f^*\|^2_{\mathcal{N}_{\Phi}(\Omega)} \leq C_6\lambda_m^{\frac{m_0-m}{m}}\|f\|_{\mathcal{N}_{\Psi}(\Omega)}^2.
\end{align*}
Next, we consider the case $f\notin H^{m_0}(\Omega)$. If $f\notin H^{m_0}(\Omega)$, by Lemma \ref{LEMMA1}, there exists a function $Q:\RR_+\mapsto\RR_+$ such that 
\begin{align*}
\int_{\RR^d} \frac{|\mathcal{F}(f)(\omega)|^2(1+\|\omega\|_2^2)^{m_0}}{Q(\|\omega\|_2)}d\omega \leq 1, \mbox{ and }\lim_{r\rightarrow +\infty} \frac{\log Q(r)}{\log r} = 0.
\end{align*}

Therefore, \eqref{eq:lem1pff1} can be changed to 
\begin{align}\label{eq:lem1pff2}
& \|f_e-f_1\|_{L_2(\RR^d)}^2+\lambda_m\|f_1\|^2_{H^m(\RR^d)} \nonumber\\
= &  \int_{\RR^d} \frac{\lambda_m(1+\|\omega\|_2^2)^{m}}{1+\lambda_m(1+\|\omega\|_2^2)^{m}}|\mathcal{F}(f)(\omega)|^2d\omega\nonumber\\
= & \int_{\Omega_2} \frac{\lambda_m(1+\|\omega\|_2^2)^{m}}{1+\lambda_m(1+\|\omega\|_2^2)^{m}}|\mathcal{F}(f)(\omega)|^2d\omega + \int_{\Omega_2^C} \frac{\lambda_m(1+\|\omega\|_2^2)^{m}}{1+\lambda_m(1+\|\omega\|_2^2)^{m}}|\mathcal{F}(f)(\omega)|^2d\omega\nonumber\\
\leq & \int_{\Omega_2} \lambda_m(1+\|\omega\|_2^2)^{m}|\mathcal{F}(f)(\omega)|^2d\omega + \int_{\Omega_2^C} |\mathcal{F}(f)(\omega)|^2d\omega\nonumber\\
\leq & C_7\left(\lambda_m^{\frac{m_0}{m}}Q_1(n)\int_{\Omega_2} (1+\|\omega\|_2^2)^{m_0}\frac{|\mathcal{F}(f)(\omega)|^2}{Q(\|\omega\|_2)}d\omega + \lambda_m^{\frac{m_0}{m}}\int_{\Omega_2^C} (1+\|\omega\|_2^2)^{m_0}\frac{|\mathcal{F}(f)(\omega)|^2}{Q(\|\omega\|_2)}d\omega\right)\nonumber\\
\leq & C_{8}\lambda_m^{\frac{m_0}{m}}Q_1(n),
\end{align}
where $\Omega_2 = \{\omega: \lambda_m(1+\|\omega\|_2^2)^{m} <Q(\|\omega\|_2)^{m/m_0}\}$, $\Omega_2^C = \RR^d\setminus \Omega_2$, and $Q_1:\RR_+ \mapsto \RR_+$ is a function such that $\sup_{\omega\in \Omega_2}Q(\omega)^{m/m_0} = Q_1(n)$, since we assume $\lambda_m\asymp n^\alpha$ for some $\alpha$. It can be seen that $Q_1$ satisfies 
\begin{align*}
    \lim_{r\rightarrow +\infty} \frac{\log Q_1(r)}{\log r} = 0.
\end{align*}
The second inequality of \eqref{eq:lem1pff2} follows from the fact that if $\omega\in \Omega_2$,
\begin{align*}
\frac{\lambda_m(1+\|\omega\|_2^2)^{m}}{Q(\|\omega\|_2)^{m/m_0}}<\bigg(\frac{\lambda_m(1+\|\omega\|_2^2)^{m}}{Q(\|\omega\|_2)^{m/m_0}}\bigg)^{m_0/m},
\end{align*}
and otherwise $$\bigg(\frac{\lambda_m(1+\|\omega\|_2^2)^{m}}{Q(\|\omega\|_2)^{m/m_0}}\bigg)^{m_0/m}  > C_9 >0.$$
Note that \eqref{eq:lem1pf1} also holds for $f\notin\mathcal{N}_{\Psi}(\Omega)$. Therefore, it can be seen that \eqref{eq:lem1pf1} and \eqref{eq:lem1pff2} imply that 
\begin{align*}
\|f^*-f\|_{L_2(\Omega)}^2 \leq C_8\lambda_m^{\frac{m_0}{m}}Q_1(n)
\mbox{ and }  \|f^*\|^2_{\mathcal{N}_{\Phi}(\Omega)} \leq C_8\lambda_m^{\frac{m_0-m}{m}}Q_1(n).
\end{align*}
This finishes the proof of Lemma \ref{LEM:FSTAR}.

\subsection{Proof of Lemma \ref{LEM:BOFKMGX}}\label{app:pflembofkMgx}
The result is implied by
\begin{align*}
& x^{d/2 + m_0 - 3/2} e^{-x} \leq Ce^{-x/2}
\Leftrightarrow  x^{d/2 + m_0 - 3/2} \leq Ce^{x/2} 
\Leftrightarrow (d/2 + m_0 - 3/2)\log x\leq \log C + \frac{x}{2}.
\end{align*}
Consider function $h(x) = \log C + \frac{x}{2} - (\frac{d}{2} + m_0 - \frac{3}{2})\log x$. The first order derivative of $h(x)$ is
\begin{align*}
\frac{dh(x)}{dx} = \frac{1}{2} - \frac{d/2 + m_0 - 3/2}{x}.
\end{align*}
If $d/2+m_0-3/2 \leq 1/2$, then $h(x)$ is an increasing function on $[1,\infty)$, which implies $h(x) \geq h(1) > 0$, and the result of Lemma \ref{LEM:BOFKMGX} holds. If $d/2+m_0-3/2 > 1/2$, then $h$ decreases on $[1,d+2m_0-3)$, and increases on $[d+2m_0-3,\infty)$, which implies that $h$ takes the minimum at $x=d+2m_0-3$. Since $h(x) \geq h(d+2m_0-3)$ for all $x\geq 1$ and $h(d+2m_0-3) > 0$, we finish the proof of Lemma \ref{LEM:BOFKMGX}.

\section{Proof of statements in Example \ref{eg:tf1}}\label{app:pfeg1}
Direct computation shows that the Fourier transform of $f$ is
\begin{align*}
    \mathcal{F}(f)(\omega) = \frac{4\sin^2(\omega/2)}{\sqrt{2\pi}\omega^2},\quad \omega\in \RR.
\end{align*}
It can be seen that
\begin{align*}
    & \int_\RR |\mathcal{F}(f)(\omega)|^2(1 + |\omega|^2)^{3/2}d\omega\\
    \geq & \int_\RR |\mathcal{F}(f)(\omega)|^2(1 + |\omega|^3) d\omega = \int_\RR \frac{2(1-\cos(\omega))^2}{\pi\omega^4}(1 + |\omega|^3) d\omega \\
    \geq & \int_{\pi/2}^\infty  \frac{2(1-\cos(\omega))^2}{\pi\omega}d\omega = \int_0^\infty \frac{2(1-\sin(t))^2}{\pi(t+\pi/2)}dt\\
    \geq & \int_0^\infty \frac{2}{\pi(t+\pi/2)}dt - \int_0^\infty \frac{4\sin(t)}{\pi(t+\pi/2)}dt\\
    \geq & \int_0^\infty \frac{2}{\pi(t+\pi/2)}dt - \int_0^\infty \frac{4\sin(t)}{\pi t}dt.
\end{align*}
Note that $\int_0^\infty \frac{2}{\pi(t+\pi/2)}dt = \infty$ and $\int_0^\infty \frac{4\sin(t)}{\pi t}dt = 2$ \citep{bartle2000introduction}, which implies $ \int_\RR |\mathcal{F}(f)(\omega)|^2(1 + |\omega|)^{3/2} d\omega = \infty$. This implies $f$ does not belong to the Sobolev space $H^{3/2}(\RR)$. By checking that
\begin{align*}
& \int_\RR |\mathcal{F}(f)(\omega)|^2(1 + |\omega|^2)^{3/2-\delta/2} d\omega \\
\leq   & \int_\RR |\mathcal{F}(f)(\omega)|^2(1 + |\omega|)^{3-\delta} d\omega = \int_\RR \frac{2(1-\cos(\omega))^2}{\pi\omega^4}(1 + |\omega|)^{3-\delta} d\omega\\
\leq &\int_{-1}^1\frac{2(1-\cos(\omega))^2}{\pi\omega^4}(1 + |\omega|)^{3-\delta} d\omega + \int_{\RR\setminus [-1,1]} \frac{8}{\pi\omega^4}(1 + |\omega|)^{3-\delta} d\omega< \infty,\forall \delta>0,
\end{align*}
we can conclude that $f$ has smoothness $3/2$. It is easily seen that, the function $Q(t):=C\log^2 (1+t)$ defined on $\RR_+$ with $C$ an appropriate constant satisfies $\int_{\RR}\frac{|\mathcal{F}(f)(\omega)|^2}{Q(|\omega|)}(1 + |\omega|^2)^{3/2}d\omega \leq 1$.

\end{document}